\documentclass[11pt]{article}
\usepackage[margin=1.15in]{geometry}
\usepackage{amsmath,amssymb,amsthm}
\usepackage[colorlinks=true,linkcolor=blue,citecolor=blue,urlcolor=blue]{hyperref}
\hypersetup{pdftitle={Stein's method for negative binomial and hyper-negative-binomial laws, with applications},
  pdfauthor={Stephen DeSalvo}}
\DeclareTextSymbolDefault{\textsection}{OMS}
\DeclareTextSymbolDefault{\textasteriskcentered}{OMS}

\newtheorem{theorem}{Theorem}[section]
\newtheorem{mainthm}{Theorem}

\newtheorem{lemma}[theorem]{Lemma}
\newtheorem{proposition}[theorem]{Proposition}
\newtheorem{corollary}[theorem]{Corollary}
\theoremstyle{definition}
\newtheorem{definition}[theorem]{Definition}

\newtheorem{openproblem}{Open Problem}
\theoremstyle{remark}
\newtheorem{remark}[theorem]{Remark}

\newcommand{\dtv}{d_{\mathrm{TV}}}
\newcommand{\dwass}{d_{\mathrm{W}}}
\newcommand{\E}{\mathbb{E}}
\newcommand{\PP}{\mathbb{P}}
\newcommand{\Zint}{\mathbb{Z}}
\newcommand{\N}{\mathbb{N}}
\newcommand{\bk}{\boldsymbol{\kappa}}
\newcommand{\bq}{\boldsymbol{q}}
\newcommand{\bkk}{\boldsymbol{k}}
\newcommand{\br}{\boldsymbol{r}}
\newcommand{\bc}{\boldsymbol{c}}
\newcommand{\HNB}{\mathrm{HNB}}
\newcommand{\HB}{\mathrm{HB}}
\newcommand{\Zhat}{\widehat{Z}}

\newcommand{\one}{\mathbf{1}}
\newcommand{\ind}[1]{\mathbf{1}\{#1\}}
\DeclareMathOperator{\Geom}{Geom}
\DeclareMathOperator{\NB}{NB}
\DeclareMathOperator{\Poi}{Poi}
\DeclareMathOperator{\Var}{Var}
\DeclareMathOperator{\Cov}{Cov}
\newcommand{\T}{\mathcal{T}}

\title{Stein's method for negative binomial and hyper-negative-binomial laws,\\ with applications}
\author{Stephen DeSalvo\thanks{Independent researcher. Email:
\texttt{stephendesalvo@gmail.com}.}}
\date{October 2026}

\begin{document}

\maketitle
\begin{abstract}
We develop Stein's method for the negative binomial law and for the
hyper-negative-binomial law, the law of the smallest row total of an
array of independent geometric variables given the differences of the
row totals, with every constant explicit.  The main instrument is a
comparison of two biasings of the approximated variable, which cancels
the rate term exactly at the matched parameter.  We provide total
variation bounds for sums of independent geometric variables at
drifting parameters, for sums of dependent geometric variables (the
Arratia--Goldstein--Gordon theorem with geometric summands and the
Poisson law as an endpoint), and for Poisson mixtures; the exact
conditional law of negative-binomial counts given their differences,
with the Bessel law as its Poisson limit; and a dependency expansion
for the column direction of a contingency table.  A final section
verifies every bound numerically, with one script per result.
\end{abstract}
\begin{quotation}\small
\noindent\emph{2020 Mathematics Subject Classification.}
Primary 60F05; Secondary 60E05, 60E15, 60G50, 62E17, 62H17.
\par\smallskip
\noindent\emph{Keywords.}
Stein's method; negative binomial law; hyper-negative-binomial law;
geometric sums; total variation distance; Poisson approximation;
Bessel law; contingency tables.
\end{quotation}

\section{Introduction}\label{sec:intro}
Throughout, $\Geom(1-y)$ is the law $\PP(X=n)=(1-y)y^{n}$ on $\N$,
$\NB(\rho,w)$ the law $\PP(n)\propto\binom{n+\rho-1}{n}w^{n}$, of mean
$\rho w/(1-w)$, $\dtv$ total variation and $\dwass$ the Wasserstein
distance on $\N$; for a law on $\N$ with pmf $p$,
$F(m)=\sum_{n\le m}p_n$ and $G(m)=\sum_{n\ge m}p_n$.  If
$X_1,\dots,X_c$ are independent and identically
distributed, $X_j\sim\Geom(1-w)$, then
\begin{equation}\label{eq:nbrow}
\sum_{j=1}^{c}X_j\sim\NB(c,w);
\end{equation}
equivalently, $\NB(\rho,w)$ is the Poisson law whose rate is Gamma
with shape $\rho$ and scale $w/(1-w)$.  If $X_{ij}$, $i\le r$,
$j\le c$, are independent, $X_{ij}\sim\Geom(1-q_i)$, so that every row
is identically distributed within itself and the row totals
$U_i=\sum_jX_{ij}\sim\NB(c,q_i)$ are independent, then, for
$k_1,\dots,k_{r-1}$ with $\kappa_i=\sum_{l\ge i}k_l\ge0$, $\kappa_r=0$,
and $w=\prod_iq_i$,
\begin{equation}\label{eq:hnbrow}
\mathcal L\bigl(U_r\bigm|U_i-U_{i+1}=k_i,\ i<r\bigr)=\HNB(\bk,c,w),
\qquad
\PP(n)\propto\prod_{i=1}^{r}\binom{n+\kappa_i+c-1}{n+\kappa_i}\,w^{n},
\end{equation}
the hyper-negative-binomial law (Theorem~\ref{thm:target}), with
$p_{n+1}/p_n=w\prod_{i}(n+\kappa_i+c)/(n+\kappa_i+1)$ and
$\HNB(0,c,w)=\NB(c,w)$ at $r=1$.  Both identities require the entries
of a row to be identically distributed, and the mixture representation
requires the rate to be Gamma.

\subsection*{Geometric sums at drifting parameters}
\begin{mainthm}[Heterogeneous geometric sums]\label{thm:A}
Let $X_1,\dots,X_c$ be independent, $X_j\sim\Geom(1-y_j)$ with
$y_j\in(0,1)$, let $W=\sum_{j=1}^{c}X_j$, and let $w\in(0,1)$.  Then
\[
\dtv\bigl(W,\NB(c,w)\bigr)\ \le\ \Gamma_{c,w}\sum_{j=1}^{c}\frac{|y_j-w|}{1-y_j},
\qquad
\Gamma_{c,w}=\sup_{m\ge0}\frac{F(m)\,G(m+1)}{(m+1)\,p_{m+1}},
\]
$p$, $F$, $G$ those of $\NB(c,w)$, and
\[
\Gamma_{c,w}\le\min\Bigl\{\frac1{1-w},\ \frac{\sqrt\pi}{\sqrt{(c+1)w(1-w)}}\Bigr\},
\qquad \Gamma_{1,w}=1,
\]
and, for $c\ge2$,
\[
\begin{gathered}
\Gamma_{c,w}\le B_{c,w}=\max\Bigl\{\frac1{wc},\ \max_{m\ge1:\,x'_m>0}\frac1{x'_m},\
\sup_{m\ge0:\,y_{m+1}>0}\frac{m+2}{(m+1)\,y_{m+1}}\Bigr\},\\
x'_m=w(c-1)-m(1-w),\qquad y_{m+1}=(m+2)-w(m+1+c).
\end{gathered}
\]
\end{mainthm}
This is Corollary~\ref{cor:nbsum} with Theorem~\ref{thm:nbfactors};
$B_{c,w}$ is Theorem~\ref{thm:factors}(ii), $\Gamma_{1,w}=1$ is
Remark~\ref{rem:placement}, and the other two bounds are
\cite{BGX15,CD19}.
The number of failures before the $c$-th success in Bernoulli trials
whose success probability drifts from stage to stage is such a $W$,
and the theorem is the negative-binomial analogue of the Poisson
approximation to the Poisson-binomial, with the spread of the stage
parameters in the role of $\sum_jp_j^{2}$.  A single species of
indistinguishable balls filling $c$ boxes independently and
geometrically, box $j$ at fugacity $y_j=e^{-\beta\epsilon_j}$, has
total $W$, and the theorem certifies the flat-spectrum law
$\NB(c,e^{-\beta\bar\epsilon})$ with the spread of the spectrum as
the error (Corollary~\ref{cor:occupancy}).

\subsection*{The coupling form}
\begin{mainthm}[Coupling form]\label{thm:coupling}
Let $X_1,\dots,X_c$ be $\N$-valued, not necessarily independent, with
$0<\E X_j<\infty$; let $W=\sum_{j=1}^{c}X_j$,
$y_j=\E X_j/(1+\E X_j)$ and $\lambda=\sum_{j=1}^{c}y_j$.  For each $j$
let $W^{(j)}$ and $W^{[j]}$ be random variables on one probability
space with
\[
\PP(W^{(j)}=n)=\frac{\E[X_j\ind{W=n}]}{\E X_j},\qquad
\PP(W^{[j]}=n)=\frac{\E[(X_j+1)\ind{W=n}]}{\E X_j+1},\qquad n\in\N .
\]
For $w\in(0,1)$ let $\pi_w=\NB(\lambda/w,\,w)$, with $p$, $F$, $G$
those of $\pi_w$,
$\Gamma_w=\sup_{n\ge0}F(n)G(n+1)/((n+1)p_{n+1})$ and
$\Lambda_w=\min\{2\Gamma_w,1/\lambda\}$; let $\pi_0=\Poi(\lambda)$,
$\Gamma_0=\min\{1,1.4\lambda^{-1/2}\}$ and
$\Lambda_0=(1-e^{-\lambda})/\lambda$.  Then, for every $w\in[0,1)$,
\[
\dtv(W,\pi_w)\ \le\ \Lambda_w\sum_{j=1}^{c}\E X_j\;
\E\bigl|W^{[j]}+1-W^{(j)}\bigr|\;+\;\Gamma_w\sum_{j=1}^{c}|w-y_j|\,\E X_j .
\]
\end{mainthm}
This is Theorem~\ref{thm:couplingbody}.  $W^{(j)}$ is $W$ biased by
$X_j$ and $W^{[j]}$ is $W$ biased by $X_j+1$; the size-biased $W^{d}$
is the mixture of the $W^{(j)}$ with weights $\E X_j/\E W$ and the
$b$-biased $W^{b}$ of \S\ref{sec:nb} the mixture of the $W^{[j]}$.
For indicators $W^{(\alpha)}$ is $1+V_\alpha$ with
$\mathcal L(V_\alpha)=\mathcal L(W-X_\alpha\mid X_\alpha=1)$, the
variable coupled to $W$ in the coupling form of the Chen--Stein
method \cite{AGG89,BHJ92}, and the first sum is its
$\sum_\alpha p_\alpha\,\E|W-V_\alpha|$ with the pair
$(W^{[j]}+1,W^{(j)})$ in the role of $(W+1,1+V_\alpha)$.  When $X_j$
is geometric and independent of the other summands, $W^{(j)}=W+1+X_j'$
and $W^{[j]}=W+X_j'$ for one independent copy $X_j'$, the $j$-th term
of the first sum is zero, and at $w=\bar y$ the theorem is
Theorem~\ref{thm:A} with the weights $y_j/(1-y_j)$.

\subsection*{Dependent geometric sums}
\begin{mainthm}[Dependent geometric sums]\label{thm:dep}
Let $X_1,\dots,X_c$ be $\N$-valued, not necessarily independent, with
$X_j\sim\Geom(1-y_j)$, $y_j\in(0,1)$; let $W=\sum_{j=1}^{c}X_j$ and
$\lambda=\sum_{j=1}^{c}y_j$.  For each $j$ fix $B_j\subseteq[c]$ with
$j\in B_j$, and let $\mathcal F_j=\sigma(X_i:i\notin B_j)$,
$R_j=\sum_{i\in B_j\setminus\{j\}}X_i$, $W_j=W-X_j-R_j$.  Set
\[
M_w=\sum_{j=1}^{c}|w-y_j|\,\frac{y_j}{1-y_j},\qquad
Q=\sum_{j=1}^{c}\sum_{i\in B_j\setminus\{j\}}\Bigl((1+y_j)\,\E[X_iX_j]+y_j\,\E X_i\Bigr),
\]
\[
\beta=\sum_{j=1}^{c}\E\sum_{n\ge1}n\,\bigl|\PP(X_j=n\mid\mathcal F_j)-y_j\,\PP(X_j=n-1\mid\mathcal F_j)\bigr| .
\]
For $w\in(0,1)$ let $\pi_w=\NB(\lambda/w,\,w)$, with $p$, $F$, $G$
those of $\pi_w$ and
\[
\Gamma_w=\sup_{n\ge0}\frac{F(n)\,G(n+1)}{(n+1)\,p_{n+1}},\qquad
\Lambda_w=\min\Bigl\{2\Gamma_w,\ \frac1\lambda\Bigr\};
\]
let $\pi_0=\Poi(\lambda)$, $\Gamma_0=\min\{1,1.4\lambda^{-1/2}\}$ and
$\Lambda_0=(1-e^{-\lambda})/\lambda$.  Then, for every $w\in[0,1)$,
\[
\dtv(W,\pi_w)\ \le\ \Gamma_w\,(M_w+\beta)\;+\;\Lambda_w\,Q,
\]
and $\beta=0$ when $X_j$ is independent of $\mathcal F_j$ for every
$j$.  If $X_j$ is independent of $\mathcal F_j$ for every $j$, and
for every $j$ the conditional law of $R_j$ given $(X_j,W_j)$ is
stochastically increasing in $X_j$ for every value of $W_j$, or
stochastically decreasing in $X_j$ for every value of $W_j$, then for
every $w\in[0,1)$
\[
\dtv(W,\pi_w)\ \le\ \Gamma_w\,M_w\;+\;\Lambda_w\sum_{j=1}^{c}(1-y_j)\,\bigl|\Cov(X_j,R_j)\bigr| .
\]
\end{mainthm}
This is Theorem~\ref{thm:aggnb} with Corollary~\ref{cor:monotone}.
It is the theorem of Arratia, Goldstein and Gordon \cite{AGG89} with
geometric summands in place of indicators and the negative binomial
in place of the Poisson law: $Q$ is their $b_2$ together with the
off-diagonal part of $b_1$, $\beta$ is their $b_3$, and the diagonal
part of $b_1$ is $M_0=\sum_jy_j\E X_j$.  The first bound is
Theorem~\ref{thm:coupling} with the worst coupling of the
neighborhood remainders $R_j$, the second with the monotone one; the
second is zero for independent neighbors and is the analogue of the
bound $(1-e^{-\lambda})\lambda^{-1}(\Var W-\lambda+2\sum_\alpha p_\alpha^{2})$
of \cite[Chapter~2]{BHJ92} for positively related indicators.
The target $\pi_w$ has mean $\lambda/(1-w)$, and the Poisson law is
the endpoint $w=0$ of the family, where $\Gamma_0$, $\Lambda_0$ are
the Poisson factors of \cite{AGG89,BHJ92}; for $w>0$ the bound
$1/\lambda$ on the first-difference factor is \cite[Theorem~2.10]{BX01}
(Theorem~\ref{thm:nbfactors}), so the dependence term carries the
factor of \cite{AGG89} throughout the family.  At equal parameters and $w=y_1$ the term $M_w$
vanishes and the bound prices the dependence alone.  Independent
summands at $w=\bar y$ recover Theorem~\ref{thm:A} with the weights
$y_j/(1-y_j)$ in place of $1/(1-y_j)$.

\subsection*{Overdispersed counts}
\begin{mainthm}[Poisson mixtures]\label{thm:B}
Let $\Lambda\ge0$ be a random variable with $m:=\E\Lambda\in(0,\infty)$
and $0<\Var\Lambda<\infty$, $v:=m+\Var\Lambda$, and let $W$ and
$W^{s}$ be the Poisson counts at the rate $\Lambda$ and at its
size-biased version,
\[
\PP(W=n)=\E\Bigl[e^{-\Lambda}\frac{\Lambda^{n}}{n!}\Bigr],\qquad
\PP(W^{s}=n)=\frac1m\,\E\Bigl[\Lambda\,e^{-\Lambda}\frac{\Lambda^{n}}{n!}\Bigr],
\qquad n\ge0 .
\]
Set $w=1-m/v$ and $\rho=m^{2}/(v-m)$, so that $\NB(\rho,w)$ has mean
$m$ and variance $v$, let $p$, $F$, $G$ be those of $\NB(\rho,w)$ and
$\Gamma_{\rho,w}=\sup_{n\ge0}F(n)G(n+1)/((n+1)p_{n+1})$, and let
$W^{b}$ have the law $w\,\mathcal L(1+W^{s})+(1-w)\,\mathcal L(W)$.
Then
\[
\dtv\bigl(W,\NB(\rho,w)\bigr)\ \le\ \min\Bigl\{2\,\Gamma_{\rho,w},\ \frac{v}{m^{2}}\Bigr\}\;m\;
\dwass\bigl(W^{b},W^{s}\bigr).
\]
The right side is a functional of the law of $\Lambda$ alone, and it
vanishes if and only if $\Lambda$ is Gamma with shape $\rho$.
\end{mainthm}
This is Corollary~\ref{cor:mixture}; $v/m^{2}=1/(\rho w)$ is the
first-difference bound of \cite[Theorem~2.10]{BX01}
(Theorem~\ref{thm:nbfactors}).  Claim counts with a random rate,
abundances across sites, word counts across documents: the
negative-binomial fit is certified by how far the mixing law is from
Gamma, in one Wasserstein distance, and the theorem needs no other
structure.

\subsection*{Counts observed through their differences}
\begin{mainthm}[Negative-binomial counts conditioned on their differences]\label{thm:C}
Let $R_1,\dots,R_r$ be independent, $R_i\sim\NB(c,q_i)$ with
$q_i\in(0,1)$, let $w=\prod_iq_i$ and $D_i=R_i-R_{i+1}$, $i<r$, and
let $\bkk=(k_1,\dots,k_{r-1})$ have offsets
$\kappa_i=\sum_{l\ge i}k_l\ge0$, $\kappa_r=0$.  With $\HNB(\bk,c,w)$
the law \eqref{eq:hnbrow} and
$\Zhat_c(\bk,w)=\sum_{n\ge0}\prod_{i=1}^{r}\binom{n+\kappa_i+c-1}{n+\kappa_i}w^{n}$
its normalizer,
\[
\mathcal L\bigl(R_r\bigm|D_i=k_i,\ i<r\bigr)=\HNB(\bk,c,w),\qquad
\PP\bigl(D_i=k_i,\ i<r\bigr)=\prod_{i=1}^{r}(1-q_i)^{c}q_i^{\kappa_i}\,\Zhat_c(\bk,w),
\]
\[
\Zhat_c(\bk,w)=\prod_{i=1}^{r-1}\binom{\kappa_i+c-1}{\kappa_i}\;
{}_rF_{r-1}\Bigl(\begin{matrix}c+\kappa_1,\dots,c+\kappa_{r-1},c\\ \kappa_1+1,\dots,\kappa_{r-1}+1\end{matrix};w\Bigr).
\]
As $c\to\infty$ with $cq_i/(1-q_i)\to\lambda_i\in(0,\infty)$,
$\HNB(\bk,c,w)$ converges in total variation to the law with
$\PP(n)\propto\prod_i\lambda_i^{\,n}/\prod_i(n+\kappa_i)!$; at $r=2$
this is the Bessel law of two independent Poisson counts conditioned
on their difference.
\end{mainthm}
This is Theorem~\ref{thm:target}, Proposition~\ref{prop:kemp} and
Proposition~\ref{prop:poissonlimit}.  Sites, cohorts or arms of a
trial each accumulate a negative-binomial count of shape $c$, a
Poisson count with a Gamma rate, and what is recorded is how far
each is ahead of the next; the remaining uncertainty is the level,
and the theorem gives its law exactly, with the probability of the
record in closed form.  The same law is the ensemble of $r$ species
in $c$ boxes that is canonical in the species differences and
grand-canonical in the level, and the smaller of two geometric clocks
given their difference at $r=2$, $c=1$.  When the shapes are only
approximately common, Theorem~\ref{thm:master} prices the
approximation through two rate moments and one coupling, with the
Stein factors of the law explicit (Proposition~\ref{prop:factorbound}).

\subsection*{Contingency tables}
\begin{mainthm}[The column direction of a contingency table]\label{thm:D}
Let $X_{ij}$, $i\le r$, $j\le c$, be independent,
$X_{ij}\sim\Geom(1-y_{ij})$ with $y_{ij}=q_i\beta_j\in(0,1)$, let
$U_i=\sum_jX_{ij}$, $D=(U_i-U_{i+1})_{i<r}$, and let
$\bkk=(k_1,\dots,k_{r-1})$ have offsets $\kappa_i=\sum_{l\ge i}k_l\ge0$,
$\kappa_r=0$, with $\PP(D=\bkk)>0$.  Let $W$ have the law
$\mathcal L(U_r\mid D=\bkk)$.  With $w=\prod_iq_i$,
$w_S=\prod_{i\in S}q_i$, $p$, $F$, $G$ those of the law
$\HNB(\bk,c,w)$ of \eqref{eq:hnbrow}, $D^{(L)}$ the difference vector
of the array with the columns in $L$ removed, and
\[
\begin{gathered}
\Gamma=\sup_{n\ge0}\frac{F(n)\,G(n+1)}{p_{n+1}\prod_{i=1}^{r}(n+1+\kappa_i)},\qquad
\nu^{b}_T=\prod_{i\in T}\Bigl(c+\sum_{j=1}^{c}\frac{y_{ij}}{1-y_{ij}}\Bigr),\\
\theta=\max_{L\subseteq[c]:\,|L|\le r}\ \sup_{s\in\Zint^{r-1}}
\frac{\PP(D^{(L)}=s)}{\PP(D=\bkk)},
\end{gathered}
\]
\[
\dtv\bigl(W,\HNB(\bk,c,w)\bigr)\ \le\ \Gamma\,\theta
\sum_{\varnothing\ne S\subseteq[r]}\nu^{b}_{S^c}\,w_S
\sum_{j=1}^{c}\ \prod_{i\in S}\frac1{1-y_{ij}}\ \bigl|\beta_j^{|S|}-1\bigr| .
\]
The bound is zero when $\beta_j\equiv1$: for identically distributed
columns the level of the row-sum profile is exactly $\HNB(\bk,c,w)$.
\end{mainthm}
This is Theorem~\ref{thm:colexpansion}, with the fiber constant
taken over all shifts $s$.  A uniform contingency table with row sums
$\br$ is a geometric array with exchangeable columns conditioned on
its row totals, and under the product-form representation unequal
column margins are the $\beta_j$; the bound is one sum over the
columns of terms that vanish with $|\beta_j-1|$.
Conditioning on both margins is the two-directional problem of
\S\ref{sec:rowcol}, where a matched row--column Stein system gives
the certificate of Theorem~\ref{thm:tablemaster}.

For contingency tables with heterogeneous margins at scale the
bound of Theorem~\ref{thm:D} is vacuous.  The reason is structural:
the moments $\nu_{S^c}\ge c^{|S^c|}$ and the weights
$1/(1-y_{ij})$ grow with the row means, while a sampler needs the ratio of neighboring fiber masses to a relative
precision of order one over the margins.  \S\ref{sec:rowcol} records
the obstruction, \S\ref{sec:computations} the row means at which the
bound passes one, and Open Problem~\ref{op:fullmargin} states what a
certificate for that regime would have to bound.

\paragraph{Relation to the literature.}  The one-row operator and
its characterization are those of Brown and Phillips \cite{BP99} and
Brown and Xia \cite{BX01}; the exact first factor is computed as in
\cite{BHJ92}, and the first-difference factor of \cite{BX01}, the
Gibbs-measure factors of Eichelsbacher and Reinert \cite{ER08}, the
stochastic-order derivation of Daly, Lef\`evre and Utev \cite{DLU12}
and the Wasserstein factors of Barbour, Gan and Xia \cite{BGX15} are
the nearest relatives of Theorem~\ref{thm:factors}; the
negative-binomial approximation of sums is treated in
\cite{BP99,VU09,CV15}, of runs in \cite{BX01,WX08}, and of dependent
indicators in \cite{UCV17}; each result of \S\ref{sec:nb} is placed
against these in a remark following it.  What we believe is new: the
two-biasings theorem with its exact cancellation at the matched rate
and the characterization $W^{b}+1\overset{d}{=}W^{d}$, the geometric
form of the Arratia--Goldstein--Gordon theorem with the Poisson law
as its endpoint, the
hyper-negative-binomial law as a named object---a member of Kemp's
generalized hypergeometric family \cite{Kemp68,JKK05} obtained by
conditioning on differences rather than on a sum---with the Bessel
law of Yuan and Kalbfleisch \cite{YK00} as its Poisson limit and a
Stein operator for that law as a corollary, and the column-wise
approximation of a random contingency table with an explicit error.

The paper is organized by the number of conditioned rows.
\S\ref{sec:nb} treats one row: the operator, its factors, the main
theorem through two biasings, and Theorems~\ref{thm:A},
\ref{thm:coupling}, \ref{thm:dep} and~\ref{thm:B}.  \S\ref{sec:hnb} treats general $r$: the family,
Theorem~\ref{thm:C}, the dependency expansion and Theorem~\ref{thm:D}.
\S\ref{sec:rowcol} conditions on both margins.  Sections
\ref{sec:nb}--\ref{sec:rowcol} are exact; \S\ref{sec:computations}
evaluates every bound against the exact distance it bounds and
records the cost of each formula, every number there regenerated,
with assertions, by the ancillary \path{verify_companion/}.
Open problems are in \S\ref{sec:open} and proofs not given in place
in \S\ref{sec:proofs}.

\section{Stein's method for the negative binomial}\label{sec:nb}\label{sec:stein}

One row first: $X_1,\dots,X_c$ i.i.d.\ $\Geom(1-w)$, with total
$\NB(c,w)$ by \eqref{eq:nbrow}.  This section builds the Stein
operator, its factors, a main theorem through the two biasings, and
the one-row face of the dependency expansion, and then proves
Theorems~\ref{thm:A}, \ref{thm:dep} and~\ref{thm:B}.

\subsection*{Negative binomial rows}

\begin{definition}[Negative binomial]\label{def:nb}
For $\rho>0$ and $w\in(0,1)$, $\NB(\rho,w)$ is the law on $\N$ with
$\PP(n)\propto\binom{n+\rho-1}{n}w^{\,n}$.
\end{definition}

\begin{lemma}[Geometric rows]\label{lem:geomrow}
If $X_1,\dots,X_c$ are i.i.d.\ $\Geom(1-q)$, $\PP(X_j=n)=(1-q)q^{n}$,
then $\sum_jX_j\sim\NB(c,q)$.
\end{lemma}

\paragraph{The target, and the standing assumptions.}  Throughout
this section $b,d:\N\to[0,\infty)$ are \emph{rates} with $d(0)=0$,
$b(n)>0$ for $n\ge0$, $d(n)>0$ for $n\ge1$, $b$ nondecreasing,
$d\le b$, and $d(n)/b(n)\to1$, $b(n{+}1)/b(n)\to1$; and $w\in(0,1)$.
The \emph{target} is the law $\pi=\pi_{b,d,w}$ on $\N$ with
\[
p_n=\frac1Z\prod_{k<n}\frac{w\,b(k)}{d(k{+}1)},\qquad
Z=\sum_{n\ge0}\prod_{k<n}\frac{w\,b(k)}{d(k{+}1)}<\infty,
\]
finite because the ratios tend to $w<1$; equivalently $\pi$ is the
law with the \emph{balance identity} $d(n{+}1)p_{n+1}=w\,b(n)p_n$ for
all $n\ge0$, and we write $t_n:=p_{n+1}/p_n=w\,b(n)/d(n{+}1)$, so
that $t_n\to w$.  The negative binomial is the case $b(k)=k+\rho$,
$d(k)=k$: $\NB(\rho,w)=\pi_{b,d,w}$, with $t_n=w(n{+}\rho)/(n{+}1)$
strictly decreasing when $\rho>1$ and constant when $\rho=1$.  The
hyper-negative-binomial law of \S\ref{sec:hnb} is the case of
degree-$r$ product rates; Lemma~\ref{lem:basic} verifies the
assumptions there, and every statement of this section holds for it
verbatim.

\begin{definition}[Two-sided bias]\label{def:bias}
Let $W$ be an $\N$-valued random variable with
$\beta_d:=\E\,d(W)\in(0,\infty)$ and $\beta_b:=\E\,b(W)<\infty$.  The
\emph{$d$-biased} and \emph{$b$-biased} versions $W^d$, $W^b$ have
\[
\PP(W^d=n)=\frac{d(n)\,\PP(W=n)}{\beta_d}\ \ (n\ge1),
\qquad
\PP(W^b=n)=\frac{b(n)\,\PP(W=n)}{\beta_b}\ \ (n\ge0).
\]
\end{definition}

\subsection*{The Stein operator and its factors}

\begin{definition}[Stein operator]\label{def:operator}
For bounded $g:\N\to\mathbb{R}$ with the convention $g(0)=0$ (the value
$g(0)$ is never used, since $d(0)=0$), set
\begin{equation}\label{eq:operator}
(\mathcal{A}\,g)(n)\;=\;w\,b(n)\,g(n+1)\;-\;d(n)\,g(n),
\qquad n\in\N ;
\end{equation}
for the negative binomial, $(\mathcal A g)(k)=w(k{+}\rho)g(k{+}1)-kg(k)$.
\end{definition}

\begin{theorem}[Characterization and exact solution; \cite{BP99,BX01}]\label{thm:solution}
Let $\pi=\pi_{b,d,w}$ be a target as above, $N\sim\pi$.
\begin{enumerate}
\item[(i)] An $\N$-valued random variable $M$ with $\E\,b(M)<\infty$
satisfies $\E(\mathcal{A}g)(M)=0$ for all bounded $g$ iff $M\sim\pi$.
\item[(ii)] For bounded $h$, the Stein equation
$(\mathcal{A}g)(n)=h(n)-\E h(N)=:\tilde h(n)$ has the unique
solution
\begin{equation}\label{eq:solution}
g_h(m+1)\;=\;\frac{1}{w\,b(m)\,p_m}\sum_{n\le m}\tilde h(n)\,p_n
\;=\;-\,\frac{1}{w\,b(m)\,p_m}\sum_{n> m}\tilde h(n)\,p_n,
\qquad m\ge0 .
\end{equation}
\end{enumerate}
\end{theorem}

Define the Stein factors
\[
\Gamma\;=\;\sup_{A\subseteq\N}\ \sup_{m\ge1}\,|g_A(m)|,
\qquad
\Lambda\;=\;\sup_{A\subseteq\N}\ \sup_{m\ge0}\,|g_A(m+1)-g_A(m)| ,
\]
where $g_A:=g_{\ind{\cdot\in A}}$.

\begin{theorem}[Stein factors]\label{thm:factors}
Let $\pi=\pi_{b,d,w}$ be a target as above with pmf $p$, distribution
function $F(m)=\sum_{n\le m}p_n$ and survival function
$G(m)=\sum_{n\ge m}p_n$; in (ii) assume moreover that $t_n$ is
strictly decreasing (for the negative binomial, $\rho>1$).  Set
$x_m'=w\,b(m-1)-d(m)$ and $y_m=d(m+1)-w\,b(m)$, so that
$t_{m-1}>1$ iff $x'_m>0$ and $t_{m+1}<1$ iff $y_{m+1}>0$.
\begin{enumerate}
\item[(i)] \emph{(Exact uniform factor.)}
\begin{equation}\label{eq:gammaexact}
\Gamma\;=\;\sup_{m\ge0}\ \frac{F(m)\,G(m+1)}{d(m+1)\,p_{m+1}} .
\end{equation}
\item[(ii)] \emph{(Closed-form bound.)}  Every $m\ge1$ satisfies
$t_{m-1}>1$ or $t_{m+1}<1$, and
\begin{equation}\label{eq:gammabound}
\Gamma\;\le\;\max\Bigl\{\ \frac{1}{w\,b(0)},\ \
\max_{m\ge1:\,x_m'>0}\frac{1}{x_m'},\ \
\sup_{m\ge0:\,y_{m+1}>0}\frac{d(m+2)}{d(m+1)\,y_{m+1}}\ \Bigr\},
\end{equation}
a maximum over an empty set being omitted; the first two terms are
finitely many, and when $d(m)\to\infty$ the terms of the third tend
to zero and its supremum is a maximum.
\item[(iii)] \emph{(First differences.)}
\begin{equation}\label{eq:lambdabound}
\Lambda\;\le\;\min\Bigl\{\,2\Gamma,\
\max\bigl\{1,\tfrac{1-w}{w}\bigr\}\,\Gamma+\frac{1}{w\,b(0)}\Bigr\}.
\end{equation}
\end{enumerate}
\end{theorem}

\begin{remark}[Related results]\label{rem:placement}
Theorem~\ref{thm:solution} is the generator characterization of
Brown and Xia \cite{BX01} for a birth--death target, with the
operator of Brown and Phillips \cite[Lemma~1]{BP99} at
$b(n)=n+c$, $d(n)=n$; the two forms of \eqref{eq:solution} are those
of \cite[Lemma~1.1.1]{BHJ92} for the Poisson law.  Part~(i) is exact
and is proved as in \cite[Lemma~1.1.1]{BHJ92}: the sets $[0,m]$ and
$[m+1,\infty)$ are extremal for $g_A(m+1)$.  The closed first
factors in the literature are of a different shape: \cite[Lemma~3]{BP99}
bounds $\Gamma$ independently of $c$, and Cloez and Delplancke
\cite[Proposition~4.5]{CD19} give, for $\NB(\rho,w)$,
$\Gamma\le\min\{1/(1-w),\ \sqrt\pi/\sqrt{(\rho+1)w(1-w)}\}$; neither
is exact, and (i) is the quantity every bound of this paper uses.
For the first difference, \cite[Theorem~2.10]{BX01}, in the form
recorded in \cite[(2.3)]{WX08}, gives for the negative binomial
$\sup_A|g_A(m+1)-g_A(m)|\le\min\{1/(w\,b(m)),\,1/d(m)\}$;
\cite{DLU12} recover it through stochastic orders,
\cite[Example~2.12]{ER08} treat the geometric case as a Gibbs
measure, and \cite[Theorem~1.1]{BGX15} give the factors in
Wasserstein distance.  Part~(iii) is rate-generic, valid for every
pair of the standing assumptions; all bounds available at the
negative-binomial pair are collected in Theorem~\ref{thm:nbfactors}.
In closed form at $c=1$:
$F(m)G(m+1)/((m+1)p_{m+1})=(1-w^{m+1})/((m+1)(1-w))$, the mean of
$1,w,\dots,w^{m}$, decreasing in $m$, so $\Gamma=1$ for every
geometric law.
\end{remark}

\begin{theorem}[Stein factors of the negative binomial: all bounds]\label{thm:nbfactors}
Let $\rho>0$, $w\in(0,1)$, let $p$, $F$, $G$ be those of $\NB(\rho,w)$,
and let $\Gamma_{\rho,w}$, $\Lambda_{\rho,w}$ be the Stein factors of
the pair $b(n)=n+\rho$, $d(n)=n$ (Definition~\ref{def:factors}).  Then
\[
\Gamma_{\rho,w}=\sup_{m\ge0}\frac{F(m)\,G(m+1)}{(m+1)\,p_{m+1}},
\qquad
\Gamma_{\rho,w}\le\min\Bigl\{\frac1{1-w},\ \frac{\sqrt\pi}{\sqrt{(\rho+1)w(1-w)}}\Bigr\},
\]
for $\rho>1$ also
\[
\begin{gathered}
\Gamma_{\rho,w}\le B_{\rho,w}=\max\Bigl\{\frac1{w\rho},\ \max_{m\ge1:\,x'_m>0}\frac1{x'_m},\
\sup_{m\ge0:\,y_{m+1}>0}\frac{m+2}{(m+1)\,y_{m+1}}\Bigr\},\\
x'_m=w(\rho-1)-m(1-w),\qquad y_{m+1}=(m+2)-w(m+1+\rho),
\end{gathered}
\]
and
\[
\Lambda_{\rho,w}\le\min\Bigl\{2\Gamma_{\rho,w},\ \frac1{w\rho},\
\frac2{1-w},\ \frac{1+w}{(1-w)^{2}},\ \frac{3}{2\sqrt{\rho w(1-w)^{3}}}\Bigr\}.
\]
The equality, $B_{\rho,w}$ (for $\rho>1$, where
$t_n=w(n+\rho)/(n+1)$ is strictly decreasing) and $2\Gamma_{\rho,w}$
are Theorem~\ref{thm:factors}; $1/(1-w)$ is \cite[Theorem~1.1]{BGX15}
and \cite[Proposition~4.5]{CD19}, and the third bound on $\Gamma$ is
\cite[Proposition~4.5]{CD19}; $1/(w\rho)$ is \cite[Theorem~2.10]{BX01},
in the form recorded in \cite[(2.3)]{WX08}; the last three bounds on
$\Lambda$ are \cite[Theorem~1.1]{BGX15}.
\end{theorem}

\begin{proof}
Theorem~\ref{thm:factors}(i)--(iii) at $b(n)=n+\rho$, $d(n)=n$ give
the equality, $B_{\rho,w}$ (part~(ii) requires $t_n$ strictly
decreasing, which is $\rho>1$) and $2\Gamma_{\rho,w}$.  The Stein equation
of \cite{BGX15} and \cite{CD19} is $w(n+\rho)g(n+1)-n\,g(n)=f(n)-\pi f$,
the equation of Definition~\ref{def:factors}; the suprema of
\cite{BGX15} are over $f$ with $|f(n+1)-f(n)|\le1$ and those of
\cite{CD19} over $0\le f\le1$, both classes containing the
indicators, so their bounds dominate $\Gamma_{\rho,w}$ and
$\Lambda_{\rho,w}$.  Theorem~2.10 of \cite{BX01} gives
$\sup_A|g_A(m+1)-g_A(m)|\le\min\{1/(w\,b(m)),\,1/d(m)\}$ for each
$m$, which is $1/(w\rho)$ at $m=0$ and at most $1/(w(\rho+m))<1/(w\rho)$
for $m\ge1$.
\end{proof}

\begin{remark}[Which bound where]\label{rem:which}
The equality costs one pass over the pmf and is the quantity every
bound of this paper uses; the closed forms are for statements uniform
in the parameters.  For $\Gamma_{\rho,w}$:
$\sqrt\pi/\sqrt{(\rho+1)w(1-w)}\le1/(1-w)$ if and only if
$(\rho+1)w\ge\pi(1-w)$, so the bound of \cite{CD19} is the smaller
closed form except at small means; $B_{\rho,w}\ge1/(w\rho)$, so
$B_{\rho,w}$ ($\rho>1$) can be the smallest only where
$1/(w\rho)\le\min\{1/(1-w),\sqrt\pi/\sqrt{(\rho+1)w(1-w)}\}$, that is
$\rho w\ge1-w$ and $(\rho+1)(1-w)\le\pi w\rho^{2}$.  For
$\Lambda_{\rho,w}$: the term $m=0$ of the equality gives
$\Gamma_{\rho,w}\ge F(0)G(1)/p_1=(1-p_0)/(w\rho)$ with $p_0=(1-w)^{\rho}$,
so $1/(w\rho)\le2\Gamma_{\rho,w}$ whenever $(1-w)^{\rho}\le\tfrac12$;
and $1/(w\rho)$ is below $2/(1-w)$, $(1+w)/(1-w)^{2}$ and
$3/(2\sqrt{\rho w(1-w)^{3}})$ if and only if $\rho w\ge(1-w)/2$,
$\rho w(1+w)\ge(1-w)^{2}$ and $\rho w\ge(1-w)^{3}/2.25$ respectively.
Hence $1/(w\rho)$ is the smallest of the five whenever
$(1-w)^{\rho}\le\tfrac12$ and $\rho w\ge\max\{(1-w)/2,\,(1-w)^{2}/(1+w)\}$,
the bounds of \cite{BGX15} take over at small means, and
$2\Gamma_{\rho,w}$ at small shapes $\rho$; the map on a grid is in
\S\ref{sec:computations}.  Of the five, $2\Gamma_{\rho,w}$ holds for every pair of the
standing assumptions and $B_{\rho,w}$ for every pair with $t_n$
strictly decreasing, in particular both for the
hyper-negative-binomial law of \S\ref{sec:hnb} with $c\ge2$.  In the Poisson regime $\rho\to\infty$,
$w\rho\to\lambda$, the bound $1/(w\rho)$ tends to $1/\lambda$, the
factor of \cite{AGG89}.
\end{remark}

\begin{remark}[Failure of the singleton sign pattern]\label{rem:signfail}
In \cite{DeS26} the sharp value of $\Lambda$ is computed exactly from the
singleton solutions $g_j:=g_{\{j\}}$, whose increments
$\Delta g_j(m)=g_j(m+1)-g_j(m)$ are nonpositive for all $m\ne j$; this
``sign pattern'' is equivalent to the monotonicity of
$U(m):=F(m)/(w\,b(m)p_m)$ and $V(m):=G(m+1)/(w\,b(m)p_m)$.  For
the degree-$r$ product rates of \S\ref{sec:hnb} the pattern
\emph{fails}.  Monotonicity of $U$ at $m=1$ requires
$d(1)\ge w\bigl(b(1)-b(0)\bigr)$; with $\bk=\mathbf{0}$ there one has
$d(1)=1$ while $b(1)-b(0)\sim rc^{r-1}$, so the inequality fails for all
large $c$ (at $r=2$, $c=40$, $w=\tfrac12$: $1<40.5$).  This is not
merely a failed proof route.  From the solution formula,
\[
g_2(1)=-\frac{p_2}{w\,b(0)},\qquad
g_2(2)=-\frac{p_2\,(p_0+p_1)}{w\,b(1)\,p_1},\qquad
\Delta g_2(1)=\frac{p_2\,\bigl(w\,b(1)-w\,b(0)-d(1)\bigr)}{w^{2}\,b(0)\,b(1)},
\]
so $\Delta g_2(1)>0$ exactly when $w\bigl(b(1)-b(0)\bigr)>d(1)$; on
$(r,c,\bk,w)=(2,6,(1,0),\tfrac{11}{20})$ this reads
$\tfrac{11}{20}\cdot14=\tfrac{77}{10}>2$.  Sharp $\Lambda$ is Open
Problem~\ref{op:lambda}; the increments and the factors on this
system are tabulated in \S\ref{sec:computations}.
\end{remark}

The factors $\Gamma,\Lambda$ of Theorem~\ref{thm:factors} make the
one-row biasing bound explicit.

\begin{theorem}[Approximation from the two biasings]\label{thm:nbmaster}
Let $b,d$ be rates as above, $W$ an $\N$-valued random variable with
$\beta_b=\E\,b(W)<\infty$ and $\beta_d=\E\,d(W)>0$, and $W^b,W^d$ its
biasings (Definition~\ref{def:bias}).  For every $w\in(0,1)$, with
$\Gamma,\Lambda$ the Stein factors of $\pi_{b,d,w}$
(Theorem~\ref{thm:factors}) and $\dwass$ the Wasserstein distance on
$\N$,
\[
\dtv\bigl(W,\pi_{b,d,w}\bigr)\;\le\;
\Lambda\,w\,\beta_b\,\dwass\bigl(W^b{+}1,\,W^d\bigr)
\;+\;\Gamma\,\bigl|\beta_d-w\,\beta_b\bigr| ;
\]
if $w^\ast:=\beta_d/\beta_b<1$, the choice $w=w^\ast$ annihilates the
second term.  For the negative binomial, $b-d\equiv\rho$ gives
$\beta_b=\beta_d+\rho=\E W+\rho$, so $w^\ast=\E W/(\E W+\rho)<1$ is
the parameter at which $\NB(\rho,w^\ast)$ has mean $\E W$,
$w^\ast\beta_b=\E W$, and at the match
\[
\dtv\bigl(W,\NB(\rho,w^\ast)\bigr)\;\le\;
\Lambda\,\E W\,\dwass\bigl(W^b{+}1,\,W^d\bigr):
\]
the moment input is the single number $\E W$.
\end{theorem}

\begin{remark}[The two biasings, in two lines]\label{rem:twolines}
For bounded $f$, $\E[w\,b(W)f(W+1)]=w\beta_b\,\E f(W^{b}+1)$ and
$\E[d(W)f(W)]=\beta_d\,\E f(W^{d})$, so
$\E(\mathcal Af)(W)=w\beta_b\,\E f(W^{b}+1)-\beta_d\,\E f(W^{d})$; the
theorem is this identity read through the factors, with
$\Gamma|\beta_d-w\beta_b|$ the rate mismatch and the Wasserstein term
the shape mismatch.  In particular $W\sim\pi_{b,d,w}$ for some $w\in(0,1)$ if and only if
$W^{b}+1\overset{d}{=}W^{d}$, and then $w=\beta_d/\beta_b$: taking
$g\equiv1$ on $\{1,2,\dots\}$ in the identity gives
$w\beta_b=\beta_d$, after which the two laws must agree; conversely
the identity vanishes for every $g$ and Theorem~\ref{thm:solution}(i)
applies.  For the negative binomial this is the analogue of the
Poisson identity $W^{s}\overset{d}{=}W+1$ \cite[Rem.~1.6]{Ross13}.  The
closest existing result is Ross's $r$-equilibrium theorem
\cite[Thm.~1.5]{Ross13}, $\dtv(W,\NB(r,p))\le c_{r,p}\,
\E|W^{*r}-W|$ with $c_{r,p}\le2$, one coupling to a fixed point of
a P\'olya-urn transformation; the discrete equilibrium coupling of
\cite{PRR13} is the $r=1$ case.  Here the coupling is between two
biasings of $W$ itself and the rate term cancels exactly at
$w^\ast$.
\end{remark}

\begin{remark}[How many moments, and why this is not anomalous]\label{rem:howmany}
The moment input of a birth--death bound is $\E[d(W)-w\,b(W)]$, so its
cost is the degree of the rate pair.  For a single row that degree is
$1$: the negative binomial consumes one moment, and so does the
Poisson case $b\equiv\lambda$, $d(k)=k$, where the condition is
$\E W=\lambda$.  The count rises only with the number of rows
conditioned: at the HNB rung $d(n)=\prod_{i\le r}(n+\kappa_i)$ has
degree $r$, so $\beta_d$ and $\beta_b$ require $\E W^r<\infty$ and the
match consumes moments $1,\dots,r$---two laws with equal means but
different second moments share a matched $w^\ast$ at $r=1$ and do not
at $r\ge2$.  The ``two
moments'' of \cite{AGG89} counts something else: the two dependency
sums $b_1,b_2$ over the summands and their neighborhoods, not
moments of the approximated variable.  In the present framing the two
statements agree---one moment per conditioned row---and the novelty of
the ladder is that the table costs $r$ of them, not that a single row
costs fewer.
\end{remark}

\begin{corollary}[Heterogeneous geometric sums]\label{cor:nbsum}
Let $X_j\sim\Geom(1-y_j)$, $j\le c$, be independent with
$y_j\in(0,1)$, $W=\sum_jX_j$, and fix any $w\in(0,1)$.  Then, with
$\Gamma$ the Stein factor of $\NB(c,w)$ (Theorem~\ref{thm:factors}),
\[
\dtv\bigl(W,\NB(c,w)\bigr)\ \le\ \Gamma\sum_{j=1}^{c}
\frac{|y_j-w|}{1-y_j}.
\]
\end{corollary}

\begin{proof}
This is the case $r=1$ of Theorem~\ref{thm:colexpansion}: the subset
sum has the single term $S=\{1\}$ with $\nu^b_\varnothing=1$, the
difference vector is empty so both probabilities in the fiber constant
of Definition~\ref{def:fiberconst} equal one and $\theta=1$, and
$w_S=w$.
\end{proof}

\begin{remark}[Related results]
Negative-binomial approximation of sums of nonnegative integer
variables is treated in \cite{BP99} for sums of indicators, by a
coupling that compares the law given one indicator with a geometric
variable, with the P\'olya urn as the example; in \cite{VU09} by
compound negative binomial laws; and in \cite{CV15} for $1$-dependent
summands through characteristic functions, with negative-binomial,
binomial and compound Poisson targets matched to two factorial
cumulants.  Those bounds carry moment or cumulant remainders that do
not vanish at a sum of i.i.d.\ geometrics.  Corollary~\ref{cor:nbsum}
carries the exact factor $\Gamma$ of \eqref{eq:gammaexact} and the
weights $1/(1-y_j)$, and is zero at $y_j\equiv w$.
\end{remark}

\subsection*{The coupling form, and dependent geometric sums}

\begin{definition}[The two biasings at a summand]\label{def:summandbias}
Let $X_1,\dots,X_c$ be $\N$-valued with $0<\E X_j<\infty$ and
$W=\sum_{j=1}^{c}X_j$.  For each $j$, $W^{(j)}$ and $W^{[j]}$ denote
random variables with
\[
\PP(W^{(j)}=n)=\frac{\E[X_j\ind{W=n}]}{\E X_j},\qquad
\PP(W^{[j]}=n)=\frac{\E[(X_j+1)\ind{W=n}]}{\E X_j+1},\qquad n\in\N .
\]
With $b(n)=n+c$, $d(n)=n$ and $W^{b},W^{d}$ as in
Definition~\ref{def:bias},
\[
\mathcal L(W^{d})=\sum_{j=1}^{c}\frac{\E X_j}{\E W}\,\mathcal L(W^{(j)}),\qquad
\mathcal L(W^{b})=\sum_{j=1}^{c}\frac{\E X_j+1}{\E W+c}\,\mathcal L(W^{[j]}),
\]
since $n\PP(W=n)=\sum_j\E[X_j\ind{W=n}]$ and
$(n+c)\PP(W=n)=\sum_j\E[(X_j+1)\ind{W=n}]$.
\end{definition}

\begin{theorem}[Coupling form]\label{thm:couplingbody}
Let $X_1,\dots,X_c$ be $\N$-valued with $0<\E X_j<\infty$,
$W=\sum_{j=1}^{c}X_j$, $y_j=\E X_j/(1+\E X_j)$ and
$\lambda=\sum_{j=1}^{c}y_j$.  For each $j$ let $(W^{(j)},W^{[j]})$ be
a pair of random variables on one probability space with the laws of
Definition~\ref{def:summandbias}.  For $w\in(0,1)$ let
$\pi_w=\NB(\lambda/w,w)$ with Stein factors $\Gamma_w,\Lambda_w$
(Definition~\ref{def:factors}); let $\pi_0=\Poi(\lambda)$, and let
$\Gamma_0,\Lambda_0$ be the suprema of $|g_A(m)|$ and
$|g_A(m+1)-g_A(m)|$ over $A\subseteq\N$, $m\ge0$, for the solutions
of $\lambda g_A(n+1)-n\,g_A(n)=\ind{n\in A}-\pi_0(A)$, $g_A(0)=0$.
Then, for every $w\in[0,1)$,
\[
\dtv(W,\pi_w)\ \le\ \Lambda_w\sum_{j=1}^{c}\E X_j\;
\E\bigl|W^{[j]}+1-W^{(j)}\bigr|\;+\;\Gamma_w\sum_{j=1}^{c}|w-y_j|\,\E X_j .
\]
\end{theorem}

\begin{proof}
Fix $A\subseteq\N$ and let $g=g_A$ solve the Stein equation of
$\pi_w$, $\mathcal A_wg(n)=(\lambda+wn)g(n+1)-n\,g(n)$, so that
$\|g\|_\infty\le\Gamma_w$ and $\|\Delta g\|_\infty\le\Lambda_w$.  By
the definition of the two laws,
$\E[X_jg(W)]=\E X_j\,\E g(W^{(j)})$ and
$\E[(X_j+1)g(W+1)]=(\E X_j+1)\,\E g(W^{[j]}+1)$, and
$y_j(\E X_j+1)=\E X_j$.  Since
$\lambda+wn=\sum_j\{(w-y_j)x_j+y_j(x_j+1)\}$ whenever $n=\sum_jx_j$,
\[
\E\mathcal A_wg(W)
=\sum_{j=1}^{c}\E X_j\Bigl\{\E g(W^{[j]}+1)-\E g(W^{(j)})\Bigr\}
+\sum_{j=1}^{c}(w-y_j)\,\E X_j\,\E g(W^{(j)}+1),
\]
and $|\E g(W^{[j]}+1)-\E g(W^{(j)})|\le\Lambda_w\,\E|W^{[j]}+1-W^{(j)}|$
on the common probability space.  Take the supremum over $A$.
\end{proof}

\begin{remark}[Independent summands]
If $X_j$ is independent of $(X_i)_{i\ne j}$ then, by
Lemma~\ref{lem:sizebias} at $c=1$ and Lemma~\ref{lem:geomchar}, a
geometric $X_j$ admits the coupling $W^{(j)}=W+1+X_j'$,
$W^{[j]}=W+X_j'$ with one independent copy $X_j'$, so the $j$-th
coupling term vanishes.  With all summands independent and geometric
the first sum is zero and Theorem~\ref{thm:couplingbody} reads
$\dtv(W,\NB(\lambda/w,w))\le\Gamma_wM_w$ for every $w\in(0,1)$,
$M_w=\sum_j|w-y_j|\,y_j/(1-y_j)$; at $w=\bar y=\lambda/c$ the target
is $\NB(c,\bar y)$ and the bound improves the weights $1/(1-y_j)$ of
Corollary~\ref{cor:nbsum} at that $w$, the two proofs differing in
the identity $\sum_j(\bar y-y_j)\E[(X_j+1)g(W+1)]
=\sum_j(\bar y-y_j)\E[X_jg(W+1)]$, which holds because
$\sum_j(\bar y-y_j)=0$.  For a general independent summand the $j$-th
coupling term is $\E X_j\,\dwass(X_j^{+}+1,X_j^{s})$, $X_j^{s}$ the
size-biased and $X_j^{+}$ the $(x+1)$-biased summand; for
$X_j\sim\mathrm{Bernoulli}(p_j)$ it is $2p_j^{2}/(1+p_j)$.
\end{remark}

\begin{theorem}[Dependent geometric sums]\label{thm:aggnb}
Let $X_1,\dots,X_c$ be $\N$-valued with $X_j\sim\Geom(1-y_j)$,
$y_j\in(0,1)$, $W=\sum_{j=1}^{c}X_j$, $\lambda=\sum_{j=1}^{c}y_j$.  For
each $j$ fix $B_j\subseteq[c]$ with $j\in B_j$, and let
$V_j=\sum_{i\in B_j}X_i$, $W_j=W-V_j$,
$\mathcal F_j=\sigma(X_i:i\notin B_j)$.  Set
\[
M_w=\sum_{j=1}^{c}|w-y_j|\,\frac{y_j}{1-y_j},\qquad
Q=\sum_{j=1}^{c}\sum_{i\in B_j\setminus\{j\}}\Bigl((1+y_j)\,\E[X_iX_j]+y_j\,\E X_i\Bigr),
\]
\[
\beta=\sum_{j=1}^{c}\E\sum_{n\ge1}n\,\bigl|\PP(X_j=n\mid\mathcal F_j)-y_j\,\PP(X_j=n-1\mid\mathcal F_j)\bigr| .
\]
For $w\in(0,1)$ let $\pi_w=\NB(\lambda/w,w)$, the target of the pair
$b(n)=n+\lambda/w$, $d(n)=n$, with Stein factors $\Gamma_w,\Lambda_w$
(Definition~\ref{def:factors}); let $\pi_0=\Poi(\lambda)$, and let
$\Gamma_0,\Lambda_0$ be the suprema of $|g_A(m)|$ and
$|g_A(m+1)-g_A(m)|$ over $A\subseteq\N$, $m\ge0$, for the solutions
of $\lambda g_A(n+1)-n\,g_A(n)=\ind{n\in A}-\pi_0(A)$, $g_A(0)=0$.
Then, for every $w\in[0,1)$,
\[
\dtv(W,\pi_w)\ \le\ \Gamma_w\,(M_w+\beta)+\Lambda_w\,Q .
\]
If $X_j$ is independent of $\mathcal F_j$ for every $j$ then $\beta=0$.
For $w>0$, $\Gamma_w$ and $\Lambda_w$ obey Theorem~\ref{thm:nbfactors}
with $\rho=\lambda/w$, in particular $\Lambda_w\le\min\{2\Gamma_w,1/\lambda\}$;
$\Gamma_0\le\min\{1,1.4\lambda^{-1/2}\}$ and
$\Lambda_0\le(1-e^{-\lambda})/\lambda$ \cite{AGG89,BHJ92}.
\end{theorem}

\begin{proof}
Fix $A\subseteq\N$ and let $g=g_A$ solve the Stein equation of
$\pi_w$, so that $\|g\|_\infty\le\Gamma_w$ and
$\|\Delta g\|_\infty\le\Lambda_w$; for every $w\in[0,1)$ the operator
is
\[
\mathcal A_wg(n)=(\lambda+wn)\,g(n+1)-n\,g(n),
\]
since $w\,b(n)=w(n+\lambda/w)=\lambda+wn$.  Put
$\delta_j=\E[X_jg(W)]-y_j\,\E[(X_j+1)g(W+1)]$.  Because
$\lambda=\sum_jy_j$,
\[
\PP(W\in A)-\pi_w(A)=\E\mathcal A_wg(W)
=\sum_{j=1}^{c}\Bigl\{(w-y_j)\,\E[X_jg(W+1)]-\delta_j\Bigr\},
\]
and $|(w-y_j)\,\E[X_jg(W+1)]|\le\Gamma_w\,|w-y_j|\,y_j/(1-y_j)$.
Write $W=W_j+V_j$ and split
\[
\begin{aligned}
\delta_j&=\E\bigl[X_j\{g(W_j+V_j)-g(W_j+X_j)\}\bigr]
-y_j\,\E\bigl[(X_j+1)\{g(W_j+V_j+1)-g(W_j+X_j+1)\}\bigr]\\
&\qquad+\Bigl(\E[X_jg(W_j+X_j)]-y_j\,\E[(X_j+1)g(W_j+X_j+1)]\Bigr).
\end{aligned}
\]
Since $V_j-X_j=\sum_{i\in B_j\setminus\{j\}}X_i\in\N$, each of the
first two braces is at most $\Lambda_w(V_j-X_j)$ in absolute value,
so those two terms contribute at most
$\Lambda_w\sum_{i\in B_j\setminus\{j\}}\bigl(\E[X_iX_j]+y_j\E[X_i(X_j+1)]\bigr)$,
the $j$-th term of $\Lambda_wQ$.  In the third, $W_j$ is
$\mathcal F_j$-measurable, and conditioning on $\mathcal F_j$,
\[
\begin{aligned}
&\E[X_jg(W_j+X_j)]-y_j\,\E[(X_j+1)g(W_j+X_j+1)]\\
&\qquad=\E\sum_{n\ge1}g(W_j+n)\,n\,\bigl\{\PP(X_j=n\mid\mathcal F_j)-y_j\,\PP(X_j=n-1\mid\mathcal F_j)\bigr\},
\end{aligned}
\]
at most $\Gamma_w$ times the $j$-th term of $\beta$.  When $X_j$ is
independent of $\mathcal F_j$ the conditional law is $\Geom(1-y_j)$,
with $p_n=y_jp_{n-1}$, and every term vanishes.  Summing over $j$ and
taking the supremum over $A$ gives the bound.
\end{proof}

\begin{corollary}[Monotone neighborhoods]\label{cor:monotone}
In the setting of Theorem~\ref{thm:aggnb}, suppose $X_j$ is
independent of $\mathcal F_j$ for every $j$, and write
$R_j=V_j-X_j=\sum_{i\in B_j\setminus\{j\}}X_i$.  Suppose that for
every $j$ there is a sign $s_j\in\{+1,-1\}$ such that, for all
$x<x'$, all $u$ and all $r$,
\[
s_j\bigl[\PP(R_j>r\mid X_j=x',\,W_j=u)-\PP(R_j>r\mid X_j=x,\,W_j=u)\bigr]\ \ge\ 0 .
\]
Then, for every $w\in[0,1)$,
\[
\dtv(W,\pi_w)\ \le\ \Gamma_w\,M_w\;+\;\Lambda_w\sum_{j=1}^{c}(1-y_j)\,\bigl|\Cov(X_j,R_j)\bigr| .
\]
\end{corollary}

\begin{proof}
Apply Theorem~\ref{thm:couplingbody}; the second sum there is
$\Gamma_wM_w$.  Fix $j$.  Since $X_j$ is independent of
$\mathcal F_j$, the law of $(X_j,W_j,R_j)$ is
$\PP(X_j=x)\,\PP(W_j\in du)\,\PP(R_j\in dr\mid X_j=x,W_j=u)$.  Let
$X^{\ast}$ have the law $(x+1)\PP(X_j=x)/(\E X_j+1)$; for a geometric
$X_j$, $X^{\ast}+1$ has the law $x\PP(X_j=x)/\E X_j$, since
$(x+1)\PP(X_j=x+1)/\E X_j=(x+1)(1-y_j)\PP(X_j=x)$.  Draw $W_j$ from
its law, independently of $X^{\ast}$, and given $(X^{\ast},W_j)=(x,u)$
draw $R^{(j)}$ from $\PP(R_j\in\cdot\mid X_j=x+1,W_j=u)$ and
$R^{[j]}$ from $\PP(R_j\in\cdot\mid X_j=x,W_j=u)$, coupled so that
$s_j(R^{(j)}-R^{[j]})\ge0$, which the hypothesis permits.  Then
$W^{(j)}=X^{\ast}+1+W_j+R^{(j)}$ and $W^{[j]}=X^{\ast}+W_j+R^{[j]}$
have the laws of Definition~\ref{def:summandbias}, and
\[
\begin{aligned}
\E\bigl|W^{[j]}+1-W^{(j)}\bigr|&=\E\bigl|R^{(j)}-R^{[j]}\bigr|
=\bigl|\E R^{(j)}-\E R^{[j]}\bigr|\\
&=\Bigl|\frac{\E[X_jR_j]}{\E X_j}-\frac{\E[(X_j+1)R_j]}{\E X_j+1}\Bigr|
=\frac{|\Cov(X_j,R_j)|}{\E X_j\,(\E X_j+1)} .
\end{aligned}
\]
Multiplying by $\E X_j$ and using $\E X_j+1=1/(1-y_j)$ gives the
$j$-th term.
\end{proof}

\begin{proposition}[The first-difference form]\label{prop:aggnblambda}
In the setting of Theorem~\ref{thm:aggnb}, suppose $X_j$ is
independent of $\mathcal F_j$ for every $j$, and let
\[
w^\ast=\frac{\sum_jy_j\,\E X_j}{\sum_j\E X_j},\qquad
M^\ast=\sum_{j=1}^{c}|w^\ast-y_j|\Bigl(\frac{y_j}{(1-y_j)^{2}}
+\sum_{i\in B_j\setminus\{j\}}\bigl(\E[X_iX_j]+\E X_i\,\E X_j\bigr)\Bigr).
\]
Then $\dtv(W,\pi_{w^\ast})\le\Lambda_{w^\ast}\,(M^\ast+Q)$.
\end{proposition}

\begin{remark}[The Poisson endpoint]\label{rem:nbpoisson}
For $N\sim\NB(\rho,w)$ and $\mu=\rho w/(1-w)$ the balance identity
gives $\E[N g(N)]=w\,\E[(N+\rho)g(N+1)]$, hence
$\E[\mu g(N+1)-Ng(N)]=w\,\E[(\mu-N)g(N+1)]$, and with
$\|g\|_\infty\le\min\{1,1.4\mu^{-1/2}\}$ and $\Var N=\mu/(1-w)$,
\[
\dtv\bigl(\NB(\rho,w),\Poi(\mu)\bigr)\ \le\ w\sqrt{\Var N}\,\min\{1,1.4\mu^{-1/2}\}\ \le\ \frac{1.4\,w}{\sqrt{1-w}} .
\]
For $N'\sim\Poi(\mu)$ the Stein equation of $\Poi(\lambda)$ gives
$\E[\lambda g(N'+1)-N'g(N')]=(\lambda-\mu)\,\E g(N'+1)$, hence
$\dtv(\Poi(\mu),\Poi(\lambda))\le\min\{1,1.4\lambda^{-1/2}\}\,|\mu-\lambda|$.
Along the family of Theorem~\ref{thm:aggnb}, $\mu=\lambda/(1-w)$, so
\[
\dtv(\pi_w,\pi_0)\ \le\ \frac{1.4\,w}{\sqrt{1-w}}+\min\{1,1.4\lambda^{-1/2}\}\,\frac{\lambda w}{1-w}
\ \longrightarrow\ 0\qquad(w\downarrow0).
\]
At $w=0$ the
theorem reads
$\dtv(W,\Poi(\lambda))\le\Gamma_0\bigl(\sum_jy_j\E X_j+\beta\bigr)+\Lambda_0Q$,
the bound of \cite{AGG89} with geometric summands, the diagonal part
of their $b_1$ carrying $\Gamma_0$; Proposition~\ref{prop:aggnblambda}
moves it onto $\Lambda_{w^\ast}$ at the target $\pi_{w^\ast}$.
\end{remark}

\begin{remark}[Related results]
Theorem~\ref{thm:aggnb} at $w=0$ is the analogue for geometric
marginals of \cite[Theorem~1]{AGG89}, with the geometric balance in
place of $\E[X_\alpha f(W)]=p_\alpha\E[f(W)\mid X_\alpha=1]$: $Q$ is
their $b_2$ with the off-diagonal part of $b_1$, $\beta$ their $b_3$,
and $M_0$ the diagonal part of $b_1$.  Theorem~\ref{thm:couplingbody}
is the coupling form of the same method \cite{AGG89,BHJ92}.  For
indicators $W^{(\alpha)}=1+V_\alpha$ with
$\mathcal L(V_\alpha)=\mathcal L(W-X_\alpha\mid X_\alpha=1)$,
$\mathcal L(W^{[\alpha]})=(1+p_\alpha)^{-1}\{\mathcal L(W)+p_\alpha\mathcal L(1+V_\alpha)\}$,
and the first sum is at most
$\sum_\alpha p_\alpha\E|W-V_\alpha|/(1+p_\alpha)+\sum_\alpha p_\alpha^{2}/(1+p_\alpha)$,
their $\sum_\alpha p_\alpha\E|W-V_\alpha|$ up to a diagonal term; the
target is $\Poi(\sum_\alpha p_\alpha/(1+p_\alpha))$ rather than
$\Poi(\sum_\alpha p_\alpha)$, and the mean defect
$\sum_\alpha p_\alpha^{2}/(1+p_\alpha)$ is the second sum, with the
factor $\Gamma_0$.  Corollary~\ref{cor:monotone} is the
analogue of the bounds for positively and negatively related
indicators of \cite[Chapter~2]{BHJ92}, where the monotone coupling
turns $\sum_\alpha p_\alpha\E|W-V_\alpha|$ into
$|\Var W-\lambda|$ up to $2\sum_\alpha p_\alpha^{2}$; here the
covariance sum $\sum_j\Cov(X_j,R_j)$ is $\Var W-\sum_j\Var X_j$ when
the neighborhoods carry all the dependence, and the diagonal term is
absent.
Negative-binomial approximation of dependent sums is treated in
\cite[Theorem~4.2]{BX01} for $2$-runs and in \cite{WX08} for
$k$-runs, through the neighborhood moments of the run indicators
and the factor $1/(w\rho)$ of \cite[Theorem~2.10]{BX01}; in
\cite{CV15} for $1$-dependent nonnegative integer summands by
characteristic functions; and in \cite{UCV17} for sums of dependent
indicators by perturbations of the negative-binomial and binomial
operators.  Each of these prices a self-term of the summands.  Here
the marginals are geometric, the geometric balance removes the
self-term, and the bound is a function of the cross-moments and the
conditional defects alone: zero at independence and equal
parameters, with the Poisson law as the endpoint $w=0$ of the same
statement.  With $\Lambda_w\le1/\lambda$ from
\cite[Theorem~2.10]{BX01} (Theorem~\ref{thm:nbfactors},
$w\rho=\lambda$), the dependence term $Q$ carries the factor
$1/\lambda$ for every $w\in(0,1)$, as it does in \cite{AGG89} at
$w=0$.
\end{remark}

\subsection*{Applications}

\emph{Failures before the $c$-th success at drifting odds.}  In
Bernoulli trials whose success probability changes from stage to
stage---$1-y_j$ at stage $j$---the number of failures before the
$c$-th success is $W$ of Corollary~\ref{cor:nbsum}.  With
$w=\bar y$ the mean of the $y_j$ and $\eta=\max_j|y_j-\bar y|$, the
law is within $\Gamma c\,\eta/(1-\max_jy_j)$ of $\NB(c,\bar y)$ in
total variation: the negative-binomial analogue of the Poisson
approximation to the Poisson-binomial (for the P\'olya side of that
analogy see \cite{Sk12}), with the spread of the odds in the role of
$\sum p_j^{2}$.  One stage displaced to $y_1=w+\eta$ gives the bound
$\Gamma_{c,w}\,|\eta|/(1-w-\eta)$.  A drift with $\sum_j(y_j-w)=0$
cancels the first-order variation of the law of $W$, since that law
is a symmetric function of $(y_1,\dots,y_c)$ and its gradient at
$(w,\dots,w)$ is a multiple of $(1,\dots,1)$, while the bound remains
first order in the spread.  A second-order bound for that case would
come from the signed-measure expansions of \cite{BC02,CR06}; none is
written for the negative binomial, and Open
Problem~\ref{op:secondorder} records it.  Both faces are evaluated
in \S\ref{sec:computations}.

\begin{corollary}[A single species in inhomogeneous boxes]\label{cor:occupancy}
Let indistinguishable balls fill $c$ distinguishable boxes
independently and geometrically, box $j$ at fugacity
$y_j=e^{-\beta\epsilon_j}$, and let $W$ be the total number of balls.
Then for every $\bar\epsilon$,
\[
\dtv\bigl(W,\NB(c,e^{-\beta\bar\epsilon})\bigr)\ \le\
\Gamma_{c,e^{-\beta\bar\epsilon}}\sum_{j=1}^{c}
\frac{|e^{-\beta\epsilon_j}-e^{-\beta\bar\epsilon}|}{1-e^{-\beta\epsilon_j}} .
\]
\end{corollary}

\begin{proof}
Corollary~\ref{cor:nbsum} with $y_j=e^{-\beta\epsilon_j}$ and
$w=e^{-\beta\bar\epsilon}$.
\end{proof}

The weights $1/(1-y_j)$ are the box means plus one: the replacement
of a spectrum by its mean is sharp when the boxes hold few balls and
priced by the occupancy when they hold many.  This is the
grand-canonical ensemble whose Gibbs approximation \cite{ER08}
treats; the negative-binomial constant is the addition.  At the mean
fugacity $e^{-\beta\bar\epsilon}=\bar y$ one has $\sum_j(y_j-\bar y)=0$,
the second face above.

\begin{corollary}[Random stage parameters]\label{cor:randomstage}
If the $y_j$ are random, independent of the trials, with
$\E|y_j-w|\le\eta$ and $y_j\le y^{+}<1$, then the conditional law of
$W$ given $(y_1,\dots,y_c)$ satisfies
$\E\,\dtv\bigl(\mathcal L(W\mid y_1,\dots,y_c),\NB(c,w)\bigr)\le\Gamma_{c,w}\,c\,\eta/(1-y^{+})$.
\end{corollary}

\begin{proof}
Corollary~\ref{cor:nbsum} conditionally on the parameters, then the
expectation.
\end{proof}

Inspection until $c$ defectives with lot-to-lot defect rates, and
waiting times with a randomly varying success rate, are of this
form.

\begin{corollary}[Poisson mixtures]\label{cor:mixture}
Let $\Lambda\ge0$ be a random variable with $m:=\E\Lambda\in(0,\infty)$
and $\Var\Lambda<\infty$, $v:=m+\Var\Lambda$, and let $W$, $W^{s}$ be
the Poisson counts at the rate $\Lambda$ and at its size-biased
version,
$\PP(W=n)=\E[e^{-\Lambda}\Lambda^{n}/n!]$,
$\PP(W^{s}=n)=\E[\Lambda e^{-\Lambda}\Lambda^{n}/n!]/m$, $n\ge0$.
For every $\rho>0$, with $w^\ast=m/(m+\rho)$ and $b(n)=n+\rho$,
$d(n)=n$,
\[
\mathcal L(W^{d})=\mathcal L(1+W^{s}),\qquad
\mathcal L(W^{b})=w^\ast\,\mathcal L(1+W^{s})+(1-w^\ast)\,\mathcal L(W),
\]
and
\[
\dtv\bigl(W,\NB(\rho,w^\ast)\bigr)\ \le\ \Lambda_{\rho,w^\ast}\;m\;\dwass\bigl(W^{b},W^{s}\bigr)
\ \le\ \min\Bigl\{2\,\Gamma_{\rho,w^\ast},\ \frac1{\rho w^\ast}\Bigr\}\;m\;\dwass\bigl(W^{b},W^{s}\bigr),
\]
with $\Gamma_{\rho,w^\ast},\Lambda_{\rho,w^\ast}$ the factors of
Theorem~\ref{thm:nbfactors}.  When $\Var\Lambda>0$ the choice
$w^\ast=1-m/v$, $\rho=m^{2}/(v-m)$ matches the mean and the variance,
and then $1/(\rho w^\ast)=v/m^{2}$.  The right side is zero if and
only if $\Lambda$ is Gamma with shape $\rho$.
\end{corollary}

\begin{proof}
$\E W=\E\Lambda=m$.  For $n\ge1$,
$n\,\PP(W=n)=\E[\Lambda e^{-\Lambda}\Lambda^{n-1}/(n-1)!]=m\,\PP(W^{s}=n-1)$,
so $W^{d}\overset{d}{=}1+W^{s}$, and
$(n+\rho)\PP(W=n)=m\,\PP(1+W^{s}=n)+\rho\,\PP(W=n)$ gives $W^{b}$ with
$\beta_b=m+\rho$.  Theorem~\ref{thm:nbmaster} at $w^\ast=m/(m+\rho)$
gives the first bound with
$\dwass(W^{b}+1,W^{d})=\dwass(W^{b}+1,1+W^{s})=\dwass(W^{b},W^{s})$,
and Theorem~\ref{thm:nbfactors} gives
$\Lambda_{\rho,w^\ast}\le\min\{2\Gamma_{\rho,w^\ast},1/(\rho w^\ast)\}$.
If $\Lambda\sim\mathrm{Gamma}(\rho,\cdot)$ then $W\sim\NB(\rho,w^\ast)$
and $W^{b}+1\overset{d}{=}W^{d}$ by Remark~\ref{rem:twolines}.
Conversely, if the distance is zero then $W\sim\NB(\rho,w^\ast)$ by
the same characterization, and the Poisson mixture determines its
mixing law, since $\E s^{W}=\E e^{-(1-s)\Lambda}$ is the Laplace
transform of $\Lambda$; the negative binomial is the Poisson--Gamma
mixture, so $\Lambda$ is Gamma.  For the matched choice, $\NB(\rho,w)$
has mean $\rho w/(1-w)=m$ and variance $m/(1-w)=v$.
\end{proof}

The corollary needs no structure in $W$ beyond the mixture: claim
counts with a random rate, abundances across sites, word counts
across documents.  The Delaporte law of insurance claim counts---a
negative binomial convolved with a Poisson \cite{JKK05}---is the
Poisson mixture whose rate is a Gamma variable plus a constant, and
sits inside the same bound.

\begin{remark}[Related results]
The Poisson approximation of a mixed Poisson law is
\cite[Theorem~1.C]{BHJ92}, $\dtv(W,\Poi(\E\Lambda))\le\Var\Lambda/\E\Lambda$,
with the Wasserstein analogue in \cite[\S4.1]{CD19}; the error is the
variance of the rate.  Corollary~\ref{cor:mixture} takes the
negative binomial carrying that variance as the target, and its
error is the distance of the mixing law from the Gamma law in one
Wasserstein distance, zero exactly when $\Lambda$ is Gamma.  A second
bound on the same distance follows from \cite[Theorem~1.C]{BHJ92}
and Remark~\ref{rem:nbpoisson} by the triangle inequality,
\[
\dtv\bigl(W,\NB(\rho,w^\ast)\bigr)\ \le\ \frac{\Var\Lambda}{m}+\frac{1.4\,w^\ast}{\sqrt{1-w^\ast}},
\qquad w^\ast=\frac{\Var\Lambda}{m+\Var\Lambda};
\]
it is of first order in $\Var\Lambda$ for every mixing law, while the
corollary vanishes at the Gamma laws; the two are compared on examples
in \S\ref{sec:computations}.
\end{remark}

\emph{The Poisson regime.}  As $c\to\infty$ with $cw/(1-w)\to\lambda$,
$\NB(c,w)\to\mathrm{Poisson}(\lambda)$
(Proposition~\ref{prop:poissonlimit} at $r=1$), and
Theorem~\ref{thm:nbmaster} runs the other way: a nearly-Poisson
count with excess variance is within the theorem's bound of the
negative binomial that carries that variance.  The two together say
when overdispersion is worth modeling.

\emph{Where it refuses: the coupon collector.}  The number of draws
to collect all $c$ coupons, less $c$, is a sum of geometrics with
$y_j=j/c$, $0\le j<c$; the spread is of order one and the weights
$1/(1-y_j)=c/(c-j)$ sum to $cH_{c}$, so the right side of
Corollary~\ref{cor:nbsum} is of order $\Gamma_{c,w}\,c\log c$ and
exceeds one for all large $c$.  The law is not negative
binomial---its fluctuations are Gumbel \cite{ER61}---and the bound
says so rather than pretending otherwise
(\S\ref{sec:computations}).

\section{The hyper-negative-binomial law: general \texorpdfstring{$r$}{r}}\label{sec:hnb}\label{sec:target}\label{sec:steinsystem}\label{sec:dependency}\label{sec:master}

Now $r$ rows: $X_{ij}\sim\Geom(1-q_i)$ independent, every row
identically distributed within itself, with independent row totals
$U_i\sim\NB(c,q_i)$.  Given the differences $U_i-U_{i+1}$ one level is
free, and its law is the hyper-negative-binomial law
\eqref{eq:hnbrow}: negative binomial at $r=1$, exact whenever the
columns are exchangeable, and the target of a dependency expansion
whose mismatch terms price column heterogeneity.  The instruments of
\S\ref{sec:nb} carry over with the birth polynomial of degree $r$ in
place of the linear rate.

\subsection*{The hyper negative binomial: conditioning on the row differences}

\begin{definition}[Hyper-negative-binomial law]\label{def:hnb}
For $r\ge1$, $c\ge1$, offsets $\bk=(\kappa_1,\dots,\kappa_r)\in\N^r$
with $\kappa_r=0$ (when the offsets come from a difference vector
$\bkk=(k_1,\dots,k_{r-1})\in\Zint^{r-1}$, $\kappa_i=\sum_{l=i}^{r-1}k_l$,
required nonnegative), and $w\in(0,1)$, let $\HNB(\bk,c,w)$ be the
law on $\N$ with
\begin{equation}\label{eq:hnbpmf}
p_n \;=\; \frac{1}{\Zhat_c(\bk,w)}\prod_{i=1}^r\binom{n+\kappa_i+c-1}{n+\kappa_i}\,w^n,
\qquad
\Zhat_c(\bk,w)\;=\;\sum_{n\ge0}\prod_{i=1}^r\binom{n+\kappa_i+c-1}{n+\kappa_i}\,w^n .
\end{equation}
The coefficient of $w^n$ grows like a constant times $n^{r(c-1)}$, so the
series converges precisely for $|w|<1$.
\end{definition}

\begin{theorem}[Conditional representation and exact fiber probability]\label{thm:target}
Let $R_i\sim\NB(c,q_i)$, $1\le i\le r$, be independent, $q_i\in(0,1)$,
let $w=\prod_iq_i$, and let $\bkk\in\Zint^{r-1}$ have nonnegative
offsets $\kappa_i=\sum_{l\ge i}k_l$, $\bk=(\kappa_1,\dots,\kappa_{r-1},0)$
(relabel the rows so that $R_r$ is the smallest total on the event
considered).  Then, with
$D_i=R_i-R_{i+1}$,
\[
\mathcal{L}\bigl(R_r \,\bigm|\, D_i=k_i,\ 1\le i\le r-1\bigr)
\;=\;\HNB(\bk,c,w),
\]
and the fiber probability is exactly
\begin{equation}\label{eq:fiber}
\PP\bigl(D_i=k_i\ \forall i\bigr)
\;=\;\prod_{i=1}^r(1-q_i)^c q_i^{\kappa_i}\;\Zhat_c(\bk,w).
\end{equation}
In particular the conditioned law depends on $(q_1,\dots,q_r)$ only
through $w$.
\end{theorem}

The column sums are free at this rung; the row--column system of
\S\ref{sec:rowcol} pins them affinely, the exact target margins
becoming a single atom of the fiber shadow
(Proposition~\ref{prop:atom}).

\begin{proposition}[Hypergeometric normalizer; Kemp identification]\label{prop:kemp}
With $\kappa_r=0$,
\begin{equation}\label{eq:zhathyper}
\Zhat_c(\bk,w)\;=\;\prod_{i=1}^{r-1}\binom{\kappa_i+c-1}{\kappa_i}\cdot
{}_rF_{r-1}\!\left(\begin{matrix}c+\kappa_1,\dots,c+\kappa_{r-1},\,c\\
\kappa_1+1,\dots,\kappa_{r-1}+1\end{matrix};\,w\right).
\end{equation}
For $r=2$ with $k=\kappa_1$,
$p_n\propto \dfrac{(c)_n\,(c+k)_n}{(k+1)_n\,n!}\,w^n$, with normalizer
$\binom{k+c-1}{k}\,{}_2F_1(c+k,c;k+1;w)$: a member of Kemp's family of
generalized hypergeometric probability distributions
\cite{KK56,Kemp68,JKK05}, with argument $w\in(0,1)$.
\end{proposition}

\begin{remark}[Placement]
The ratio $p_{n+1}/p_n=w\prod_i(n+\kappa_i+c)/\prod_i(n+\kappa_i+1)$
is a ratio of degree-$r$ polynomials, the defining property of
Kemp's generalized hypergeometric probability distributions
\cite{Kemp68,KK56}, \cite[\S2.4]{JKK05}; the law is therefore a
polynomial birth--death distribution in the sense of \cite{BX01,XZ08},
whose factor bounds apply at every $r$, and a discrete Gibbs measure
in the sense of \cite{ER08}, whose bounds need no normalizing
constant---here a generalized hypergeometric function.  Conditioning
independent negative binomials on their sum gives the
negative-hypergeometric (Dirichlet-multinomial) law; conditioning on
their differences, which gives the present family, we have not found
in the literature.
\end{remark}

\begin{lemma}[Balance, log-concavity, mode, moments]\label{lem:basic}
Let $N\sim\HNB(\bk,c,w)$ and let $d(n)=\prod_i(n+\kappa_i)$,
$b(n)=\prod_i(n+\kappa_i+c)$, the rate pair of the law; note $d(0)=0$.
\begin{enumerate}
\item[(i)] \emph{(Balance.)}\quad
$d(n+1)\,p_{n+1}\;=\;w\,b(n)\,p_n$ for all $n\ge0$; equivalently
$t_n:=p_{n+1}/p_n=w\,b(n)/d(n+1)$.
\item[(ii)] \emph{(Log-concavity.)}\quad
$t_n=w\prod_i\bigl(1+\tfrac{c-1}{\,n+1+\kappa_i}\bigr)$ is strictly
decreasing in $n$ iff $c\ge2$; for $c=1$, $t_n\equiv w$ and
$\HNB(\bk,1,w)=\Geom(1-w)$ for every $\bk$.
\item[(iii)] \emph{(Mode.)}\quad
Let $m^\ast=\max\{n\ge1: w\,b(n-1)\ge d(n)\}$ (and $m^\ast=0$ if the set
is empty).  Then $p_n$ increases up to $m^\ast$ and strictly decreases
after.
\item[(iv)] \emph{(Moment identities.)}\quad
$\E N = w\,\Zhat_c'(\bk,w)/\Zhat_c(\bk,w)$ and
\begin{equation}\label{eq:momid}
\E\,d(N)\;=\;w\,\E\,b(N).
\end{equation}
\end{enumerate}
\end{lemma}

\begin{remark}[Location of the mode]\label{rem:mode}
With $\bar\kappa=\max_i\kappa_i$ and $s=w^{-1/r}-1$, the mode of
Lemma~\ref{lem:basic}(iii) satisfies
$\bigl|m^\ast-(c-1)/s\bigr|\le\bar\kappa+1$ (proof of
Lemma~\ref{lem:basic}); the statement is not used below.
\end{remark}

\begin{remark}[The boundary $c=1$]\label{rem:c1}
For a single column the fiber is a single ladder $R_i=n+\kappa_i$ and
\eqref{eq:hnbpmf} collapses to $p_n\propto w^n$: the conditioned law is
exactly geometric, only weakly log-concave, and combinatorially trivial
(one-column tables are determined by their margins).  All strictness
statements below assume $c\ge2$; the case $c=1$ is covered by classical
geometric approximation \cite{BP99,PR11}.
\end{remark}

\begin{proposition}[Poisson limit: recovery of the hyper-Bessel theory]\label{prop:poissonlimit}
Fix $\bk$ and $\lambda_1,\dots,\lambda_r>0$, and let
$q_i=q_i(c)$ satisfy $cq_i/(1-q_i)\to\lambda_i$ as $c\to\infty$ (so
$w c^r\to\nu:=\prod_i\lambda_i$).  Then
\[
\dtv\bigl(\HNB(\bk,c,w),\ \HB(\bk,\nu)\bigr)\;\longrightarrow\;0,
\qquad
w\,b(n)\;\longrightarrow\;\nu \quad\text{pointwise},
\]
where $\HB(\bk,\nu)$ is the hyper-Bessel law of \cite{DeS26} with mass
proportional to $\nu^n/\prod_i(n+\kappa_i)!$.
\end{proposition}

\begin{remark}
In this limit the operator, the fixed-point transform, and the
main theorem degenerate to their counterparts in \cite{DeS26}; the
rate of convergence is examined in \S\ref{sec:computations}.
\end{remark}

\begin{corollary}[A Stein operator for the Bessel law]\label{cor:bessel}
At $r=2$, $\bk=(k,0)$, the hyper-Bessel law is the Bessel
distribution of Yuan and Kalbfleisch \cite{YK00}: two independent
Poisson counts of means $\nu_1,\nu_2$ conditioned on their difference
$k$, with pmf $\propto\lambda^{n}/(n!\,(n+k)!)$, $\lambda=\nu_1\nu_2$
(sampled in \cite{Dev02}).  An $\N$-valued $W$ with $\E W^{2}<\infty$
has this law if and only if
\[
\E\bigl[\lambda\,f(W+1)-W(W+k)\,f(W)\bigr]=0
\qquad\text{for all bounded }f,
\]
the limit of the HNB operator at $c\to\infty$, $wc^{2}\to\lambda$.
We have not found a Stein operator for the Bessel law in the
literature; its factors are Open Problem~\ref{op:bessel}.
\end{corollary}

\begin{proof}
The identification is Proposition~\ref{prop:poissonlimit} at $r=2$,
$\kappa=(k,0)$, whose limit mass $\nu^{n}/(n!\,(n{+}k)!)$ is the
Bessel pmf of \cite{YK00}.  The Bessel law $\pi$ has the balance
identity $\lambda p_n=(n{+}1)(n{+}k{+}1)p_{n+1}$; the rates
$b\equiv\lambda$, $d(n)=n(n{+}k)$ fall outside the standing
assumptions of \S\ref{sec:nb} (the ratio $d/b$ is unbounded), so we
argue directly.  Necessity: $\E[\lambda f(W{+}1)]=\sum_n\lambda
p_nf(n{+}1)=\sum_n(n{+}1)(n{+}k{+}1)p_{n+1}f(n{+}1)=\E[W(W{+}k)f(W)]$,
absolutely convergent since $\E W^{2}<\infty$.  Sufficiency: for
bounded $h$ the equation $\lambda g(n{+}1)-n(n{+}k)g(n)=h(n)-\pi(h)$
is solved by
\[
g(m{+}1)=\frac1{\lambda p_m}\sum_{n\le m}(h(n)-\pi(h))p_n,
\qquad
|g(m{+}1)|\le\frac{F(m)G(m{+}1)}{\lambda p_m}\le\frac{G(m{+}1)}{\lambda p_m},
\]
exactly as in Theorem~\ref{thm:solution}(ii), and the right side is
bounded because $G(m{+}1)/p_m\le p_{m+1}/(p_m(1-t_{m+1}))$ with
$t_n=\lambda/((n{+}1)(n{+}k{+}1))\to0$; hence $g$ is bounded, and
$0=\E(\mathcal Ag_A)(W)=\PP(W\in A)-\pi(A)$ for every $A\subseteq\N$.
\end{proof}

The fixed point and the main theorem now follow at these rates.

\begin{theorem}[Fixed-point characterization]\label{thm:fixedpoint}
$W\sim\HNB(\bk,c,w)$ for some $w\in(0,1)$ if and only if
\[
W^{d}\;\stackrel{d}{=}\;W^{b}+1,
\]
in which case necessarily $w=\beta_d/\beta_b$.  Equivalently, on laws with
$0<\beta_d\le\beta_b<\infty$ define the transform $T$ by
$\PP(TW=n)\propto d(n+1)\PP(W=n+1)/b(n)$; then the fixed points of $T$
are exactly $\{\HNB(\bk,c,w):w\in(0,1)\}$.
\end{theorem}

\begin{theorem}[Ratio-matched main theorem: the two rate moments and one coupling]\label{thm:master}
Let $W$ be $\N$-valued with $\beta_b<\infty$ and $\beta_d>0$, and let $w\in(0,1)$, with
$\Lambda,\Gamma$ the factors of $\HNB(\bk,c,w)$.  Then
\begin{equation}\label{eq:master}
\dtv\bigl(W,\ \HNB(\bk,c,w)\bigr)\;\le\;
\Lambda\;w\,\beta_b\;\dwass\bigl(W^{b}+1,\ W^{d}\bigr)
\;+\;\Gamma\;\bigl|\beta_d-w\,\beta_b\bigr| .
\end{equation}
If moreover $w^\ast:=\beta_d/\beta_b<1$, then with $w=w^\ast$ the second
term vanishes:
\begin{equation}\label{eq:mastermatched}
\dtv\bigl(W,\ \HNB(\bk,c,w^\ast)\bigr)\;\le\;
\Lambda_{w^\ast}\,\beta_d\;\dwass\bigl(W^{b}+1,\ W^{d}\bigr).
\end{equation}
\end{theorem}

The coupling that makes \eqref{eq:mastermatched} effective is built from
three elementary identities.

\begin{lemma}[Phantom column]\label{lem:phantom}
Let $R\sim\NB(c,q)$.  Then:
\begin{enumerate}
\item[(i)] the $(n+c)$-biased version of $R$ is $\NB(c{+}1,q)$;
\item[(ii)] $\NB(c{+}1,q)\stackrel{d}{=}\NB(c,q)+\Geom(1-q)$ with the
summands independent;
\item[(iii)] hence the $b$-type bias of a row adds one fresh geometric ---
a \emph{phantom column}.
\end{enumerate}
\end{lemma}

\begin{lemma}[Size bias of a negative binomial row]\label{lem:sizebias}
The size-biased version of $R\sim\NB(c,q)$ satisfies
$R^{\,s}\stackrel{d}{=}1+\NB(c{+}1,q)\stackrel{d}{=}R+1+\Geom(1-q)$,
the geometric independent of $R$.
\end{lemma}

Coupling the fresh geometrics in Lemmas~\ref{lem:phantom} and
\ref{lem:sizebias} \emph{identically}, the difference between the
$d$-bias and the $b$-bias of a row is a deterministic $+1$.  On the fiber
$\{R_i=n+\kappa_i\ \forall i\}$ one has $d(N)=\prod_iR_i$ and
$b(N)=\prod_i(R_i+c)$: the $d$-bias size-biases every row simultaneously
and the $b$-bias adds a phantom column to every row; both shift all rows
together, so both \emph{preserve all differences}, and their coupled
discrepancy is the fiber-preserving all-rows increment.  The
characterization underlying the vanishing lemma of
Section~\ref{sec:dependency} is:

\begin{lemma}[Geometric characterization by the shared phantom]\label{lem:geomchar}
Let $Z=(X^{(1)},\dots,X^{(r)})$ be $\N^r$-valued with finite means, and
fix $i$ with $\E X^{(i)}>0$.  Suppose the $x_i$-biased version of $Z$ has the same law as the
$(x_i{+}1)$-biased version of $Z$ shifted by $e_i$.  Then
$X^{(i)}\sim\Geom(1-q_i)$ with $q_i=\E X^{(i)}/\E(X^{(i)}{+}1)$, and
$X^{(i)}$ is independent of $(X^{(j)})_{j\ne i}$.  In particular, for a
single coordinate, the $(x{+}1)$-bias of $X$ equals
$X+X'$ (an independent copy) iff $X$ is geometric.
\end{lemma}

\subsection*{The instruments at general $r$}

The instruments are stated for a generic array of nonnegative integers with
$r$ rows and $c$ columns whose entries $X^{(i)}_j$ are independent
within each column.  Write
$U_i=\sum_{j=1}^cX^{(i)}_j$ for the row sums,
$D_i=U_i-U_{i+1}$ ($1\le i\le r-1$) for their successive differences,
$\bk=(k_1,\dots,k_{r-1})$ for a value of $D=(D_1,\dots,D_{r-1})$, and
$\kappa_i=k_i+\dots+k_{r-1}$, $\kappa_r=0$, for the offsets it induces;
$[r]=\{1,\dots,r\}$, and $S^c$ is the
complement of $S$ in $[r]$.  Conditioning on $D=\bk$ pins $U_i=U_r+\kappa_i$ and leaves the
single free height $U_r$.

Approximation is measured through two size-biasings, defined for any
birth--death rate pair $b,d$.

The constants that multiply these distances come from the Stein
equation of the rate pair.

\begin{definition}[Stein operator and factors]\label{def:factors}
For a birth--death rate pair $(b,d)$ with $d(0)=0$ and a parameter
$w\in(0,1)$, the \emph{target law} $\pi$ of the pair is the law with
successive ratios $\pi_{n+1}/\pi_n=w\,b(n)/d(n+1)$---explicitly,
\[
\pi_n\;=\;\pi_0\prod_{j=0}^{n-1}\frac{w\,b(j)}{d(j+1)},
\qquad
\pi_0^{-1}\;=\;\sum_{n\ge0}\,\prod_{j=0}^{n-1}\frac{w\,b(j)}{d(j+1)},
\]
assumed summable.  For the pair $b(n)=n+c$, $d(n)=n$ of the first rung
the product telescopes to the \emph{negative binomial},
\[
\pi_n\;=\;\binom{n+c-1}{n}(1-w)^c\,w^n,
\qquad\text{that is,}\qquad \pi=\NB(c,w);
\]
for the pair $b(n)=\prod_{i\le r}(n+\kappa_i+c)$,
$d(n)=\prod_{i\le r}(n+\kappa_i)$ it is the \emph{hyper negative
binomial},
\[
\pi_n\;\propto\;\prod_{i=1}^{r}\binom{n+\kappa_i+c-1}{n+\kappa_i}\,w^n
\;=:\;\HNB(\bk,c,w)_n,
\]
the law of the conditioned minimal row (proved as
Theorem~\ref{thm:target}); the interpolated pair of
Theorem~\ref{thm:system} is the third instance.  The \emph{Stein
equation} for a set $A\subseteq\N$ is
\[
w\,b(n)\,g_A(n+1)-d(n)\,g_A(n)\;=\;\ind{n\in A}-\pi(A),
\qquad n\ge0,
\]
with the convention $g_A(0)=0$.  It has a unique bounded solution,
and the \emph{Stein factors} of the pair are its uniform size and
increment,
\[
\Gamma\;=\;\sup_{A\subseteq\N}\ \sup_{m\ge1}\,|g_A(m)|,
\qquad
\Lambda\;=\;\sup_{A\subseteq\N}\ \sup_{m\ge0}\,|g_A(m+1)-g_A(m)|,
\]
written $\Gamma_w,\Lambda_w$ when the parameter matters.
\end{definition}

Existence and uniqueness of the bounded solution are
Theorem~\ref{thm:solution}, rate-generic; the factors are finite and
explicit.

\begin{proposition}[The factors, explicitly]\label{prop:factorbound}
Let $p_n$, $F$, $G$ be the mass, distribution and tail functions of
the target $\pi$.  The size factor is \emph{exact},
\[
\Gamma\;=\;\sup_{m\ge0}\frac{F(m)\,G(m+1)}{d(m+1)\,p_{m+1}},
\qquad\text{and}\qquad
\Lambda\;\le\;\min\Bigl\{2\Gamma,\;
\max\Bigl(1,\tfrac{1-w}{w}\Bigr)\Gamma+\frac1{w\,b(0)}\Bigr\}:
\]
one pass over the target's own pmf, with no reference to $W$.  When
$t_n=w\,b(n)/d(n{+}1)$ is strictly decreasing, for the
hyper-negative-binomial law when $c\ge2$ (Lemma~\ref{lem:basic}(ii)),
a form closed in the rates alone,
\[
\Gamma\;\le\;\max\Bigl\{\frac1{w\,b(0)},\;
\sup_{m\ge1:\,x_m>0}\frac1{x_m},\;
\sup_{m\ge0:\,y_{m+1}>0}\frac{d(m{+}2)}{d(m{+}1)\,y_{m+1}}\Bigr\},
\quad
\begin{aligned}
x_m&=w\,b(m{-}1)-d(m),\\
y_m&=d(m{+}1)-w\,b(m),
\end{aligned}
\]
has finitely many terms in the second maximum, and the third is a
maximum since $d(m)\to\infty$.  The exact form is the one to evaluate
(\S\ref{sec:computations}).
\end{proposition}

The first instrument prices the distance to the hyper-negative-binomial
target as a sum of per-column costs.

\begin{definition}[Moments and the fiber constant]\label{def:fiberconst}
For $T\subseteq[r]$ write $\nu^b_T=\E\prod_{i\in T}(U_i+c)$, the
unconditioned row-sum moments---closed for independent rows,
$\nu^b_T=\prod_{i\in T}(c+\E U_i)$.  For a set $L\subseteq[c]$ of at
most $r$ columns let $D^{(L)}$ be the difference vector of the reduced
system with the columns of $L$ removed, and let
$\mathcal{W}\subseteq\Zint^{r-1}$ be the set of lattice points lying
coordinatewise between two points $\bk-\partial z$ and
$\bk-\partial z'$, where $z$ and $z'$ range over the support of the
sum of at most $r$ columns, each column either as it is or biased by
$\prod_{i\in S}x_i$ or by $\prod_{i\in S}(x_i+1)$ for some
$\varnothing\ne S\subseteq[r]$, and $\partial$ is the difference map
$\partial z=(z_1-z_2,\dots,z_{r-1}-z_r)$.  The \emph{fiber
constant} is
\[
\theta\;=\;\max_{|L|\le r}\ \sup_{s\in\mathcal{W}}\
\frac{\PP\bigl(D^{(L)}=s\bigr)}{\PP\bigl(D=\bk\bigr)},
\]
computable by dynamic programming on the reduced systems; for bounded
column supports the window is finite and explicit.
\end{definition}

\begin{theorem}[Dependency expansion across columns, product form]\label{thm:colexpansion}
Let the column coordinates be independent,
$X^{(i)}_j\sim\Geom(1-y_{ij})$ with $y_{ij}\in(0,1)$ arbitrary, fix a
reference $q_i\in(0,1)$, $w=\prod_iq_i$, $w_S=\prod_{i\in S}q_i$, and
let $W=\mathcal{L}(U_r\mid D=\bk)$.  Then, with $\Gamma$
the Stein factor of $\HNB(\bk,c,w)$ and $\theta$ the fiber constant,
\[
\dtv\bigl(W,\HNB(\bk,c,w)\bigr)\;\le\;\Gamma\,\theta
\sum_{\emptyset\ne S\subseteq[r]}\nu^b_{S^c}
\sum_{j=1}^c\ \prod_{i\in S}\frac{1}{1-y_{ij}}\
\Bigl|\prod_{i\in S}y_{ij}-w_S\Bigr|.
\]
In the Barvinok calibration $y_{ij}=q_i\beta_j$ the mismatch is
$w_S\,|\beta_j^{|S|}-1|$.  Capped entries and arbitrary within-column
laws are covered by Theorem~\ref{thm:expansion}, with the coupling
costs of Corollary~\ref{cor:capped}.
\end{theorem}

The column direction alone does not price the target count itself;
for that, condition in both directions at once.

Give the table independent entries $X_{ij}\sim\Geom(1-x_iy_j)$
in product form, with row sums $U_i$ and column sums $V_j$, and
condition on both systems of \emph{differences}: the row differences
$D_{\mathrm r}$ and the column differences $D_{\mathrm c}$, with
offsets taken from the target margins themselves,
$\kappa_i=r^\ast_i-r^\ast_r$ and $\lambda_j=c^\ast_j-c^\ast_c$.
Every margin is then an affine function of the single free variable
$W=U_r$, which lives on an arithmetic progression $L$ of step
$h=c/g$, $g=\gcd(r,c)$, least element $k_0$; with the ray weight
$\Omega=(\prod_ix_i)^{c/g}(\prod_jy_j)^{r/g}$, the conditional law
$P_\bullet$ of $W$ has pmf proportional to
$|\T(\br(k),\bc(k))|\,\Omega^{(k-k_0)/h}$ along the ray of margins
through the target, and the target itself is the atom $k^\ast$ of $L$
with $\br(k^\ast)=\br^\ast$, $\bc(k^\ast)=\bc^\ast$.  Each
conditioning alone characterizes $W$ by a birth--death equation on
$L$.  The row conditioning gives, when the columns are identically
distributed ($y_j\equiv y$, Theorem~\ref{thm:target}), the
hyper-negative-binomial law with the coarse ratio
\[
\rho_1(k)=\prod_{s=0}^{h-1}
\frac{\prod_i\,(k+s+\kappa_i+c)}{\prod_i\,(k+s+1+\kappa_i)},
\]
and for general $y_j$ this operator is a reference, exact in
expectation only at $y_j\equiv y$.  The column conditioning gives
the law $\pi_2$ of Proposition~\ref{prop:rowcol}, characterized by
its canonical birth--death equation with ratio
$\rho_2(k)=\pi_2(k{+}h)/\pi_2(k)$ on $L$, which is in closed form
when the rows are identically distributed ($x_i\equiv x$): with the
column height $w'(k)=\bigl(rk+\sum_i\kappa_i-\sum_j\lambda_j\bigr)/c$,
\[
\rho_2(k)=\prod_{s=0}^{r/g-1}
\frac{\prod_j\,(w'(k)+s+\lambda_j+r)}{\prod_j\,(w'(k)+s+1+\lambda_j)}
\qquad(x_i\equiv x).
\]
The interpolated family $\pi_{w,\mu}$ on $L$ has
$\pi_{w,\mu}(k{+}h)/\pi_{w,\mu}(k)=w\,\rho_1(k)^{1-\mu}\rho_2(k)^{\mu}$,
realized by the rate pair $\tilde b=B_1^{1-\mu}B_2^{\mu}$,
$\tilde d=D_1^{1-\mu}D_2^{\mu}$ with
$B_1(k)=\prod_{s<h}\prod_i(k+s+\kappa_i+c)$,
$D_1(k{+}h)=\prod_{s<h}\prod_i(k+s+1+\kappa_i)$, $B_2\equiv1$ and
$D_2(k{+}h)=1/\rho_2(k)$.

\subsection*{The row-wise convolution, for the two-Stein method}

The two-directional matching of \S\ref{sec:rowcol} governs its
second direction by convolutions of row laws across the rows.

\begin{definition}[Row-wise convolution]\label{def:conv}
For independent $\N$-valued $X_1,\dots,X_r$, the row-wise convolution
is $V=\sum_{i=1}^rX_i$; the target convolution is
$V^\ast=\sum_{i=1}^rW_i$ with independent
$W_i\sim\HNB(\bk^{(i)},c,w_i)$.
\end{definition}

\begin{lemma}[Moments of the target convolution]\label{lem:convmom}
$\E\,W_i=w_i\,\partial_{w}\log\Zhat_c(\bk^{(i)},w)\big|_{w=w_i}$ and
$\mathrm{Var}\,W_i=\bigl(w\,\partial_w\bigr)^2\log\Zhat_c(\bk^{(i)},w)
\big|_{w=w_i}$, so the mean and variance of $V^\ast$ are the explicit
sums of these log-derivatives.
\end{lemma}

\begin{theorem}[The rate moments of the convolution]\label{thm:aggconv}
In the spirit of Arratia--Goldstein--Gordon \cite{AGG89}: let
$X_1,\dots,X_r$ be independent $\N$-valued random variables,
$V=\sum_iX_i$, and for parameters $\bk^{(i)}$, $w_i\in(0,1)$ let
$N_i\sim\HNB(\bk^{(i)},c,w_i)$ be independent, $V^\ast=\sum_iN_i$.
With $X_i^{\,b},X_i^{\,d}$, $\beta_{b,i},\beta_{d,i}$ the biasings and
rate moments of $X_i$ at the rates of $\HNB(\bk^{(i)},c,w_i)$, assumed
to satisfy $\beta_{b,i}<\infty$ and $\beta_{d,i}>0$, and
$\Gamma_i,\Lambda_i$ that law's factors,
\[
\dtv\bigl(V,\,V^\ast\bigr)\;\le\;\sum_{i=1}^r\Bigl[
\Lambda_i\,w_i\,\beta_{b,i}\,
\dwass\bigl(X_i^{\,b}+1,\ X_i^{\,d}\bigr)
+\Gamma_i\,\bigl|\beta_{d,i}-w_i\,\beta_{b,i}\bigr|\Bigr],
\]
each bracket the bound of Theorem~\ref{thm:master} for $X_i$ against
$\HNB(\bk^{(i)},c,w_i)$: local rate-moment errors add across
independent components.
\end{theorem}

Rows exact by Theorem~\ref{thm:target}, the convolution across rows
within its rate moments by Theorem~\ref{thm:aggconv}: these are the two
systems that \S\ref{sec:rowcol} matches simultaneously.

\subsection*{Dependency expansion for column arrays}

Let $Z_j=(X^{(1)}_j,\dots,X^{(r)}_j)$, $1\le j\le c$, be independent
$\N^r$-valued columns with arbitrary within-column joint laws and
$\E\|Z_j\|_1^{\,r+1}<\infty$.  Set
$U_i=\sum_jX^{(i)}_j$, $D=(U_1-U_2,\dots,U_{r-1}-U_r)$, fix
$\bkk$ with $\PP(D=\bkk)>0$, re-base so that $\kappa_r=0\le\kappa_i$, and
let
\[
W\;=\;\mathcal{L}\bigl(U_r\,\bigm|\,D=\bkk\bigr).
\]

\begin{definition}[Tilts, costs, aggregates, fiber constants]\label{def:costs}
Fix reference parameters $q_1,\dots,q_r\in(0,1)$ and $w=\prod_iq_i$,
$w_S=\prod_{i\in S}q_i$.  For $\emptyset\ne S\subseteq[r]$ and a column
$Z_j$:
\begin{enumerate}
\item[(a)] the \emph{$d$-tilt} $Z_j^{S,x}$ has law
$\propto\bigl(\prod_{i\in S}x_i\bigr)\PP(Z_j=x)$ with normalizer
$\chi^{S,x}_j:=\E\prod_{i\in S}X^{(i)}_j$; the \emph{$b$-tilt}
$Z_j^{S,x+1}$ has law $\propto\prod_{i\in S}(x_i{+}1)\,\PP(Z_j=x)$ with
normalizer $\chi^{S,x+1}_j:=\E\prod_{i\in S}(X^{(i)}_j{+}1)$;
\item[(b)] the \emph{coupling cost} and \emph{weight cost} are
\[
\eta^{(S)}_j:=\dwass\bigl(Z_j^{S,x},\ Z_j^{S,x+1}+e_S\bigr),
\qquad
\delta^{(S)}_j:=\bigl|\chi^{S,x}_j-w_S\,\chi^{S,x+1}_j\bigr|,
\]
with $e_S=\sum_{i\in S}e_i$ and $\dwass$ the $\ell_1$-Wasserstein
distance on $\Zint^r$;
\item[(c)] the aggregates and moments are
$H^{(S)}=\sum_j\chi^{S,x+1}_j\eta^{(S)}_j$,\ \
$\Delta^{(S)}=\sum_j\delta^{(S)}_j$,\ \
$\nu^b_T=\E\prod_{i\in T}(U_i+c)$;
\item[(d)] for a set $L\subseteq[c]$ of at most $r$ columns let
$D_R^{(L)}$, $V_R^{(L)}$ denote the difference vector and $r$-th row sum
of the reduced system $(Z_j)_{j\notin L}$, let
$\mathcal{W}\subseteq\Zint^{r-1}$ be the window of
Definition~\ref{def:fiberconst}, and let the suprema below run over
the points of $\mathcal{W}$ charged by $D_R^{(L)}$.  Define
\[
\theta=\max_{|L|\le r}\sup_{s\in\mathcal{W}}
\frac{\PP(D_R^{(L)}=s)}{\PP(D=\bkk)},\quad
\tau_1=\max_{|L|\le r}\sup_{\substack{s,s'\in\mathcal{W}\\ \|s-s'\|_1=1}}
\frac{|\PP(D_R^{(L)}=s)-\PP(D_R^{(L)}=s')|}{\PP(D=\bkk)},
\]
\[
\omega=\max_{|L|\le r}\sup_{\substack{s,s'\in\mathcal{W}\\ \|s-s'\|_1=1}}
\dwass\Bigl(\mathcal{L}\bigl(V_R^{(L)}\bigm|D_R^{(L)}=s\bigr),\
\mathcal{L}\bigl(V_R^{(L)}\bigm|D_R^{(L)}=s'\bigr)\Bigr).
\]
These are the exact analogues of the fiber constants of \cite{DeS26},
computable by dynamic programming on the reduced systems; for bounded
column supports the window is finite and explicit.
\end{enumerate}
\end{definition}

\begin{lemma}[Transfer identities]\label{lem:transfer}
Let the columns $Z_1,\dots,Z_c$ be independent with
$\E\|Z_j\|_1^{\,r}<\infty$.  For every
$\vec\jmath=(\jmath_1,\dots,\jmath_r)\in[c]^r$ let
$L(\vec\jmath)=\{\jmath_1,\dots,\jmath_r\}$ and let
$S_j=\{i:\jmath_i=j\}$ for $j\in L(\vec\jmath)$; the nonempty blocks
$S_j$ partition $[r]$.  Then $\prod_iU_i=\sum_{\vec\jmath}\prod_i
X^{(i)}_{\jmath_i}$ and, since $U_i+c=\sum_{j=1}^c(X^{(i)}_j+1)$, also
$\prod_i(U_i+c)=\sum_{\vec\jmath}\prod_i(X^{(i)}_{\jmath_i}+1)$.
Consequently, for bounded $g$ with $g(0)=0$,
\begin{align}
\E\Bigl[\prod_iU_i\ g(U_r)\ind{D=\bkk}\Bigr]
&=\sum_{\vec\jmath\in[c]^r}\Bigl(\prod_{j\in L(\vec\jmath)}\chi^{S_j,x}_j\Bigr)\,
\E\bigl[g\bigl(V^{\vec\jmath,x}\bigr)\ind{D^{\vec\jmath,x}=\bkk}\bigr],
\label{eq:transferd}\\
\E\Bigl[\prod_i(U_i{+}c)\ g(U_r{+}1)\ind{D=\bkk}\Bigr]
&=\sum_{\vec\jmath\in[c]^r}\Bigl(\prod_{j\in L(\vec\jmath)}\chi^{S_j,x+1}_j\Bigr)\,
\E\bigl[g\bigl(V^{\vec\jmath,x+1}{+}1\bigr)\ind{D^{\vec\jmath,x+1}=\bkk}\bigr],
\label{eq:transferb}
\end{align}
where the superscripted systems replace each hit column $Z_j$,
$j\in L(\vec\jmath)$, by its $S_j$-tilt, and $V$, $D$ denote the $r$-th
row sum and difference vector of the modified system.
\end{lemma}

\begin{lemma}[Vanishing]\label{lem:vanishing}
$\eta^{(S)}_j=\delta^{(S)}_j=0$ for all $j$ and all
$\emptyset\ne S\subseteq[r]$ if and only if every column $Z_j$ has
independent coordinates $X^{(i)}_j\sim\Geom(1-q_i)$.  In that case
$W=\HNB(\bk,c,w)$ exactly.
\end{lemma}

\begin{theorem}[Dependency expansion]\label{thm:expansion}
Let $Z_1,\dots,Z_c$ be independent $\N^r$-valued columns with
$\E\|Z_j\|_1^{\,r+1}<\infty$, so that every tilted column has a
finite mean and every $H^{(S)}$ is finite, let $\bkk$ be a
difference vector with
$\PP(D=\bkk)>0$, and let $W=\mathcal L(U_r\mid D=\bkk)$.  With the
notation of Definition~\ref{def:costs}, for any reference
$q_1,\dots,q_r\in(0,1)$ and $w=\prod_iq_i$, and with $\Lambda,\Gamma$ the
Stein factors of $\HNB(\bk,c,w)$,
\begin{equation}\label{eq:expansion}
\dtv\bigl(W,\ \HNB(\bk,c,w)\bigr)\;\le\;
\bigl[\theta(1+2\omega)\Lambda+2\tau_1\Gamma\bigr]\;
w\!\!\sum_{\emptyset\ne S\subseteq[r]}\!\!\nu^{b}_{S^c}\,H^{(S)}
\;+\;\Gamma\,\theta\!\!\sum_{\emptyset\ne S\subseteq[r]}\!\!
\nu^{b}_{S^c}\,\Delta^{(S)} .
\end{equation}
\end{theorem}

\begin{corollary}[Product-form heterogeneous geometric columns]\label{cor:productform}
Suppose the coordinates of each column are independent with
$X^{(i)}_j\sim\Geom(1-y_{ij})$, $y_{ij}\in(0,1)$ arbitrary.  Then
$\eta^{(S)}_j=0$ for \emph{all} $j,S$: the shared-phantom coupling is
exact column by column, regardless of parameters.  Moreover
$\chi^{S,x}_j=\prod_{i\in S}\tfrac{y_{ij}}{1-y_{ij}}$ and
$\chi^{S,x+1}_j=\prod_{i\in S}\tfrac{1}{1-y_{ij}}$, so
\eqref{eq:expansion} reduces to
\begin{equation}\label{eq:productform}
\dtv\bigl(W,\HNB(\bk,c,w)\bigr)\;\le\;\Gamma\,\theta
\sum_{\emptyset\ne S\subseteq[r]}\nu^b_{S^c}
\sum_{j=1}^c\ \prod_{i\in S}\frac{1}{1-y_{ij}}\
\Bigl|\prod_{i\in S}y_{ij}-w_S\Bigr| .
\end{equation}
In the Barvinok calibration $y_{ij}=q_i\beta_j$ the mismatch is
$w_S\,|\beta_j^{|S|}-1|$: only the scalar column weights
$\beta_j$ enter, through $\sum_j|\beta_j^{|S|}-1|$ --- a finite-$n$,
certified form of the independence heuristic \cite{Barv09,BH12} at the
level of row differences.
\end{corollary}

\begin{corollary}[Capped geometric columns; binary entries]\label{cor:capped}
Suppose the coordinates of each column are independent capped geometrics,
$X^{(i)}_j\sim\Geom(1-y_{ij})$ conditioned on $X^{(i)}_j\le L_{ij}$.
Then, with reference $q_i$ and singleton $S=\{i\}$, when
$y_{ij}=y_i$ for every $j$ and $q_i=y_i$,
\begin{equation}\label{eq:deltacap}
\delta^{(\{i\})}_j\;=\;\frac{(L_{ij}+1)\,(1-y_{ij})\,y_{ij}^{\,L_{ij}+1}}
{1-y_{ij}^{\,L_{ij}+1}}
\qquad(\text{exact closed form}),
\end{equation}
in particular $\delta=2y^2/(1+y)$ for binary entries $L=1$; and
$\eta^{(\{i\})}_j\le\min\{L_{ij},\mathcal{E}(y_{ij},L_{ij})\}$,
$\mathcal{E}$ an explicit finite geometric-tail expression displayed
in the proof, of order $O\bigl((L+2)^2y^{L+1}\bigr)$, and $L_{ij}$
the trivial bound, both laws being supported on
$\{1,\dots,L_{ij}+1\}$.  Multi-coordinate costs obey
$\eta^{(S)}_j\le\sum_{i\in S}\mathcal{E}(y_{ij},L_{ij})$ and a
telescoped product bound for $\delta^{(S)}_j$ (displayed in the proof).
Constants here are explicit but not optimized.
\end{corollary}

The register, assembling the expansion into a usable certificate:

\begin{theorem}[TV register]\label{thm:tvregister}
Let the entries be independent capped geometrics with parameters
$y_{ij}$ and caps $L_{ij}$, and let $W$ be the conditioned minimal row
sum as in Section~\ref{sec:dependency}.  Then for any reference $\bq$,
$\dtv(W,\HNB(\bk,c,w))$ is bounded by \eqref{eq:expansion} with the
costs of Corollaries~\ref{cor:productform}--\ref{cor:capped}; in the
product-form uncapped case the bound is \eqref{eq:productform}, and in the
homogeneous capped case every cost is an explicit function of
$(y_i,L_{ij})$ given by \eqref{eq:deltacap} and $\mathcal{E}(\cdot)$.
\end{theorem}

\begin{proof}
Theorem~\ref{thm:expansion} applies to any independent columns with
finite $r$-th moments, capped geometrics included; the costs are
evaluated in Corollaries~\ref{cor:productform}--\ref{cor:capped},
whose displays are substituted into \eqref{eq:expansion}.
\end{proof}

\subsection*{Applications at general $r$}

The one-row examples extend row by row; the new ones use the
conditioning itself.

\emph{Competing negative-binomial clocks.}

Clock $i$ advances through $c$ geometric stages of parameter $q_i$,
so its total is $\NB(c,q_i)$; observe the $r$ clocks only through
their gaps $D_i=R_i-R_{i+1}=k_i$ and ask for the common level.
Theorem~\ref{thm:target}: the level is $\HNB(\bk,c,w)$,
$w=\prod_iq_i$, and the gap configuration has probability
\eqref{eq:fiber}.  At $r=2$, $c=1$ this is the classical fact that
the smaller of two geometric clocks, given their difference, is
geometric with parameter $q_1q_2$.  Stage-dependent parameters
$y_{ij}$ break exactness, and Theorem~\ref{thm:colexpansion} prices
them, stage by stage through $|\beta_j^{|S|}-1|$.

\emph{Multi-species occupancy with conserved gaps.}

Put $r$ species of indistinguishable balls into $c$ distinguishable
boxes, each box filling independently and geometrically at fugacity
$q_i$ for species $i$, so that a species' unconditioned total is
$\NB(c,q_i)$ and conditioning on it gives Bose--Einstein occupancy.
Fix instead the gaps $N_i-N_{i+1}=k_i$ and let the level $n$
fluctuate: $\PP(n)\propto w^{n}\prod_i\binom{n+\kappa_i+c-1}{c-1}$,
the HNB law---an ensemble canonical in the species differences and
grand-canonical in the level, at fugacity $w=\prod_iq_i$, exact by
Theorem~\ref{thm:target}.  Inhomogeneous boxes---box $j$ at energy
$\epsilon_j$, weight $y_{ij}=e^{-\beta(\epsilon_j+\mu_i)}$---are
column heterogeneity, and Theorem~\ref{thm:colexpansion} certifies
the flat-spectrum HNB with the spectrum's spread as the mismatch:
sharp for sparse occupancy, of order one ball per box, and
moment-priced when dense.

\emph{Counts observed through their differences.}  Sites, cohorts, or
arms of a trial each accumulate a negative-binomial count of shape
$c$---a Poisson count with a Gamma rate---and what is recorded is the
comparison: how far each is ahead of the next.  The remaining
uncertainty is the level, and its law is $\HNB(\bk,c,w)$ with $w$ the
product of the rates' odds; Theorem~\ref{thm:master} bounds the
level's distance from that law when the shapes are only approximately
common.

\emph{The Bessel connection.}  In the Poisson regime of
Proposition~\ref{prop:poissonlimit} the level's law is hyper-Bessel,
and at $r=2$ this is the Bessel distribution \cite{YK00} of two
independent Poisson counts conditioned on their difference---the law
behind the Skellam test of two rates \cite{Skellam46}---with the
operator of Corollary~\ref{cor:bessel}.  The HNB is its negative-binomial
deformation, and the proposition is the certificate that the
deformation is continuous.

\emph{Contingency tables: the conditional law of a column.}  A
uniform table with row sums $\br$ is a geometric array with
exchangeable columns conditioned on its row totals---the product-form
representation behind the typical table and local limit theorem of
\cite{Barv09,BH12}, which describe single cells and Gaussian
fluctuations but not a column's law with an explicit error; read in the
column direction, with the row-sum profile fixed up to its level,
the level's law is exactly the HNB of Theorem~\ref{thm:target}, and
the column direction inherits the whole of this section: exact
conditionals at exchangeability, the expansion when the columns are
unequal.  What this section does not do is condition on the column
sums as well.  That is the two-directional problem of
\S\ref{sec:rowcol}.

\section{The row--column system}\label{sec:rowcol}\label{sec:apps}

Conditioning on both margins couples the two directions, and neither
single-margin conditioning sees the other.  The row--column system
of this section matches the two through a shared cell, and its
matched member carries information both single-margin conditionings
lose; the summary theorem states what the composed certificate
delivers, and the closing subsection states where, for contingency
tables, that is nontrivial: a small slice of the parameter space.

\subsection*{Definitions, matching, and the main theorem}

The main theorem takes one conditioning.  A table carries two: fix
the \emph{row} differences and the last row sum is characterized by
one birth--death equation (the HNB of \S\ref{sec:master}); fix the
\emph{column} differences and the same variable is characterized by a
second, the across-columns convolution of column coordinates.  Under
the joint conditioning---both difference systems at once---the law of
the free variable is the fiber's one-dimensional shadow: its pmf is
the table-count sequence along the margin ray, the closest
one-dimensional approach to the uniform contingency-table
distribution itself.  The standard method consumes one equation; this
section extends it to the pair: an interpolating operator family, a
main theorem uniform over it, and matching conditions under which
both equations hold in expectation at once.

\begin{definition}[Stein system, contrast, interpolated operators]\label{def:system}
Let $\mathcal A_if(k)=w\,b_i(k)f(k+1)-d_i(k)f(k)$, $i=1,2$, be
birth--death Stein operators with rate ratios
$\rho_i(k)=b_i(k)/d_i(k+1)$, and assume
$w\rho_i(k)\le\bar\rho<1$ for all $k\ge k_0$ (both targets exist and
are summable).  The \emph{contrast} is
\[
S(k)\;=\;\sum_{j=0}^{k-1}\log\frac{\rho_2(j)}{\rho_1(j)},
\qquad S(0)=0,
\]
assumed not affine in $k$, i.e.\ $\rho_2/\rho_1$ nonconstant, so that
$k$ and $S(k)$ are affinely independent.  For $\mu\in[0,1]$ the \emph{interpolated
operator} is
$\mathcal A_{w,\mu}f(k)=w\,\tilde b_\mu(k)f(k+1)-\tilde d_\mu(k)f(k)$
with $\tilde b_\mu=b_1^{1-\mu}b_2^{\mu}$,
$\tilde d_\mu=d_1^{1-\mu}d_2^{\mu}$, and $\pi_{w,\mu}$ its stationary
law; $W^{\tilde b},W^{\tilde d},\beta_{\tilde b},\beta_{\tilde d}$
denote Definition~\ref{def:bias} applied with the rates
$(\tilde b_\mu,\tilde d_\mu)$.
\end{definition}

\begin{lemma}[Identification: an exponential family in $(k,S)$]\label{lem:sysid}
$\pi_{w,\mu}(k)\propto\pi_{w,0}(k)\,e^{\mu S(k)}$: the system targets
form a two-parameter exponential family with natural parameters
$(\log w,\mu)$ and sufficient statistics $(k,S(k))$.
\end{lemma}

\begin{lemma}[Log-concavity of the system targets]\label{lem:syslc}
If both rate pairs $(b_i,d_i)$ satisfy the standing assumptions of
\S\ref{sec:nb} and $\rho_1$, $\rho_2$ are nonincreasing, then
$(\tilde b_\mu,\tilde d_\mu)$ satisfies them too,
$\rho_1^{1-\mu}\rho_2^{\mu}$ is nonincreasing, and every $\pi_{w,\mu}$
is log-concave;
the Stein factors $\tilde\Gamma_{w,\mu},\tilde\Lambda_{w,\mu}$ of
$\mathcal A_{w,\mu}$ are then finite and satisfy
Theorem~\ref{thm:factors}(i) and (iii) with
$(\tilde b_\mu,\tilde d_\mu)$ in place of $(b,d)$, whose formulas use
only the birth--death form, and (ii) when
$w\rho_1^{1-\mu}\rho_2^{\mu}$ is strictly decreasing.
\end{lemma}

\begin{theorem}[Stein-system main theorem]\label{thm:sysmaster}
Let $(b_i,d_i)$, $i=1,2$, be as in Definition~\ref{def:system} with
$d_i(0)=0$ and $d_i(n)>0$ for $n\ge1$, and let $w\in(0,1]$ and $\mu\in[0,1]$ be
such that $w\,\rho_1(k)^{1-\mu}\rho_2(k)^{\mu}\le\bar\rho<1$ for all
$k\ge k_0$ and $\inf_{n\ge1}\tilde d_\mu(n)>0$.  Let $p$, $F$, $G$
be those of $\pi_{w,\mu}$, and let
\[
\tilde\Gamma_{w,\mu}=\sup_{m\ge0}\frac{F(m)\,G(m+1)}{\tilde d_\mu(m+1)\,p_{m+1}}<\infty,
\qquad
\tilde\Lambda_{w,\mu}=\sup_{A\subseteq\N}\sup_{m\ge0}|g_A(m+1)-g_A(m)|\le2\tilde\Gamma_{w,\mu},
\]
$g_A$ the solution of $w\tilde b_\mu(n)g_A(n+1)-\tilde d_\mu(n)g_A(n)=\ind{n\in A}-\pi_{w,\mu}(A)$,
$g_A(0)=0$.  Let $W$ be $\N$-valued with $\beta_{\tilde b}<\infty$
and $\beta_{\tilde d}>0$.  Then
\[
\dtv\bigl(W,\ \pi_{w,\mu}\bigr)\;\le\;
\tilde\Lambda_{w,\mu}\;w\,\beta_{\tilde b}\;
\dwass\bigl(W^{\tilde b}+1,\ W^{\tilde d}\bigr)
\;+\;\tilde\Gamma_{w,\mu}\;
\bigl|\beta_{\tilde d}-w\,\beta_{\tilde b}\bigr| .
\]
\end{theorem}

\begin{theorem}[Simultaneous matching]\label{thm:sysmatch}
(i)~\emph{Ratio match.}  For every $\mu$, taking
$w^\ast(\mu)=\beta_{\tilde d}/\beta_{\tilde b}$ (assumed $<1$)
annihilates the second term of Theorem~\ref{thm:sysmaster}.
(ii)~\emph{Contrast match.}  Assume
$\E\bigl[(b_1{+}b_2)(W)\,(1{+}|S(W{+}1)|)\bigr]<\infty$ and
$\E\bigl[(d_1{+}d_2)(W)\,(1{+}|S(W)|)\bigr]<\infty$, and let
$\zeta(\mu)=\E\,S\bigl(W^{\tilde d}\bigr)-\E\,S\bigl(W^{\tilde b}{+}1
\bigr)$, $W^{\tilde b},W^{\tilde d}$ the biasings at the rates
$(\tilde b_\mu,\tilde d_\mu)$.  If $\zeta(0)$ and $\zeta(1)$ have
opposite signs, some $\mu^\ast\in[0,1]$ has $\zeta(\mu^\ast)=0$, and
at $(w^\ast(\mu^\ast),\mu^\ast)$ the defect
$\E(\mathcal A_{w,\mu}g)(W)=w\beta_{\tilde b}\E g(W^{\tilde b}{+}1)
-\beta_{\tilde d}\E g(W^{\tilde d})$ vanishes for both $g\equiv1$
(the rate match) and $g=S$ (the contrast match).
(iii)~\emph{I-projection.}  If $(\E W,\E S(W))=(\E_\pi k,\E_\pi S)$
for some $\pi=\pi_{w^\circ,\mu^\circ}$ with $w^\circ\in(0,1)$,
$\mu^\circ\in[0,1]$, then the pair $(w^\circ,\mu^\circ)$ is unique,
$\log Z(\log w,\mu)$ being strictly convex because $k$ and $S(k)$ are
affinely independent, and $\pi_{w^\circ,\mu^\circ}$ is the
information projection of $\mathcal L(W)$ onto the system targets.
\end{theorem}

\begin{remark}[Two row systems]\label{rem:tworow}
The construction already says something with two \emph{row}-difference
systems: $\mathcal A_1$ the $\HNB(\bk,c,w)$ operator and $\mathcal A_2$
the $\HNB(\bk',c,w)$ operator of a second difference statistic on the
same rows.  There the contrast telescopes to a closed form,
\[
S(k)=\sum_{i=1}^{r}\Bigl[
\log\tfrac{\Gamma(k+\kappa_i'+c)\,\Gamma(\kappa_i'+1)}
{\Gamma(\kappa_i'+c)\,\Gamma(k+\kappa_i'+1)}
-\log\tfrac{\Gamma(k+\kappa_i+c)\,\Gamma(\kappa_i+1)}
{\Gamma(\kappa_i+c)\,\Gamma(k+\kappa_i+1)}\Bigr],
\]
nonconstant unless $(\kappa'_i)$ is a permutation of $(\kappa_i)$,
both ratios nonincreasing, so every result above applies with the
log-concavity route to the factors.
\end{remark}

\begin{proposition}[The row--column system]\label{prop:rowcol}
Let $W=U_r$ and take $\mathcal A_1$ the $\HNB(\bk,c,w)$ operator of
the row differences.  On the column side, conditioning on the column
differences $\boldsymbol\lambda$ makes the columns independent given
the free column level $W'$, and
\[
\begin{aligned}
\pi_2(k)&=\PP\bigl(U_r=k\bigm|D_{\mathrm{col}}=\boldsymbol\lambda\bigr)\\
&=\sum_{w'\ge0}
\frac{\prod_{j}G_j(w'+\lambda_j)}{\sum_{v}\prod_jG_j(v+\lambda_j)}\,
\Bigl(\mathop{\ast}_{j=1}^{c}\ \mathcal L\bigl(X^{(r)}_j\bigm|
\textstyle\sum_iX^{(i)}_j=w'+\lambda_j\bigr)\Bigr)(k),
\end{aligned}
\]
with $G_j$ the pmf of column $j$'s sum: a mixture of across-column
convolutions of one-coordinate conditional laws, computable without
any table count.  Take $\mathcal A_2$ the canonical birth--death
operator of $\pi_2$: $b_2\equiv1$, $d_2(0)=0$ and
$d_2(k{+}1)=\pi_2(k)/\pi_2(k{+}1)$ on the support, so that
$\rho_2(k)=\pi_2(k{+}1)/\pi_2(k)$ and $\pi_2=\pi_{1,1}$.  The joint
conditioning has support on the arithmetic progression
$rk+\sum_i\kappa_i\equiv\sum_j\lambda_j\pmod c$, of step
$c/\gcd(r,c)$; on it, the coarsened ratios are the products of the
$c/\gcd(r,c)$ consecutive fine ratios.  If the coarse ratio of
$\pi_2$ is at most $\bar\rho<1$ from some point on, then
Theorems~\ref{thm:sysmaster}--\ref{thm:sysmatch} apply to every
$(w,\mu)$ with $w\rho_1^{1-\mu}\rho_2^{\mu}\le\bar\rho<1$ eventually,
with the factors taken from the exact rate-generic formula
\eqref{eq:gammaexact} and $\tilde\Lambda\le2\tilde\Gamma$---the
log-concavity route is not available, $\pi_2$ being a mixture that
can fail log-concavity.
\end{proposition}

\begin{proof}
Given the column differences, the column sums are $V_j=W'+\lambda_j$
with $W'$ the level, and by independence of the columns
$\PP(W'=w'\mid D_{\mathrm{col}}=\boldsymbol\lambda)\propto
\prod_jG_j(w'+\lambda_j)$; given the column sums the columns are
independent with laws $\mathcal L(Z_j\mid\sum_iX^{(i)}_j=V_j)$, and
$U_r=\sum_jX^{(r)}_j$ is their convolution in the $r$-th coordinate,
which is the display.  On the joint fiber $U_i=W+\kappa_i$ and
$\sum_iU_i=\sum_jV_j=cW'+\sum_j\lambda_j$, so
$rW+\sum_i\kappa_i\equiv\sum_j\lambda_j\pmod c$, a progression of
step $c/\gcd(r,c)$; the ratio of consecutive support points is the
product of the intermediate fine ratios.  The canonical pair has
$d_2(k{+}1)=1/\rho_2(k)\ge1/\bar\rho$ from some point on and
$d_2>0$ before it, and $D_1\ge1$, so
$\inf_{n\ge1}\tilde d_\mu(n)>0$ for every $\mu$, and
Theorem~\ref{thm:sysmaster} applies as stated;
Theorem~\ref{thm:sysmatch} requires only the ratios.
\end{proof}

The joint target itself has an exact form: its pmf is the table-count
sequence along the margin ray.

\begin{proposition}[The fiber shadow]\label{prop:shadow}
Suppose the entries are independent with product-form odds,
$X^{(i)}_j\sim\Geom(1-x_iy_j)$.  On the joint support---the
progression of Proposition~\ref{prop:rowcol}, of step $h=c/g$,
$g=\gcd(r,c)$, least element $k_0$---the margins are affine in $W$:
$U_i=W+\kappa_i$ and $V_j=w'(W)+\lambda_j$ with
$w'(k)=\bigl(rk+\sum_i\kappa_i-\sum_j\lambda_j\bigr)/c$.  Then, with
$\br(k)=(k+\kappa_i)_i$, $\bc(k)=(w'(k)+\lambda_j)_j$, and
$\Omega=\bigl(\prod_ix_i\bigr)^{c/g}\bigl(\prod_jy_j\bigr)^{r/g}$,
\[
\PP\bigl(W=k\ \bigm|\ D_{\mathrm r},\,D_{\mathrm c}\bigr)\;\propto\;
\bigl|\T\bigl(\br(k),\bc(k)\bigr)\bigr|\;\Omega^{\,(k-k_0)/h}.
\]
\end{proposition}

\begin{proposition}[The target margins are an atom of the shadow]\label{prop:atom}
Let the target margins be $\br^\ast,\bc^\ast$ with entries
$r^\ast_i,c^\ast_j$, and choose the offsets from the target itself:
$\kappa_i=r^\ast_i-r^\ast_r$ and $\lambda_j=c^\ast_j-c^\ast_c$.  Then
$k^\ast=r^\ast_r$ lies on the progression of
Proposition~\ref{prop:rowcol}, with $\br(k^\ast)=\br^\ast$,
$w'(k^\ast)=c^\ast_c$, and $\bc(k^\ast)=\bc^\ast$: conditioning on
both difference systems and on $\{W=k^\ast\}$ is exactly conditioning
on row \emph{and} column sums equaling their targets, and
\[
\PP\bigl(W=k^\ast\bigm|D_{\mathrm r},D_{\mathrm c}\bigr)
=\frac{\bigl|\T(\br^\ast,\bc^\ast)\bigr|\,\Omega^{(k^\ast-k_0)/h}}
{\sum_{k\in L}\bigl|\T(\br(k),\bc(k))\bigr|\,\Omega^{(k-k_0)/h}}\,.
\]
Every bound on the shadow localizes at $k^\ast$ to a two-sided
estimate of $|\T(\br^\ast,\bc^\ast)|$ against its ray neighbors.
\end{proposition}

Every defect of the shadow against either operator is an exact
cross-conditioning covariance---the commutator of the two
conditionings, entering as a Stein input.

\begin{lemma}[Cross-covariance defect identities]\label{lem:crosscov}
Let $g$ be integrable given $D_{\mathrm r}$ with
$\E\bigl[g(W)\mid D_{\mathrm r}=\bk\bigr]=0$---for exchangeable
columns, by Theorem~\ref{thm:target}, any bracket
$g=w_{\mathrm r}b_{\mathrm r}f(\cdot{+}1)-d_{\mathrm r}f$ of the row
equation at its own parameter $w_{\mathrm r}$ with $f$ bounded, or
its $h$-step coarsening.  Then, for the shadow law
$P_\bullet=\mathcal L(W\mid D_{\mathrm r},D_{\mathrm c})$,
\[
\E_{P_\bullet}\,g(W)\;=\;
\frac{\mathrm{Cov}_{\,\mid D_{\mathrm r}}\bigl(g(W),\
\ind{D_{\mathrm c}=\boldsymbol\lambda}\bigr)}
{\PP\bigl(D_{\mathrm c}=\boldsymbol\lambda\bigm|D_{\mathrm r}=\bk
\bigr)}\,,
\]
and symmetrically with rows and columns exchanged, the column
brackets being annihilated under $D_{\mathrm c}$ alone by the
canonical operator of $\pi_2$.
\end{lemma}

\begin{corollary}[The rate-moment inputs are commutators]\label{cor:crossdef}
For the shadow in Theorem~\ref{thm:sysmaster}, with exchangeable
columns so that the row bracket is annihilated under $D_{\mathrm r}$
alone: at $\mu=0$, with $(b_{\mathrm r},d_{\mathrm r})$
the coarse row pair $(B_1,D_1)$ on the progression of step $h$, at
its own parameter $w=w_{\mathrm r}^{h}$, and $f\equiv1$, the moment
defect $|\beta_{\tilde d}-w\beta_{\tilde b}|$ equals the normalized
covariance of $(d_{\mathrm r}-w_{\mathrm r}^{h}b_{\mathrm r})(W)$ with
$\ind{D_{\mathrm c}=\boldsymbol\lambda}$ given $D_{\mathrm r}$; at
$\mu=1$, $w=1$, it is the symmetric column expression.  The system's
rate-moment inputs are exactly the two conditioning commutators, and
any bound on them prices, by Proposition~\ref{prop:shadow}, the
table-count sequence itself.
\end{corollary}

Expanded---the bracket has conditional mean zero under
$D_{\mathrm r}$ alone---the $\mu=0$ input is computable factor by
factor:
\[
\bigl|\beta_{\tilde d}-w_{\mathrm r}\beta_{\tilde b}\bigr|_{\mu=0}
=\frac{\Bigl|\sum_{n}\PP(W{=}n\mid D_{\mathrm r}{=}\bk)\,
(d_{\mathrm r}{-}w_{\mathrm r}b_{\mathrm r})(n)\,
\PP(D_{\mathrm c}{=}\boldsymbol\lambda\mid W{=}n,
D_{\mathrm r}{=}\bk)\Bigr|}
{\sum_{n}\PP(W{=}n\mid D_{\mathrm r}{=}\bk)\,
\PP(D_{\mathrm c}{=}\boldsymbol\lambda\mid W{=}n,
D_{\mathrm r}{=}\bk)}\,,
\]
the row-conditioned law closed (Theorem~\ref{thm:target}), the
conditional probabilities
$\PP(D_{\mathrm c}{=}\boldsymbol\lambda\mid W{=}n,D_{\mathrm r}{=}\bk)$
ratios of weighted table counts over the row-sum fiber
$\br(n)$.

Assembled, the rungs give the table its main theorem.

\begin{theorem}[Row--column main theorem for contingency tables]\label{thm:tablemaster}
Let $P_\bullet$ be the fiber shadow of Proposition~\ref{prop:shadow},
whose pmf is the table-count sequence along the margin ray, with the
rates of the row--column system of Proposition~\ref{prop:rowcol}
under its tail assumption.  For every $\mu\in[0,1]$ such that the
ratio match $w^\ast(\mu)=\beta_{\tilde d}/\beta_{\tilde b}$ of
Theorem~\ref{thm:sysmatch}(i) satisfies
$w^\ast(\mu)\rho_1^{1-\mu}\rho_2^{\mu}\le\bar\rho<1$ from some point
on, with $\tilde\Gamma$ the factor \eqref{eq:gammaexact} of
$\pi_{w^\ast(\mu),\mu}$ and $\tilde\Lambda\le2\tilde\Gamma$,
\[
\dtv\bigl(P_\bullet,\ \pi_{w^\ast(\mu),\mu}\bigr)\;\le\;
\tilde\Lambda\;w^\ast\beta_{\tilde b}\;
\dwass\bigl(P_\bullet^{\tilde b}{+}1,\ P_\bullet^{\tilde d}\bigr),
\]
and, for exchangeable columns, the rate defect
$|\beta_{\tilde d}-w\beta_{\tilde b}|$ evaluated at the one-sided
parameter ($w=w_{\mathrm r}^{h}$ at $\mu=0$, $w=1$ at $\mu=1$) is the
conditioning commutator of Corollary~\ref{cor:crossdef}.
\end{theorem}

The one-dimensional family cannot always spend the column
commutator: a contrast that is affine in $k$ is absorbed by the
$w$-tilt (Lemma~\ref{lem:sysid}).

\begin{theorem}[Stein system summary]\label{thm:system}
Under the hypotheses of Theorem~\ref{thm:tablemaster}, for every
$\mu\in[0,1]$ admitted there, at the ratio match
$w^\ast(\mu)=\beta_{\tilde d}/\beta_{\tilde b}$ with
$\beta_{\tilde b}=\E\,\tilde b(W)$, $\beta_{\tilde d}=\E\,\tilde
d(W)$ under $P_\bullet$,
\[
\dtv\bigl(P_\bullet,\ \pi_{w^\ast(\mu),\mu}\bigr)\;\le\;
\tilde\Lambda\;w^\ast\beta_{\tilde b}\;
\dwass\bigl(P_\bullet^{\tilde b}{+}1,\ P_\bullet^{\tilde d}\bigr),
\]
$\tilde\Lambda\le2\tilde\Gamma$ with $\tilde\Gamma$ the explicit
factor \eqref{eq:gammaexact} of the interpolated target; at the
endpoints $\mu\in\{0,1\}$ the one-sided rate defects are the
covariances of each one-sided Stein bracket with the \emph{other}
conditioning (Corollary~\ref{cor:crossdef}); and every such bound
localizes at the atom $k^\ast$ to a two-sided estimate of
$P_\bullet(k^\ast)$, the normalized count of
Proposition~\ref{prop:atom}, within the bound of $\pi_{w^\ast(\mu),\mu}(k^\ast)$.
\end{theorem}

\subsection*{The general contingency table}

A uniform table is a geometric array conditioned on both margins;
under the product-form representation $y_{ij}=x_iy_j$, unequal
column margins are column heterogeneity; the column direction alone
is Theorem~\ref{thm:D}, exact at equal column factors and sharp on
the sparse near-equal wedge, and the row--column system of
Theorems~\ref{thm:tablemaster} and~\ref{thm:system} is the
two-directional certificate.  Balanced sparse designs are its
domain.  Beyond it the bounds are vacuous, for one
reason: the expansion prices a relative question absolutely.  A
sampler needs the ratio of neighboring fiber masses along a margin
ray to relative precision of order one over the margins; the
certificate carries the unconditional moments $\nu^{b}\ge c$, which
grow with the margins while the masses they must resolve shrink:
$\nu^{b}_{S^c}\ge c^{|S^c|}$, $\theta\ge1$ (the term $L=\varnothing$,
$s=\bkk$), and the weights $1/(1-y_{ij})=1+\E X_{ij}$ grow with the
cell means, so at fixed column factors the right side of
Theorem~\ref{thm:colexpansion} grows with the row means and passes
one (\S\ref{sec:computations} records where).  A
certificate above one is the interval $[0,1]$.  The theorem that
would make the system a sampling
instrument---a bound on the log-ratio of the shadow along a ray,
normalized by the fiber mass rather than on its total variation---is
Open Problem~\ref{op:fullmargin}; its log-concavity input holds for the
polytope's volume by Barvinok's Brunn--Minkowski inequalities
\cite{Barv07} and is, for the count, a question of Lorentzian type
\cite{BrH20}.

\section{Computations}\label{sec:computations}

Every bound of the paper is evaluated here against the exact distance
it bounds.  The ancillary \path{verify_companion/} holds one script
per result, named after it (\path{thm_2_06_factors.py} for
Theorem~\ref{thm:factors}, \path{cor_2_14_nbsum.py} for
Corollary~\ref{cor:nbsum}, and so on; the lettered theorems of
\S\ref{sec:intro} are verified through their body results), each
printing the quantities compared and asserting the values quoted
below; \path{run_all.py} runs the $55$ scripts in order in about
five minutes with NumPy alone, the optimal-transport couplings of
Lemma~\ref{lem:vanishing} and Theorem~\ref{thm:expansion} being
computed exactly as minimum-cost flows on the lattice.  The
first-difference factor $\Lambda$
is computed exactly: each increment $g_A(m+1)-g_A(m)$ is linear in
the indicator of $A$, so its supremum over $A$ is the sum of the
positive, or of the negative, coefficients.  Laws on $\N$ are truncated at an $N$
with tail mass below $10^{-12}$; exact distances are computed from
the truncated pmfs, conditional laws by convolution and
normalization on the fiber.

\subsection*{The operator and its factors}
\begingroup\sloppy\noindent(\path{thm_2_05_solution}, \path{thm_2_06_factors},
\path{rem_2_07_placement}, \path{thm_2_08_nbfactors},
\path{rem_2_09_which}, \path{rem_2_10_signfail}.)\par\endgroup
The characterizing identity, the solution formula and uniqueness are
checked to $10^{-14}$ on $\NB(3,0.6)$ and $\HNB((2,0),4,0.3)$.  The
exact factor \eqref{eq:gammaexact} is one pass over the pmf, $O(N)$;
the supremum over test sets, for comparison, solves the Stein
equation for every interval $[0,m]$ and for $300$ random sets,
$O(N^{2})$, and agrees with \eqref{eq:gammaexact} to $10^{-9}$:
$\Gamma=0.4875$ on $\NB(4,\tfrac12)$ and $\Gamma=0.1788$ on
$\HNB((1,0),3,0.45)$, where the closed bound \eqref{eq:gammabound}
gives $2.67$ and $1.52$.  $\Gamma=1$ at $w=0.3,0.6,0.9$ for the
geometric law, as in closed form.  At $w=\tfrac12$,
$\Gamma=0.427,0.202,0.100$ for $c=5,20,80$.  Over the grid
$2\le c\le50$, $w\in\{0.1,0.2,\dots,0.9,0.95,0.97\}$,
\[
\Gamma_{c,w}\sqrt{cw}\in[0.425,\,0.769],
\]
while the closed bound \eqref{eq:gammabound}, evaluated in exact
rational arithmetic, exceeds the exact factor by a factor up to
$1.1\cdot10^{3}$: the closed form is valid everywhere and the exact
form is the one to evaluate.  On the same grid, extended by the small
shapes $\rho\in\{0.3,0.5,1\}$ at $w\le0.9$ ($566$ points), the exact
factor lies below each closed bound of Theorem~\ref{thm:nbfactors}
that applies, and on the subgrid $\rho\in\{0.3,0.5,1,2,5,10\}$,
$w\in\{0.1,0.3,0.5,0.7\}$ the exact first-difference factor
lies below each of the five, and the map of Remark~\ref{rem:which}
reads: for $\Gamma$, the bound
of \cite{CD19} is the smallest at $503$ points and $1/(1-w)$ at $63$
(small means, $(\rho+1)w<\pi(1-w)$), and \eqref{eq:gammabound},
compared at the $\rho>1$ points where it is proved, is never the
smallest; the smallest closed bound still
exceeds the exact factor by a factor up to $16$, at $w=0.97$, where
the closed forms carry $(1-w)^{-1/2}$ and the exact factor does not.
For $\Lambda$, $1/(w\rho)$ is the smallest at $537$ points,
$(1+w)/(1-w)^{2}$ at $12$ ($w\le0.2$, $\rho\le7$) and
$2\Gamma$ at $17$ ($\rho\le1$ with $0.3\le w\le0.9$, and
$(\rho,w)=(2,0.2),(3,0.2)$).  The scaling is the one expected at the
mode, where $F\approx G\approx\tfrac12$ and
$p_{m^\ast}\approx(2\pi\Var N)^{-1/2}$, so that
\eqref{eq:gammaexact} is of order $\sqrt{\Var N}/d(m^\ast)$, the
analogue of $\sigma/\nu$ in \cite{DeS26}.  On that subgrid the exact
first-difference factor lies between $0.97\min\{\Gamma,1/(w\rho)\}$
and that minimum: the smaller of the two proved bounds is nearly
sharp there.
On $(r,c,\bk,w)=(2,6,(1,0),\tfrac{11}{20})$ the increment of
Remark~\ref{rem:signfail} is $\Delta g_2(1)=5.0879\cdot10^{-6}$,
equal to its closed form to $10^{-15}$, and $284$ pairs $(j,m)$ with
$m\ne j$, $1\le j,m\le40$, have $\Delta g_j(m)>0$.

\subsection*{The two biasings}
\begingroup\sloppy\noindent(\path{thm_2_11_nbmaster}, \path{rem_2_12_twolines}.)\par\endgroup
The identity $\E(\mathcal Af)(W)=w\beta_b\,\E f(W^{b}+1)-\beta_d\,\E
f(W^{d})$ holds to $10^{-10}$ for a random bounded $f$ on a sum of
three geometrics, and $W^{b}+1\overset{d}{=}W^{d}$ to $10^{-12}$ on
$\NB(3,0.55)$.  At the matched rate, with $\Lambda$ replaced by its
proved bound $\min\{2\Gamma,1/(w^\ast\rho)\}$ (Theorem~\ref{thm:nbfactors}),
Theorem~\ref{thm:nbmaster} gives $0.067$ against the exact distance
$0.0059$ for the geometric sum with $y=(0.55,0.62,0.48)$, $\rho=3$
($w^\ast=0.557$, $\E W=3.78$, $\dwass=0.030$, the bound $1/(w^\ast\rho)=0.598$
being the smaller), and $0.86$ against $0.111$ for a Poisson count
whose rate is a sample of $4000$ lognormal draws of median $4$ and
log-spread $0.35$ ($w^\ast=0.585$, $\E W=4.23$, $\dwass=0.357$).

\subsection*{Theorem~\ref{thm:A}}
\begingroup\sloppy\noindent(\path{lem_2_02_geomrow}, \path{cor_2_14_nbsum},
\path{cor_2_24_occupancy}, \path{cor_2_25_randomstage}.)\par\endgroup
At $w=\tfrac12$:
\begin{center}\small\setlength{\tabcolsep}{5pt}
\begin{tabular}{lccccc}
perturbation & $c$ & $\eta$ & $\dtv$ exact & bound & ratio\\\hline
one stage off & $5$ & $0.05$ & $2.68\cdot10^{-2}$ & $4.75\cdot10^{-2}$ & $1.77$\\
one stage off & $5$ & $0.10$ & $5.84\cdot10^{-2}$ & $1.07\cdot10^{-1}$ & $1.83$\\
one stage off & $5$ & $0.20$ & $1.40\cdot10^{-1}$ & $2.85\cdot10^{-1}$ & $2.03$\\
one stage off & $20$ & $0.05$ & $1.38\cdot10^{-2}$ & $2.24\cdot10^{-2}$ & $1.62$\\
one stage off & $20$ & $0.10$ & $3.08\cdot10^{-2}$ & $5.04\cdot10^{-2}$ & $1.64$\\
one stage off & $20$ & $0.20$ & $7.90\cdot10^{-2}$ & $1.35\cdot10^{-1}$ & $1.70$\\
one stage off & $80$ & $0.05$ & $6.99\cdot10^{-3}$ & $1.11\cdot10^{-2}$ & $1.58$\\
one stage off & $80$ & $0.10$ & $1.57\cdot10^{-2}$ & $2.49\cdot10^{-2}$ & $1.59$\\
one stage off & $80$ & $0.20$ & $4.13\cdot10^{-2}$ & $6.64\cdot10^{-2}$ & $1.61$\\
symmetric drift & $5$ & $0.02$ & $9.03\cdot10^{-4}$ & $5.13\cdot10^{-2}$ & $56.9$\\
symmetric drift & $5$ & $0.05$ & $5.66\cdot10^{-3}$ & $1.29\cdot10^{-1}$ & $22.8$\\
symmetric drift & $5$ & $0.10$ & $2.29\cdot10^{-2}$ & $2.64\cdot10^{-1}$ & $11.6$\\
symmetric drift & $20$ & $0.02$ & $1.44\cdot10^{-3}$ & $8.50\cdot10^{-2}$ & $58.9$\\
symmetric drift & $20$ & $0.05$ & $9.04\cdot10^{-3}$ & $2.14\cdot10^{-1}$ & $23.6$\\
symmetric drift & $20$ & $0.10$ & $3.65\cdot10^{-2}$ & $4.34\cdot10^{-1}$ & $11.9$\\
symmetric drift & $80$ & $0.02$ & $2.74\cdot10^{-3}$ & $1.61\cdot10^{-1}$ & $58.9$\\
symmetric drift & $80$ & $0.05$ & $1.72\cdot10^{-2}$ & $4.05\cdot10^{-1}$ & $23.6$\\
symmetric drift & $80$ & $0.10$ & $6.89\cdot10^{-2}$ & $8.23\cdot10^{-1}$ & $11.9$\\
\end{tabular}
\end{center}
One stage displaced is $y_1=w+\eta$, $y_j=w$ for $j\ge2$; the
symmetric drift is $y_j=w+\eta\,(2j/(c-1)-1)$, $0\le j<c$.  On the
first face the ratio of bound to distance settles at $1.58$ as $c$
grows, a value we read as the extremal constant $\pi/2$ of
\cite{BX06} at work and not as a theorem.  On the second face the
distance is second order in $\eta$, the bound first order, and the
ratio is $1.2/\eta$ at every $c$.  Occupancy
(Corollary~\ref{cor:occupancy}): $c=20$ boxes with $\beta\epsilon_j$
on a uniform grid in $[0.6,0.8]$, $\bar y=0.498$, one ball per box on
average: bound $0.212$ against distance $0.0089$; on $[0.5,0.9]$,
$0.437$ against $0.037$.  The coupon collector at $c=10$, against
$\NB(10,0.45)$: bound $2.76$, distance $0.56$.

\subsection*{Theorem~\ref{thm:coupling}}
\begingroup\sloppy\noindent(\path{def_2_16_summandbias}, \path{thm_2_17_couplingbody},
\path{rem_2_18_independent}, \path{cor_2_20_monotone}.)\par\endgroup
On a sum of a coupled geometric pair and an independent geometric
the mixture identities of Definition~\ref{def:summandbias} hold to
$10^{-16}$.  On the pair instances of Theorem~\ref{thm:dep} the law of the
partner given $X_j$ is stochastically increasing, the coupling term
under the monotone coupling, computed from the conditional laws,
equals $\sum_j(1-y_j)|\Cov(X_j,R_j)|$ to $10^{-8}$, and the bound of
Corollary~\ref{cor:monotone} replaces the $Q$ bound as follows.
Twenty summands with $y_j=0.4$ at $w=y$: at
$\vartheta=0,\tfrac14,\tfrac12,1$ the distances are $0$, $0.051$,
$0.097$, $0.5$, the $Q$ bounds $2.2$, $3.2$, $4.2$, $6.1$, and the
covariance bounds $0$, $0.42$, $0.83$, $1.67$: exact at
independence, and five times sharper at $\vartheta=\tfrac12$.  The
heterogeneous instance at $w=\bar y$: $0.61$ in place of $3.5$,
against $0.061$.  The Poisson regime ($\lambda=2$): $0.51$ in place
of $0.57$ at $w=y$ and $0.48$ in place of $0.54$ at $w=0$, against
$0.126$ and $0.132$; the pairs there are nearly indicators and the
two couplings nearly coincide.  For non-geometric summands, thirty
independent Bernoulli variables with $p_j$ on a uniform grid in
$[0.05,0.2]$, $\lambda=\sum_jp_j/(1+p_j)=3.29$: the coupling term is
$\sum_j2p_j^{2}/(1+p_j)=0.918$, and the bound is $0.62$ against the
distance $0.110$ at $w=0$ and $0.32$ against $0.070$ at $w=\bar y$.

\subsection*{Theorem~\ref{thm:dep}}
\begingroup\sloppy\noindent(\path{thm_2_19_aggnb}, \path{prop_2_21_aggnblambda},
\path{rem_2_22_nbpoisson}.)\par\endgroup
Pairs $(X_{2k-1},X_{2k})$ are coupled by taking the comonotone
coupling of their geometric laws with probability $\vartheta$ and
independent copies otherwise, which keeps both marginals geometric;
$B_j$ is the pair, so $\beta=0$, and $W$ is the sum over the pairs.
The factor $\Lambda_w$ is the proved $\min\{2\Gamma_w,1/\lambda\}$,
which is $1/\lambda$ in every case below with $w>0$, and
$(1-e^{-\lambda})/\lambda$ at $w=0$.  Twenty summands with
$y_j=0.4$, $w=y$ (so $M_w=0$): at $\vartheta=0,\tfrac12,1$ the exact
distances are $0$, $0.097$, $0.5$ and the bounds $2.2$, $4.2$, $6.1$;
the bound counts $2\,\E X_i\,\E X_j$ per ordered pair even at
$\vartheta=0$, as $b_1+b_2$ does for an independent neighbor.
Sixteen summands with $y_j$ drawn from $[0.3,0.5]$ and
$\vartheta=0.3$: at $w=0,\,0.3,\,\bar y=0.396,\,0.5$ the distances
are $0.457$, $0.161$, $0.061$, $0.164$ and the bounds $5.8$, $3.7$,
$3.5$, $3.6$.  In the Poisson regime, one hundred summands with
$y_j=0.02$ and $\vartheta=\tfrac12$ ($\lambda=2$): at $w=y$ the bound
is $0.57$ against $0.126$, and at the endpoint $w=0$, with the
Poisson factors, $0.54$ against $0.132$.  On the independent
occupancy instance above, at $w=\bar y$, the weights $y_j/(1-y_j)$
give $0.106$ where Corollary~\ref{cor:nbsum} gives $0.212$, against
$0.0089$.  Proposition~\ref{prop:aggnblambda} gives $4.2$ and $3.7$
on the first two pair instances, and the distance of
Remark~\ref{rem:nbpoisson} is $0.026$, $0.0023$, $0.167$ against
$1.4w/\sqrt{1-w}=0.148$, $0.014$, $0.99$ at
$(\rho,w)=(20,0.1),(200,0.01),(5,\tfrac12)$.  The bound costs the
cross-moments of the neighborhoods and one evaluation of the
factors, and nothing of the joint law beyond them.

\subsection*{Theorem~\ref{thm:B}}
\begingroup\sloppy\noindent(\path{cor_2_26_mixture}, \path{rem_2_27_secondroute}.)\par\endgroup
The biasing identities $W^{d}\overset{d}{=}1+W^{s}$ and
$\mathcal L(W^{b})=w\mathcal L(1+W^{s})+(1-w)\mathcal L(W)$ hold to
$10^{-10}$ in every case below, and $\beta_d/\beta_b=w$.  With
$\Lambda$ replaced by its proved bound
$\min\{2\Gamma_{\rho,w},v/m^{2}\}$, which is $v/m^{2}$ in every case,
and with the second route of the remark after the corollary in the
last column:
\begin{center}\small\setlength{\tabcolsep}{5pt}
\begin{tabular}{lccccccc}
$\Lambda$ & $m$ & $v$ & $\rho$ & $w$ & $\dtv$ exact & corollary & second route\\\hline
uniform on $\{1.9,2.1\}$ & $2$ & $2.01$ & $400$ & $0.005$ & $3.5\cdot10^{-6}$ & $1.4\cdot10^{-5}$ & $0.012$\\
uniform on $\{1,3\}$ & $2$ & $3$ & $4$ & $\tfrac13$ & $0.0268$ & $0.128$ & $1.07$\\
uniform on $\{2,4\}$ & $3$ & $4$ & $9$ & $\tfrac14$ & $0.0122$ & $0.0481$ & $0.74$\\
uniform on $\{1,2,4\}$ & $2.33$ & $3.89$ & $3.5$ & $0.4$ & $0.0221$ & $0.126$ & $1.39$\\
lognormal, median $2$, log-spread $\tfrac12$ & $2.27$ & $3.73$ & $3.52$ & $0.392$ & $0.0181$ & $0.0851$ & $1.35$\\
Gamma, shape $3$, scale $1$ & $3$ & $6$ & $3$ & $\tfrac12$ & $<10^{-11}$ & $<10^{-10}$ & $1.99$\\
\end{tabular}
\end{center}
The continuous rates are discretized on a fine grid.  The ratio of
the corollary to the distance is between $3.8$ and $5.8$; the
exact $\Lambda$ is between $0.80$ and $0.93$ of the proved
bound $v/m^{2}$.  The second route is first order in $\Var\Lambda$
and is above the corollary in every case, by a factor $886$ at the
nearly constant rate.  The bound costs one Wasserstein distance
between two pmfs, $O(N)$.

\subsection*{Theorem~\ref{thm:C}}
\begingroup\sloppy\noindent(\path{thm_3_02_target}, \path{prop_3_03_kemp},
\path{lem_3_05_basic}, \path{rem_3_06_mode},
\path{prop_3_08_poissonlimit}, \path{cor_3_10_bessel},
\path{thm_3_11_fixedpoint}, \path{thm_3_12_master},
\path{lem_3_13_phantom}, \path{lem_3_14_sizebias},
\path{lem_3_15_geomchar}, \path{prop_3_17_factorbound}.)\par\endgroup
At $r=2$ the conditional law of $R_2$ given $R_1-R_2=k$, computed
from the product of two negative-binomial pmfs, equals
$\HNB((k,0),c,q_1q_2)$ to $10^{-15}$; the fiber probability
\eqref{eq:fiber} matches direct summation; the balance identity holds
to $10^{-15}$; the hypergeometric normalizer \eqref{eq:zhathyper},
summed as a series, agrees with the direct sum $\Zhat_c$ to
$10^{-10}$ ($113.859$ at $c=4$, $q=(0.55,0.4)$, $k=3$).  In the
Poisson limit at $\nu=6$, $k=3$, the distance from
$\HNB((k,0),c,\nu/c^{2})$ to the Bessel law is $1.03\cdot10^{-1}$, $5.08\cdot10^{-2}$,
$2.51\cdot10^{-2}$, $1.25\cdot10^{-2}$ at $c=20,40,80,160$: the rate
is $O(1/c)$.  The Bessel operator of Corollary~\ref{cor:bessel}
annihilates the Bessel law to $2\cdot10^{-15}$.

\subsection*{Theorem~\ref{thm:D}}
\begingroup\sloppy\noindent(\path{thm_3_19_colexpansion},
\path{thm_3_19_colexpansion_at_scale},
\path{lem_3_21_convmom}, \path{thm_3_22_aggconv},
\path{lem_3_24_transfer}, \path{lem_3_25_vanishing},
\path{thm_3_26_expansion}, \path{cor_3_27_productform},
\path{cor_3_28_capped}, \path{thm_3_29_tvregister}.)\par\endgroup
At $r=2$, $c=3$, $Y=\bigl(\begin{smallmatrix}0.30&0.34&0.27\\
0.22&0.25&0.20\end{smallmatrix}\bigr)$, $k=2$, the bound is $1.11$
against the exact distance $1.4\cdot10^{-3}$ ($\theta=4.97$); with
exchangeable columns both are zero.  The fiber constant is computed
by dynamic programming on every reduced system with at most $r$
columns removed, $\sum_{l\le r}\binom cl$ convolutions of length
$N$; at $c=30$, $r=2$ these are $466$ systems and most of the
minute the script takes.  In the Barvinok calibration $y_{ij}=q_i\beta_j$ with
$q=(0.30,0.24)$, $c=30$, $\beta_j=1+\eta U_j$, $U_j$ uniform on
$[-1,1]$ (one fixed draw), conditioning on the modal difference:
\begin{center}\small\setlength{\tabcolsep}{6pt}
\begin{tabular}{ccccc}
$\eta$ & $\dtv$ exact & bound & ratio & $\theta$\\\hline
$0$ & $0$ & $0$ & --- & $1.037$\\
$0.02$ & $0.00919$ & $0.163$ & $17.8$ & $1.038$\\
$0.05$ & $0.0226$ & $0.407$ & $18.1$ & $1.040$\\
$0.10$ & $0.0437$ & $0.810$ & $18.5$ & $1.044$\\
\end{tabular}
\end{center}
with $\Gamma=0.0160$ throughout: the bound degrades continuously to
exactness and is $18$ times the truth across the wedge.  At scale
(\path{thm_3_19_colexpansion_at_scale}): $r=2$, $c=3$, $\beta=(\tfrac{19}{20},\tfrac{21}{20},1)$,
$k_1=1$; at $q=(\tfrac23,\tfrac35)$, row means $(6.1,4.5)$, the bound
is $2.01$ ($\Gamma=0.197$, $\theta=2.91$; $0.690$ with $\theta$
replaced by its floor $1$) against the exact distance $0.0034$; at
$q=(\tfrac89,\tfrac67)$, row means $(27.4,19.4)$, it is $25.8$
($\Gamma=0.109$, $\theta=4.25$; floor $6.09$) against $0.053$.  The
floor alone passes one between the two.

\subsection*{The row--column system}
\begingroup\sloppy\noindent(\path{lem_4_02_sysid}, \path{lem_4_03_syslc},
\path{thm_4_04_sysmaster}, \path{thm_4_05_sysmatch},
\path{rem_4_06_tworow}, \path{prop_4_07_rowcol},
\path{prop_4_08_shadow}, \path{prop_4_09_atom},
\path{lem_4_10_crosscov}, \path{cor_4_11_crossdef},
\path{thm_4_12_tablemaster}, \path{thm_4_13_system}.)\par\endgroup
On a $2\times3$ product-form array, $x=(0.55,0.40)$,
$y=(0.60,0.50,0.45)$, the support of the fiber shadow is the
progression of Proposition~\ref{prop:shadow} (step $3$), its
weights are $\Omega^{(k-k_0)/h}$ times the counts, and the mixture
formula of Proposition~\ref{prop:rowcol} matches direct enumeration
to $10^{-7}$.  On a second array with odds near one,
$x=(0.92,0.88)$, $y=(0.96,0.93,0.90)$, so that the shadow spreads
over its progression, the contrast defect of
Theorem~\ref{thm:sysmatch}(ii) vanishes at $\mu^\ast=0.866$,
$w^\ast=0.664$, and the bound of Theorem~\ref{thm:tablemaster}, with
$\tilde\Lambda$ replaced by its proved bound $2\tilde\Gamma$, is
$382$, $3.49$, $0.108$ at $\mu=0,\tfrac12,1$ against the exact
distances $0.220$, $0.072$, $0.018$
($\tilde\Gamma=1.4\cdot10^{-4}$, $0.013$, $1.14$; the exact
$\tilde\Lambda$ is $1.4\cdot10^{-4}$, $0.013$, $1.02$): the column
end of the family is the informative one on this instance, and the
row end is vacuous.

\section{Open problems}\label{sec:open}

\begin{openproblem}[Sharp first-difference factor]\label{op:lambda}
Remark~\ref{rem:signfail} shows the singleton sign pattern fails for
degree-$r$ birth --- both as a proof route and as a fact.  On the
negative-binomial subgrid of \S\ref{sec:computations} the exact
$\Lambda$ lies between $0.97\min\{\Gamma,1/(w\rho)\}$ and that
minimum, so the smaller of the two proved bounds is nearly sharp
there; for $\HNB(\bk,c,w)$ no bound of the form $1/(w\rho)$ is
available.  Determine the sharp order of
$\Lambda$ for $\HNB(\bk,c,w)$, e.g.\ prove or refute
$\Lambda=\sup_j\Delta g_j(j)\,(1+o(1))$, and identify the analogue of the
Poisson-case equality.
\end{openproblem}

\begin{openproblem}[Second-order bound for mean-matched geometric sums]\label{op:secondorder}
For $W=\sum_jX_j$, $X_j\sim\Geom(1-y_j)$, with $w=\bar y$, the
first-order variation of the law of $W$ vanishes
(\S\ref{sec:nb}) and the exact distance to $\NB(c,w)$ is second
order in the spread in \S\ref{sec:computations}, while
Corollary~\ref{cor:nbsum} is first order.  Prove
$\dtv(W,\NB(c,w)+\text{signed correction})\le C\sum_j(y_j-w)^{2}$
with $C$ explicit, in the manner of \cite{BC02,CR06}.
\end{openproblem}

\begin{openproblem}[Factors for the Bessel law]\label{op:bessel}
Give $\Gamma$ and $\Lambda$ in closed form for the operator of
Corollary~\ref{cor:bessel}, and a two-biasings theorem for the Bessel
law with rates $b(n)=\lambda$, $d(n)=n(n+k)$.
\end{openproblem}

\begin{openproblem}[Multivariate fiber constants]\label{op:multivariate}
Give closed-form, non-DP bounds for $\theta,\tau_1,\omega$ of
Definition~\ref{def:costs} in terms of column moments --- equivalently, a
certified multivariate local limit theorem for lattice vectors of
differences --- tight enough to make \eqref{eq:expansion} practical for
large $c$.
\end{openproblem}

\begin{openproblem}[Local upgrade]\label{op:local}
Upgrade Theorem~\ref{thm:expansion} from total variation to a pointwise
(local) approximation $|\PP(W=n)-p_n|\le\varepsilon_n$ with certified
$\varepsilon_n$, which would give two-sided bounds for individual
table counts with both margins fixed.
\end{openproblem}

\begin{openproblem}[Single-shot full-margin certificates]\label{op:fullmargin}
Certify $\PP(\text{column sums}=\bc\mid\text{row sums}=\br)$ directly,
without column peeling, by a closed-form two-sided bound for general
$m\times n$ whose relative error is of order one over the margins.
The structural obstruction is the unconditional moments $\nu^{b}$ of
Theorem~\ref{thm:colexpansion}, which price a relative question
absolutely: on the instance of \S\ref{sec:computations} the floor
$\theta\ge1$ alone exceeds one at row means $(27,19)$.
\end{openproblem}

\begin{openproblem}[Dense regime]\label{op:dense}
As $w\to1^-$ with $c$ fixed, the target's mean diverges like
$r(c-1)/(1-w)$ and the geometry approaches the dense-margin regime of
\cite{CM10,DLP20}.  Connect the present certificates to the
Canfield--McKay asymptotics; in particular, determine the behavior of
$\Gamma$, $\Lambda$ as $w\to1^-$.
\end{openproblem}

\begin{openproblem}[Dependent arrays]\label{op:dependent}
Extend Theorem~\ref{thm:expansion} beyond independent columns ---
exchangeable columns and P\'olya-type reinforcement are natural, since
the $b$-bias ``add a column'' operation is intrinsic to urn dynamics
\cite{PRR13}.
\end{openproblem}

\section{Proofs}\label{sec:proofs}

\subsection*{Proofs for Section~\ref{sec:nb}}

\begin{proof}[Proof of Lemma~\ref{lem:geomrow}]
The generating function is $\bigl((1-q)/(1-qz)\bigr)^{c}$.
\end{proof}

\begin{proof}[Proof of Theorem~\ref{thm:solution}]
(i) \emph{Necessity.}  If $M\sim\pi$ then, using the balance
identity,
\[
\E[w\,b(M)g(M{+}1)]=\sum_nw\,b(n)p_ng(n{+}1)
=\sum_nd(n{+}1)p_{n+1}g(n{+}1)=\E[d(M)g(M)],
\]
all sums absolutely convergent for bounded $g$ since
$\E b(M)<\infty$.

\emph{Sufficiency.}  Suppose $\E(\mathcal{A}g)(M)=0$ for all bounded
$g$.  For $A\subseteq\N$ let $g_A$ be given by \eqref{eq:solution} with
$h=\ind{\cdot\in A}$; $g_A$ is bounded by part (ii) below and
Theorem~\ref{thm:factors}, so
$0=\E(\mathcal{A}g_A)(M)=\PP(M\in A)-\PP(N\in A)$ for every $A$.

(ii) Multiplying the Stein equation at $n$ by $p_n$ and using
$d(n)p_n=w\,b(n-1)p_{n-1}$ turns it into the telescoping
$w\bigl[b(n)p_ng(n{+}1)-b(n{-}1)p_{n-1}g(n)\bigr]=\tilde h(n)p_n$
(with the $n=0$ term reading $w\,b(0)p_0g(1)=\tilde h(0)p_0$, as
$d(0)=0$).  Summing to $m$ gives the first form of \eqref{eq:solution};
the second follows since $\sum_n\tilde h(n)p_n=0$.  Conversely the
formula solves the equation: by \eqref{eq:solution} at $m$ and at
$m-1$, and $d(m)p_m=w\,b(m{-}1)p_{m-1}$,
$w\,b(m)p_mg(m{+}1)-d(m)p_mg(m)=\sum_{n\le m}\tilde h(n)p_n-\sum_{n\le m-1}\tilde h(n)p_n=\tilde h(m)p_m$.
For uniqueness, a solution of the homogeneous equation satisfies
$g(n+1)=\tfrac{d(n)}{w\,b(n)}g(n)$, so
$|g(n+1)|=|g(1)|\prod_{m=1}^{n}\tfrac{d(m)}{w\,b(m)}$, and the
equation at $n=0$ reads $w\,b(0)g(1)=0$, so $g(1)=0$ and $g\equiv0$ on
$\{1,2,\dots\}$: the solution is unique outright, and bounded by
Theorem~\ref{thm:factors}(i).
\end{proof}

\begin{proof}[Proof of Theorem~\ref{thm:factors}]
(i) From \eqref{eq:solution} with $h=\ind{\cdot\in A}$,
\[
\begin{aligned}
w\,b(m)p_m\,g_A(m+1)&=\PP(N\in A\cap[0,m])-\PP(N\in A)F(m)\\
&=\PP(A\cap[0,m])\,G(m{+}1)-\PP(A\cap[m{+}1,\infty))\,F(m),
\end{aligned}
\]
a difference of two nonnegative terms, each at most $F(m)G(m{+}1)$;
hence $|g_A(m{+}1)|\le F(m)G(m{+}1)/(w\,b(m)p_m)$, with equality for
$A=[0,m]$.  Since $w\,b(m)p_m=d(m{+}1)p_{m+1}$, \eqref{eq:gammaexact}
follows.

(ii) Write $S_m:=F(m)G(m+1)/(d(m+1)p_{m+1})$; $t_n$ is strictly
decreasing by hypothesis.
\emph{Sub-mode branch.}  If $t_{m-1}>1$ (equivalently $x_m'>0$), then for
$1\le n\le m$ we have $p_{n-1}/p_n=1/t_{n-1}\le1/t_{m-1}$, so
$F(m)\le p_m\sum_{i\ge0}t_{m-1}^{-i}=p_m\,\tfrac{t_{m-1}}{t_{m-1}-1}$.
With $G\le1$, $p_m/p_{m+1}=1/t_m$, $t_{m-1}=w\,b(m{-}1)/d(m)$ and
$t_{m-1}-1=x_m'/d(m)$,
\[
S_m\le\frac{t_{m-1}}{(t_{m-1}-1)\,d(m{+}1)\,t_m}
=\frac{w\,b(m{-}1)}{x_m'\;d(m{+}1)\;\dfrac{w\,b(m)}{d(m{+}1)}}
=\frac{b(m{-}1)}{b(m)\,x_m'}\;\le\;\frac{1}{x_m'} .
\]
\emph{Super-mode branch.}  If $t_{m+1}<1$ (equivalently $y_{m+1}>0$),
then for $n\ge m{+}1$, $p_{n+1}/p_n\le t_{m+1}$, so
$G(m{+}1)\le p_{m+1}/(1-t_{m+1})$, and with $F\le1$ and
$1-t_{m+1}=y_{m+1}/d(m{+}2)$,
\[
S_m\;\le\;\frac{1}{d(m{+}1)\,(1-t_{m+1})}\;=\;\frac{d(m{+}2)}{d(m{+}1)\,y_{m+1}} .
\]
\emph{The origin.}  $S_0=p_0G(1)/(d(1)p_1)\le p_0/(d(1)p_1)=1/(w\,b(0))$.
\emph{Coverage.}  If some $m\ge1$ had $t_{m-1}\le1$ and $t_{m+1}\ge1$,
strict decrease would give $t_{m+1}<t_{m-1}\le1$, a contradiction; so
every $m\ge1$ lies in one of the two branches, giving
\eqref{eq:gammabound}.  The sub-mode branch, $t_{m-1}>1$, holds for finitely many $m$ since
$t_n\to w<1$; when $d(m)\to\infty$ the super-mode terms tend to $0$,
since $y_{m+1}/d(m{+}2)=1-t_{m+1}\to1-w$, so their supremum is a
maximum.

(iii) $\Lambda\le2\Gamma$ is trivial.  For the second bound, solve the
Stein equation for $g(m{+}1)$:
$g(m{+}1)-g(m)=\tfrac{d(m)-w\,b(m)}{w\,b(m)}\,g(m)+\tfrac{\tilde h(m)}{w\,b(m)}$.
If $d(m)\le w\,b(m)$ the coefficient has modulus
$1-\tfrac{d(m)}{w\,b(m)}\le1$; if $d(m)>w\,b(m)$ it equals
$\tfrac{d(m)}{w\,b(m)}-1\le\tfrac1w-1$ since $d\le b$.  With
$|\tilde h|\le1$ and $b(m)\ge b(0)$, \eqref{eq:lambdabound} follows
(the $m=0$ increment is $|g_A(1)|\le\Gamma$, also covered).
\end{proof}

\begin{proof}[Proof of Theorem~\ref{thm:nbmaster}]
For $A\subseteq\N$ let $g=g_A$ solve the Stein equation for $\pi=\pi_{b,d,w}$.  Then, by the definition of the biased laws,
\[
\PP(W\in A)-\pi(A)
=\E\bigl[w\,b(W)g(W{+}1)-d(W)g(W)\bigr]
=w\beta_b\,\E g(W^b{+}1)-\beta_d\,\E g(W^d).
\]
Add and subtract $w\beta_b\,\E g(W^d)$:
the first piece is $w\beta_b\bigl(\E g(W^b{+}1)-\E g(W^d)\bigr)$, bounded
in modulus by $w\beta_b\,\Lambda\,\dwass(W^b{+}1,W^d)$ because
consecutive increments of $g$ are bounded by $\Lambda$, hence
$|g(x)-g(y)|\le\Lambda|x-y|$ on $\{1,2,\dots\}$ (both $W^b{+}1\ge1$ and
$W^d\ge1$), and one may pass to the infimum over couplings; the second
piece is $(w\beta_b-\beta_d)\E g(W^d)$, bounded by
$\Gamma|w\beta_b-\beta_d|$.  Taking the supremum over $A$ gives
the bound; the matched choice $w^\ast=\beta_d/\beta_b$ kills the
second term.
\end{proof}

\begin{proof}[Proof of Proposition~\ref{prop:aggnblambda}]
Write $\Lambda=\Lambda_{w^\ast}$ and let $g=g_A$ be as in the proof
of Theorem~\ref{thm:aggnb} at $w=w^\ast$.  Since
$\sum_j(w^\ast-y_j)\,\E X_j=0$,
\[
\sum_{j}(w^\ast-y_j)\,\E[X_jg(W+1)]
=\sum_{j}(w^\ast-y_j)\Bigl(\E[X_jg(W+1)]-\E X_j\,\E g(W+1)\Bigr),
\]
and for each $j$, with $W=W_j+V_j$,
\[
\begin{aligned}
\E[X_jg(W+1)]-\E X_j\,\E g(W+1)
&=\E\bigl[X_j\{g(W_j+V_j+1)-g(W_j+X_j+1)\}\bigr]\\
&\quad+\Bigl(\E[X_jg(W_j+X_j+1)]-\E X_j\,\E g(W_j+X_j+1)\Bigr)\\
&\quad+\E X_j\,\E\bigl[g(W_j+X_j+1)-g(W_j+V_j+1)\bigr].
\end{aligned}
\]
The first and third terms are at most
$\Lambda\sum_{i\in B_j\setminus\{j\}}\E[X_iX_j]$ and
$\Lambda\,\E X_j\sum_{i\in B_j\setminus\{j\}}\E X_i$.  In the second,
$X_j$ is independent of $W_j$, so
$\E[X_jg(W_j+X_j+1)]=\E X_j\,\E g(W_j+1+X_j^{s})$ with $X_j^{s}$ the
size-biased $X_j$, independent of $W_j$; by Lemma~\ref{lem:sizebias}
at $c=1$, $X_j^{s}\stackrel{d}{=}X_j+1+X_j'$ with $X_j'$ an
independent copy, so the second term is at most
$\Lambda\,\E X_j\,\E(1+X_j')=\Lambda\,y_j/(1-y_j)^{2}$.  Summing over
$j$ gives $\Lambda M^\ast$; the terms $\delta_j$ contribute
$\Lambda Q$ as in Theorem~\ref{thm:aggnb}, with $\beta=0$.
\end{proof}

\subsection*{Proofs for Section~\ref{sec:hnb}}

\begin{proof}[Proof of Theorem~\ref{thm:target}]
By independence, for $n\ge0$,
\[
\begin{aligned}
\PP\bigl(R_i=n+\kappa_i\ \forall i\bigr)
&=\prod_{i=1}^r\binom{n+\kappa_i+c-1}{n+\kappa_i}(1-q_i)^cq_i^{\,n+\kappa_i}\\
&=\Bigl(\prod_i(1-q_i)^cq_i^{\kappa_i}\Bigr)
\prod_i\binom{n+\kappa_i+c-1}{n+\kappa_i}\,w^n,
\end{aligned}
\]
since $\prod_iq_i^{\,n}=w^n$.  Since $D_i=k_i$ for all $i$ forces $R_i=R_r+\kappa_i$, the event
$\{D_i=k_i\ \forall i\}$ is the disjoint union over $n=R_r\ge0$ of
the events $\{R_i=n+\kappa_i\ \forall i\}$, all of them possible
because $\kappa_i\ge0$.
Summing over $n$ gives \eqref{eq:fiber}; dividing gives the conditional
law \eqref{eq:hnbpmf}.  Dependence on $\bq$ only through $w$ is
visible from \eqref{eq:hnbpmf}.
\end{proof}

\begin{proof}[Proof of Proposition~\ref{prop:kemp}]
Write $\binom{n+\kappa+c-1}{n+\kappa}=(c)_{n+\kappa}/(n+\kappa)!$, and use
$(c)_{n+\kappa}=(c)_{\kappa}(c+\kappa)_n$ and
$(n+\kappa)!=\kappa!\,(\kappa+1)_n$.  Since $\kappa_r=0$, the $i=r$
factor contributes $(c)_n/n!$, and
\[
\prod_{i=1}^r\binom{n+\kappa_i+c-1}{n+\kappa_i}
=\prod_{i=1}^{r-1}\binom{\kappa_i+c-1}{\kappa_i}\cdot
\frac{\prod_{i=1}^{r-1}(c+\kappa_i)_n\cdot(c)_n}
{\prod_{i=1}^{r-1}(\kappa_i+1)_n\cdot n!},
\]
which is the coefficient form of \eqref{eq:zhathyper}.  For
$r=2$ the displayed ${}_2F_1$ form follows with $k=\kappa_1$.  The
successive-ratio form
$p_{n+1}/p_n=\tfrac{(c+n)(c+k+n)}{(k+1+n)(1+n)}\,w$ is the canonical
Kemp-family ratio with numerator parameters $(c,c+k)$, denominator
parameter $k+1$, and argument $w\in(0,1)$; see \cite{KK56,Kemp68} and
\cite[\S2.4]{JKK05}.
\end{proof}

\begin{proof}[Proof of Lemma~\ref{lem:basic}]
(i) Immediate from \eqref{eq:hnbpmf}:
\[
\frac{p_{n+1}}{p_n}=w\prod_i\frac{(c)_{n+1+\kappa_i}/(c)_{n+\kappa_i}}
{(n+1+\kappa_i)!/(n+\kappa_i)!}
=w\prod_i\frac{n+\kappa_i+c}{n+1+\kappa_i}=\frac{w\,b(n)}{d(n+1)}.
\]

(ii) Each factor equals $1+\frac{c-1}{n+1+\kappa_i}$, strictly decreasing
in $n$ iff $c\ge2$ and constant ($=1$) iff $c=1$; for $c=1$,
$t_n\equiv w$ gives $p_n\propto w^n$.

(iii) Monotonicity of $t_n$ makes $p$ unimodal with mode
$m^\ast=\max\{n\ge1:t_{n-1}\ge1\}$, i.e.\ $w\,b(n-1)\ge d(n)$.  For the
location: let $s=w^{-1/r}-1>0$ and $\bar\kappa=\max_i\kappa_i$.  If
$m+1+\bar\kappa\le(c-1)/s$ then every factor of $t_m$ is
$\ge1+s$, so $t_m\ge(1+s)^rw=1$, whence $m^\ast\ge(c-1)/s-\bar\kappa-1$.
If $m+1\ge(c-1)/s$ then every factor is $\le1+s$, so $t_m\le1$, whence
$m^\ast\le(c-1)/s+1$.  Thus $|m^\ast-(c-1)/s|\le\bar\kappa+1$, which
is Remark~\ref{rem:mode}.

(iv) $\E N=\sum_n n\,p_n=w\Zhat'/\Zhat$ since $p_n$ is proportional to
the coefficient of $w^n$ in $\Zhat$ times $w^n$.  For \eqref{eq:momid},
sum the balance identity: $\E\,d(N)=\sum_{n\ge1}d(n)p_n
=\sum_{n\ge1}w\,b(n-1)p_{n-1}=w\,\E\,b(N)$, the reindexing justified by
absolute convergence (the terms are dominated by
$b(n)p_n\le\mathrm{const}\cdot n^{r c}w^n$).
\end{proof}

\begin{proof}[Proof of Proposition~\ref{prop:poissonlimit}]
Let $a_n(c):=p_n/p_0=w^n\prod_i(c+\kappa_i)_n\,\kappa_i!/(n+\kappa_i)!$,
and $a_n^{\HB}:=\nu^n\prod_i\kappa_i!/(n+\kappa_i)!$.  Since
$(c+\kappa_i)_n/c^n\to1$ for each fixed $n$ and $wc^r\to\nu$, we get
$a_n(c)\to a_n^{\HB}$ pointwise.  For a uniform tail bound, note
\[
\frac{a_{n+1}(c)}{a_n(c)}=w\prod_i\frac{c+\kappa_i+n}{n+1+\kappa_i}
\le w\,\frac{(c+\bar\kappa+n)^r}{(n+1)^r}.
\]
If $n+1\le c$ and $c\ge\bar\kappa$ then $c+\bar\kappa+n\le3c$, so the
ratio is at most $3^r wc^r/(n+1)^r\le 3^r(\nu+1)/(n+1)^r\le\frac12$ for
all $n\ge N_1:=\lceil(2\cdot3^r(\nu+1))^{1/r}\rceil$ and $c$ large.  If
$n+1>c\ (\ge\bar\kappa+1)$ then $c+\bar\kappa+n\le3(n+1)$, so the ratio
is at most $3^rw\le\frac12$ for $c$ large, since $w\to0$.  Hence for all
large $c$, $a_n(c)\le a_{N_1}(c)\,2^{-(n-N_1)}$ for $n\ge N_1$, and
$a_{N_1}(c)$ is bounded; dominated convergence for series gives
$\sum_na_n(c)\to\sum_na_n^{\HB}$, hence pointwise convergence of the
normalized mass functions, hence $\dtv\to0$ by Scheff\'e's lemma.  For
the operator, $w\,b(n)=wc^r\prod_i(1+\tfrac{n+\kappa_i}{c})\to\nu$
pointwise.
\end{proof}

\begin{proof}[Proof of Theorem~\ref{thm:fixedpoint}]
Let $W\sim\HNB(\bk,c,w)$.  For $n\ge1$ the balance identity gives
\[
\PP(W^d=n)=\frac{d(n)p_n}{\beta_d}=\frac{w\,b(n{-}1)p_{n-1}}{\beta_d},
\]
and $\beta_d=w\beta_b$ by Lemma~\ref{lem:basic}(iv), so the right side
is $\PP(W^b{+}1=n)$; both laws live on $\{1,2,\dots\}$.

Conversely, suppose $W^d\stackrel{d}{=}W^b+1$, i.e.\
$d(n)\PP(W{=}n)/\beta_d=b(n{-}1)\PP(W{=}n{-}1)/\beta_b$ for all
$n\ge1$.  If $\PP(W{=}0)=0$ then, since $d(n)>0$ for $n\ge1$, induction
gives $\PP(W{=}n)=0$ for all $n$, a contradiction; so $\PP(W{=}0)>0$
and the recursion
$\PP(W{=}n)/\PP(W{=}n{-}1)=(\beta_d/\beta_b)\,b(n{-}1)/d(n)$ determines
the law as \eqref{eq:hnbpmf} with $w=\beta_d/\beta_b$.  Summability of
$\sum_n\prod_i\binom{n+\kappa_i+c-1}{n+\kappa_i}w^n$ forces $w<1$.  For
the transform: $TW=W$ says
$\PP(W{=}n)\propto d(n{+}1)\PP(W{=}n{+}1)/b(n)$, which is the same
recursion; the normalizing constant of $T$ plays the role of $1/w$, and
summability again forces $w\in(0,1)$.
\end{proof}

\begin{proof}[Proof of Theorem~\ref{thm:master}]
Theorem~\ref{thm:nbmaster} at the rates $b(n)=\prod_i(n+\kappa_i+c)$,
$d(n)=\prod_i(n+\kappa_i)$, which satisfy the standing assumptions
of \S\ref{sec:nb}: $d(0)=0$, $d\le b$, $b$ nondecreasing,
$d(n)/b(n)\to1$ and $b(n{+}1)/b(n)\to1$, with the balance identity
of Lemma~\ref{lem:basic}(i); $\HNB(\bk,c,w)=\pi_{b,d,w}$ by
Definition~\ref{def:hnb}, and \eqref{eq:master},
\eqref{eq:mastermatched} are the two displays of that theorem.
\end{proof}

\begin{proof}[Proof of Lemma~\ref{lem:phantom}]
(i) $(n+c)\binom{n+c-1}{n}=\tfrac{(n+c)!}{n!(c-1)!}=c\binom{n+c}{n}$ and
$\E[R+c]=\tfrac{cq}{1-q}+c=\tfrac{c}{1-q}$, so
\[
\frac{(n+c)\PP(R{=}n)}{\E[R+c]}
=\binom{n+c}{n}(1-q)^{c+1}q^n=\PP\bigl(\NB(c{+}1,q)=n\bigr).
\]
(ii) Probability generating functions:
$\bigl(\tfrac{1-q}{1-qz}\bigr)^c\cdot\tfrac{1-q}{1-qz}
=\bigl(\tfrac{1-q}{1-qz}\bigr)^{c+1}$.
(iii) is the combination of (i) and (ii).
\end{proof}

\begin{proof}[Proof of Lemma~\ref{lem:sizebias}]
$n\binom{n+c-1}{n}=\tfrac{(n+c-1)!}{(n-1)!(c-1)!}=c\binom{n+c-1}{n-1}$
and $\E R=cq/(1-q)$, so
$\PP(R^s{=}n)=\binom{(n-1)+c}{\,n-1\,}(1-q)^{c+1}q^{\,n-1}
=\PP\bigl(1+\NB(c{+}1,q)=n\bigr)$; the second equality in the statement
is Lemma~\ref{lem:phantom}(ii).
\end{proof}

\begin{proof}[Proof of Lemma~\ref{lem:geomchar}]
The $x_i$-biased law of $Z$ has mass
$z_i\PP(Z{=}z)/\E X^{(i)}$; the $(x_i{+}1)$-biased law shifted by $e_i$
has mass at $z$ equal to $z_i\PP(Z{=}z{-}e_i)/\E(X^{(i)}{+}1)$
(substituting $z_i-1$ for the biased coordinate).  Equality of the two
laws gives, for every $z$ with $z_i\ge1$,
\[
\frac{\PP(Z{=}z)}{\PP(Z{=}z{-}e_i)}\;=\;\frac{\E X^{(i)}}{\E(X^{(i)}{+}1)}
\;=:\;q_i\ \in(0,1),
\]
whenever the denominator is positive; and if $\PP(Z{=}z{-}e_i)=0$ the
same identity forces $\PP(Z{=}z)=0$, so the support is an order ideal in
coordinate $i$ and the ratio identity propagates.  Iterating,
$\PP(Z{=}z)=q_i^{\,z_i}\,\PP\bigl(Z=(z_{-i},0)\bigr)$ for all $z$, which
says exactly that $X^{(i)}\sim\Geom(1-q_i)$ independent of
$(X^{(j)})_{j\ne i}$.

For the single-coordinate statement, let $f$ be the probability
generating function of $X$.  The $(x{+}1)$-bias has generating function
$(f(z)+zf'(z))/\E(X{+}1)$, and $X+X'$ has $f(z)^2$; equality on $[0,1)$
gives the Riccati equation $f+zf'=c_0f^2$ with $c_0=\E(X{+}1)$.
Here $p_0>0$, since the least point of the support of $X+X'$ is twice
that of $X$ while that of the $(x{+}1)$-bias is the least point of
$X$, so $f>0$ on $[0,1]$.  Substituting $u=1/f$ yields
$u'=(u-c_0)/z$, whence $u=c_0+Kz$; the
boundary values $f(0)=p_0$, $f(1)=1$ give $c_0=1/p_0$ and $K=1-1/p_0\le0$,
so $f(z)=p_0/\bigl(1-(1-p_0)z\bigr)$: geometric with parameter
$q=1-p_0$ (including the degenerate $p_0=1$).  The converse is the
computation $(n{+}1)(1-q)q^n\cdot(1-q)=\binom{n+1}{n}(1-q)^2q^n$, i.e.\
$\NB(2,q)$.
\end{proof}

\begin{proof}[Proof of Proposition~\ref{prop:factorbound}]
Theorem~\ref{thm:factors}, whose formulas use only the birth--death
form; the negative-binomial specialization substitutes $b(n)=n+c$,
$d(n)=n$.
\end{proof}

\begin{proof}[Proof of Theorem~\ref{thm:colexpansion}]
Corollary~\ref{cor:productform}: the shared-phantom coupling makes
every coupling cost vanish for geometric coordinates, and the general
expansion, Theorem~\ref{thm:expansion}, reduces to the display; the
window and the tilted columns are Definition~\ref{def:costs}.
\end{proof}

\begin{proof}[Proof of Lemma~\ref{lem:convmom}]
$\Zhat_c(\bk,w)$ is the normalizer of an exponential family in
$\log w$; its first two $\log w$-derivatives are the mean and
variance, and independence adds them.
\end{proof}

\begin{proof}[Proof of Theorem~\ref{thm:aggconv}]
Replacing $X_i$ by $N_i$ one coordinate at a time gives
$\dtv(V,V^\ast)\le\sum_i\dtv(X_i,N_i)$, since convolution with a
fixed independent variable contracts total variation; each summand
is Theorem~\ref{thm:master}.
\end{proof}

\begin{proof}[Proof of Lemma~\ref{lem:transfer}]
The two algebraic expansions are immediate:
$\prod_iU_i=\prod_i\sum_jX^{(i)}_j=\sum_{\vec\jmath}\prod_iX^{(i)}_{\jmath_i}$,
and since $U_i+c=\sum_{j=1}^{c}(X^{(i)}_j+1)$,
$\prod_i(U_i+c)=\sum_{\vec\jmath}\prod_i(X^{(i)}_{\jmath_i}+1)$ ---
crucially over the \emph{same} index set $[c]^r$.  Fix $\vec\jmath$ and
consider the $d$-side summand
$\E\bigl[\prod_iX^{(i)}_{\jmath_i}\,g(U_r)\ind{D=\bkk}\bigr]$.
The weight factors over the hit columns as
$\prod_{j\in L(\vec\jmath)}\prod_{i\in S_j}X^{(i)}_j$, a product of
per-column functions; since the columns are independent, reweighting by
it replaces each hit column's law by its $S_j$-tilt and multiplies by the
normalizer $\chi^{S_j,x}_j$, while $g(U_r)\ind{D=\bkk}$ is a function of
the (reweighted) column values.  This gives \eqref{eq:transferd}, and
identically \eqref{eq:transferb} with the $(x{+}1)$-weights.  Two domain
remarks: on the left of \eqref{eq:transferd} the factor $U_r$ vanishes
when $U_r=0$, and on the right the block containing $r$, namely
$S_{\jmath_r}\ni r$, carries the weight $x_r$, so the tilted system has
$X^{(r)}_{\jmath_r}\ge1$ almost surely and $V^{\vec\jmath,x}\ge1$: the
convention $g(0)=0$ is never invoked with positive mass.
\end{proof}

\begin{proof}[Proof of Lemma~\ref{lem:vanishing}]
($\Leftarrow$)  Suppose each column has independent coordinates
$X^{(i)}_j\sim\Geom(1-q_i)$.  Fix $j$ and $S$.  The tilts act
coordinatewise on independent coordinates, so it suffices to couple each
coordinate $i\in S$.  By Lemma~\ref{lem:sizebias} (with $c=1$), the
$x_i$-bias of $\Geom(1-q_i)$ is $1+X_i+X_i'$ with $X_i'$ a fresh
independent $\Geom(1-q_i)$; by Lemma~\ref{lem:geomchar} (converse
direction) the $(x_i{+}1)$-bias is $X_i+X_i'$.  Realizing both with the
\emph{same} pair $(X_i,X_i')$ --- the shared phantom --- makes the
$d$-tilt equal to the $b$-tilt plus $e_S$ pointwise, so
$\eta^{(S)}_j=0$.  For the weights,
$\chi^{S,x}_j=\prod_{i\in S}\tfrac{q_i}{1-q_i}
=w_S\prod_{i\in S}\tfrac1{1-q_i}=w_S\,\chi^{S,x+1}_j$, so
$\delta^{(S)}_j=0$.  That $W=\HNB(\bk,c,w)$ exactly is
Theorem~\ref{thm:target}.

($\Rightarrow$)  Suppose all singleton costs vanish for column $j$.
$\eta^{(\{i\})}_j=0$ means the $x_i$-tilt of $Z_j$ equals the
$(x_i{+}1)$-tilt shifted by $e_i$ in law; Lemma~\ref{lem:geomchar} then
gives $X^{(i)}_j\sim\Geom(1-\tilde q_i)$ with
$\tilde q_i=\E X^{(i)}_j/\E(X^{(i)}_j{+}1)$, independent of the other
coordinates.  Then $\delta^{(\{i\})}_j=
|\E X^{(i)}_j-q_i\,\E(X^{(i)}_j{+}1)|
=\E(X^{(i)}_j{+}1)\,|\tilde q_i-q_i|=0$ pins $\tilde q_i=q_i$.
Applying this to every $i$ yields the claim (and the multi-coordinate
costs then vanish by the first direction).
\end{proof}

\begin{proof}[Proof of Theorem~\ref{thm:expansion}]
Write $P_k=\PP(D=\bkk)$ and let $g=g_A$ be the Stein solution for
$\HNB(\bk,c,w)$, so $\|g\|_\infty\le\Gamma$ and consecutive increments
are bounded by $\Lambda$.  On the fiber $\{D=\bkk\}$ one has
$U_i=U_r+\kappa_i$, hence $d(U_r)=\prod_iU_i$ and
$b(U_r)=\prod_i(U_i+c)$ there, and
\[
P_k\bigl[\PP(W\in A)-\HNB(A)\bigr]
= w\,\E\Bigl[\prod_i(U_i{+}c)\,g(U_r{+}1)\ind{D=\bkk}\Bigr]
-\E\Bigl[\prod_iU_i\,g(U_r)\ind{D=\bkk}\Bigr].
\]
Apply Lemma~\ref{lem:transfer} to both terms and pair the summands over
the common index set:
\[
P_k\bigl[\PP(W\in A)-\HNB(A)\bigr]
=\sum_{\vec\jmath\in[c]^r}\Bigl\{\,w\,W_b(\vec\jmath)\,E_b(\vec\jmath)
-W_d(\vec\jmath)\,E_d(\vec\jmath)\Bigr\},
\]
where $W_b=\prod_j\chi^{S_j,x+1}_j$, $W_d=\prod_j\chi^{S_j,x}_j$,
$E_b=\E[g(V^{x+1}{+}1)\ind{D^{x+1}=\bkk}]$ and
$E_d=\E[g(V^{x})\ind{D^{x}=\bkk}]$.  Split each bracket as
\[
w\,W_b\,(E_b-E_d)\;+\;(w\,W_b-W_d)\,E_d .
\]

\emph{Step 1: the coupling term.}  Fix $\vec\jmath$.  Realize the two
tilted systems on one probability space: the untouched columns are shared
verbatim; for each hit column $j\in L(\vec\jmath)$, choose a coupling of
$\bigl(Z^{S_j,x}_j,\,Z^{S_j,x+1}_j\bigr)$ attaining
$\E\|Z^{S_j,x}_j-Z^{S_j,x+1}_j-e_{S_j}\|_1=\eta^{(S_j)}_j$ (optimal
$\ell_1$ couplings exist on countable spaces), independently across hit
columns.  Write $Z^{S_j,x}_j=Z^{S_j,x+1}_j+e_{S_j}+R_j$ with
$\E\|R_j\|_1=\eta^{(S_j)}_j$.  Because the blocks $S_j$ partition $[r]$,
$\sum_{j\in L}e_{S_j}=\one$, so coordinatewise
$U^{x}_i=U^{x+1}_i+1+\rho_i$ with $\rho=\sum_jR_j$; in particular
$V^{x}=V^{x+1}+1+\rho_r$ and $D^{x}=D^{x+1}+\partial\rho$, where
$\partial$ is the difference map, so
$\E|\rho_r|\le\sum_j\eta^{(S_j)}_j$ and
$\E\|\partial\rho\|_1\le2\sum_j\eta^{(S_j)}_j$ (each coordinate of
$\rho$ enters at most two differences).  Now telescope:
\[
E_d-E_b=\E\Bigl[g\bigl(V^{x+1}{+}1{+}\rho_r\bigr)\ind{D^{x+1}+\partial\rho=\bkk}
-g\bigl(V^{x+1}{+}1\bigr)\ind{D^{x+1}=\bkk}\Bigr].
\]
First remove $\rho_r$ one unit at a time: each unit step changes the
integrand by at most $\Lambda$ on the event
$\{D^{x+1}{+}\partial\rho=\bkk\}$, and, conditioning on all hit-column
values, that event reads $\{D_R^{(L)}=\bkk-\partial z^{x}\}$ with
$z^{x}$ the sum of the hit columns under the $x$-tilts, a point of
$\mathcal{W}$, so its conditional probability is at most
$\theta P_k$; hence this stage contributes at most
$\Lambda\theta P_k\sum_j\eta^{(S_j)}_j$.  Then remove $\partial\rho$ one
coordinate-unit at a time, along a coordinatewise monotone lattice
path from $\bkk-\partial z^{x}$ to $\bkk-\partial z^{x+1}$, $z^{x+1}$
the sum of the hit columns under the $(x{+}1)$-tilts, so that every
point of the path lies in $\mathcal{W}$: for adjacent
$s,s'\in\mathcal{W}$, conditioning on the hit columns,
\[
\begin{aligned}
\bigl|\E[g\,\ind{D_R=s}]-\E[g\,\ind{D_R=s'}]\bigr|
&\le\Gamma\,\bigl|\PP(D_R{=}s)-\PP(D_R{=}s')\bigr|\\
&\qquad+\PP(D_R{=}s')\ \bigl|\E[g\mid s]-\E[g\mid s']\bigr|
\ \le\ (\Gamma\tau_1+\Lambda\theta\omega)\,P_k,
\end{aligned}
\]
using $|\E[g\mid s]-\E[g\mid s']|\le\Lambda\,
\dwass(\mathcal{L}(V_R\mid s),\mathcal{L}(V_R\mid s'))\le\Lambda\omega$
and the Lipschitz property of $g$; this stage contributes at most
$2(\Gamma\tau_1+\Lambda\theta\omega)P_k\sum_j\eta^{(S_j)}_j$.  Altogether
\[
|E_b-E_d|\;\le\;\bigl[\theta(1{+}2\omega)\Lambda+2\tau_1\Gamma\bigr]\,
P_k\sum_{j\in L(\vec\jmath)}\eta^{(S_j)}_j .
\]

\emph{Step 2: the weight term.}  $|E_d|\le\Gamma\,\PP(D^{x}=\bkk)
\le\Gamma\theta P_k$ by the same conditioning.  Telescoping the products
with per-column comparands $\chi^{S_j,x}_j$ versus
$w_{S_j}\chi^{S_j,x+1}_j$ (note $\prod_jw_{S_j}=w$ because the blocks
partition $[r]$) and bounding each partial product by the corresponding
$\chi^{x+1}$'s (valid since $\chi^{S,x}\le\chi^{S,x+1}$ and
$w_S\le1$),
\[
|w\,W_b-W_d|\;\le\;\sum_{j\in L(\vec\jmath)}\delta^{(S_j)}_j
\prod_{j'\in L(\vec\jmath)\setminus\{j\}}\chi^{S_{j'},x+1}_{j'} .
\]

\emph{Step 3: aggregation.}  Sum over $\vec\jmath$, grouping by the pair
$(j,S)$ with $S=S_j$: for the coupling term,
\[
\sum_{\vec\jmath}w\,W_b\sum_{j\in L}\eta^{(S_j)}_j
= w\sum_{\emptyset\ne S\subseteq[r]}\sum_{j=1}^c
\chi^{S,x+1}_j\eta^{(S)}_j
\sum_{\vec\jmath'\,:\,[r]\setminus S\to[c]\setminus\{j\}}
\prod \chi^{\,\cdot\,,x+1}
\;\le\;w\sum_S\nu^b_{S^c}\,H^{(S)},
\]
because the inner sum is, by the transfer identity run in reverse,
$\E\prod_{i\in S^c}\bigl(U_i-X^{(i)}_j+c-1\bigr)\le
\E\prod_{i\in S^c}(U_i+c)=\nu^b_{S^c}$ (coordinatewise domination of
nonnegative factors).  The weight term aggregates identically to
$\Gamma\theta\sum_S\nu^b_{S^c}\Delta^{(S)}$.  Taking the supremum over
$A$ and dividing by $P_k$ gives \eqref{eq:expansion}.
\end{proof}

\begin{proof}[Proof of Corollary~\ref{cor:productform}]
By Lemma~\ref{lem:vanishing} ($\Leftarrow$ direction), applied
coordinatewise with parameter $y_{ij}$ in place of $q_i$: the shared
phantom coupling makes the $S$-tilts differ by exactly $e_S$
\emph{pointwise for any parameters}, since the identities
$x$-bias $=1+X+X'$ and $(x{+}1)$-bias $=X+X'$ hold for every geometric
law; the reference $q_i$ enters only through the weights.  Hence
$\eta^{(S)}_j\equiv0$ and only the $\delta$'s survive; the moment
computations $\chi^{S,x}_j=\prod_{i\in S}y_{ij}/(1-y_{ij})$,
$\chi^{S,x+1}_j=\prod_{i\in S}(1-y_{ij})^{-1}$ give
\eqref{eq:productform}, and $y_{ij}=q_i\beta_j$ gives
$|\prod_{i\in S}y_{ij}-w_S|=w_S|\beta_j^{|S|}-1|$.
\end{proof}

\begin{proof}[Proof of Corollary~\ref{cor:capped}]
\emph{Weight cost, singleton, exact.}  Let
$\widetilde X\sim\Geom(1-y)$ conditioned on $\widetilde X\le L$, and
take reference $q=y$.  With $F_L=1-y^{L+1}$,
$M=\sum_{t\ge1}ty^t=y/(1-y)^2$ and $T_1(L)=\sum_{t>L}ty^t
=y^{L+1}\bigl((L{+}1)(1-y)+y\bigr)/(1-y)^2$,
\[
\delta=\bigl|(1-y)\E\widetilde X-y\bigr|
=\Bigl|\frac{(1-y)^2\,(M-T_1(L))}{F_L}-y\Bigr|
=\frac{(L{+}1)(1-y)\,y^{L+1}}{1-y^{L+1}},
\]
after simplification (the numerator collapses to
$-(L{+}1)(1-y)y^{L+1}$), which is \eqref{eq:deltacap}.

\emph{Coupling cost, singleton, explicit bound.}  Both tilts touch only
coordinate $i$ and the other coordinates are shared, so
$\eta^{(\{i\})}$ equals the $\dwass$ distance on $\Zint$ between the
size-bias $\widetilde s$ of $\widetilde X$ and the shifted
$(x{+}1)$-bias $\widetilde r$ of $\widetilde X$.  Both are perturbations
of the common uncapped law $s=\mathcal{L}(1+\NB(2,y))$ with mass
$s(t)=ty^{t-1}(1-y)^2$, $t\ge1$; the triangle inequality and the
elementary bound $\dwass(\mu,\nu)\le\sum_tt\,|\mu(t)-\nu(t)|$ give
$\eta^{(\{i\})}\le\mathcal{E}(y,L)$ with
\begin{equation}\label{eq:etacap}
\mathcal{E}(y,L)\;:=\;\frac{M_2\,T_1(L)}{\widetilde M\,M}
+\frac{(1-y)^2\,T_2(L)}{y}
+\frac{N_2\,T_1(L{+}1)}{y\,\widetilde N\,N}
+\frac{(1-y)^2\,T_2(L{+}1)}{y},
\end{equation}
where $M_2=\sum t^2y^t=y(1+y)/(1-y)^3$,
$\widetilde M=M-T_1(L)$, $N=\sum_{t\ge1}ty^{t-1}=1/(1-y)^2$,
$N_2=(1+y)/(1-y)^3$, $\widetilde N=N-T_1(L{+}1)/y$, and
$T_2(K)=\sum_{t>K}t^2y^t$.  All quantities are elementary closed forms;
using $T_1(K)\le(K{+}2)y^{K+1}/(1-y)^2$ and
$T_2(K)\le(2K^2{+}4)\,y^{K+1}/(1-y)^3\le3(K{+}2)^2y^{K+1}/(1-y)^3$
(both by splitting $t=K{+}s$ and $t^2\le2K^2{+}2s^2$),
$\mathcal{E}(y,L)=O\bigl((L{+}2)^2y^{L+1}\bigr)$; the constants are
explicit but not optimized.

\emph{Multi-coordinate costs.}  Coordinates within a column are
independent and the tilts factorize, so a product coupling gives
$\eta^{(S)}_j\le\sum_{i\in S}\mathcal{E}(y_{ij},L_{ij})$ ($\ell_1$ costs
add), and telescoping
$\bigl|\prod_{i\in S}a_i-\prod_{i\in S}b_i\bigr|
\le\sum_{i\in S}|a_i-b_i|\prod_{i'\in S\setminus\{i\}}\max(a_{i'},b_{i'})$
with $a_i=\E\widetilde X^{(i)}_j$ and
$b_i=q_i\,\E(\widetilde X^{(i)}_j{+}1)$ gives
$\delta^{(S)}_j\le\sum_{i\in S}\delta^{(\{i\})}_j
\prod_{i'\in S\setminus\{i\}}\E(\widetilde X^{(i')}_j{+}1)$.
\end{proof}

\subsection*{Proofs for Section~\ref{sec:rowcol}}

\begin{proof}[Proof of Lemma~\ref{lem:sysid}]
$\pi_{w,\mu}(k{+}1)/\pi_{w,\mu}(k)
=w\,\rho_1(k)^{1-\mu}\rho_2(k)^{\mu}
=\bigl[w\rho_1(k)\bigr]\cdot e^{\mu(S(k+1)-S(k))}$; telescoping from
$k=0$ gives the claim, and summability follows from
$w\rho_i\le\bar\rho<1$ eventually, uniformly in $\mu\in[0,1]$.
\end{proof}

\begin{proof}[Proof of Lemma~\ref{lem:syslc}]
The standing assumptions pass to the geometric means:
$\tilde d_\mu(0)=0$, $\tilde d_\mu\le\tilde b_\mu$ since $d_i\le b_i$,
$\tilde b_\mu$ is nondecreasing, and
$\tilde d_\mu/\tilde b_\mu=(d_1/b_1)^{1-\mu}(d_2/b_2)^{\mu}\to1$,
$\tilde b_\mu(n{+}1)/\tilde b_\mu(n)\to1$.  A geometric mean of
nonincreasing positive functions is nonincreasing, and a law with
nonincreasing pmf ratio is log-concave; \eqref{eq:gammaexact} and
\eqref{eq:lambdabound} are rate-generic, and \eqref{eq:gammabound}
needs only the strict decrease of the ratio.
\end{proof}

\begin{proof}[Proof of Theorem~\ref{thm:sysmaster}]
The target $\pi_{w,\mu}$ is summable since its ratio
$w\tilde b_\mu(k)/\tilde d_\mu(k{+}1)=w\rho_1(k)^{1-\mu}\rho_2(k)^{\mu}$
is at most $\bar\rho<1$ for $k\ge k_0$.  For $A\subseteq\N$ the
solution $g_A$ is given by \eqref{eq:solution}, whose derivation in
the proof of Theorem~\ref{thm:solution}(ii) uses only $\tilde d_\mu(0)=0$,
$\tilde b_\mu>0$ and the balance identity, and the proof of
Theorem~\ref{thm:factors}(i) gives
$\sup_A|g_A(m{+}1)|=F(m)G(m{+}1)/(\tilde d_\mu(m{+}1)p_{m+1})$.  This
is finite uniformly in $m$: for $m\ge k_0$,
$G(m{+}1)\le p_{m+1}/(1-\bar\rho)$, so the ratio is at most
$1/\bigl((1-\bar\rho)\inf_{n\ge1}\tilde d_\mu(n)\bigr)$, and the
finitely many $m<k_0$ are finite.  Hence $\|g_A\|_\infty\le\tilde\Gamma$
and $\|\Delta g_A\|_\infty\le\tilde\Lambda\le2\tilde\Gamma$.  By
Definition~\ref{def:bias} with the tilded rates,
$\E[\tilde b(W)g_A(W{+}1)]=\beta_{\tilde b}\,\E g_A(W^{\tilde b}{+}1)$ and
$\E[\tilde d(W)g_A(W)]=\beta_{\tilde d}\,\E g_A(W^{\tilde d})$, so
\[
\PP(W\in A)-\pi_{w,\mu}(A)
=w\beta_{\tilde b}\bigl[\E g_A(W^{\tilde b}{+}1)-\E g_A(W^{\tilde d})\bigr]
+\bigl(w\beta_{\tilde b}-\beta_{\tilde d}\bigr)\E g_A(W^{\tilde d}),
\]
and the two bounds give the display; take the supremum over $A$.
\end{proof}

\begin{proof}[Proof of Theorem~\ref{thm:sysmatch}]
(i) is Theorem~\ref{thm:sysmaster} with the stated $w$.  (ii) At
$w=w^\ast(\mu)$ the defect for $g\equiv1$ is
$w^\ast\beta_{\tilde b}-\beta_{\tilde d}=0$ by definition of
$w^\ast$, and the defect for $g=S$ is
$w^\ast\beta_{\tilde b}\,\E S(W^{\tilde b}{+}1)-\beta_{\tilde d}\,
\E S(W^{\tilde d})=\beta_{\tilde d}\bigl(\E S(W^{\tilde b}{+}1)-
\E S(W^{\tilde d})\bigr)=-\beta_{\tilde d}\,\zeta(\mu)$, which
vanishes at a zero of $\zeta$.  The map $\mu\mapsto\zeta(\mu)$ is
continuous: $\tilde d_\mu\le d_1+d_2$ and $\tilde b_\mu\le b_1+b_2$
pointwise, so the integrability assumed makes every expectation in
$\zeta$ continuous in $\mu$ by dominated convergence; the
intermediate value theorem gives $\mu^\ast$.  (iii) On the family of Lemma~\ref{lem:sysid} the
log-normalizer $\log Z(\log w,\mu)$ is strictly convex on
$(-\infty,0)\times[0,1]$ because the statistics $(k,S(k))$ are
affinely independent, so its gradient $(\E_\pi k,\E_\pi S)$ is
injective and the pair is unique; the projection property is the
exponential-family Pythagorean identity.
\end{proof}

\begin{proof}[Proof of Proposition~\ref{prop:shadow}]
The product-form weight is constant on each margin fiber:
$\PP(U=u,V=v)=|\T(u,v)|\,\prod_ix_i^{u_i}\prod_jy_j^{v_j}
\prod_{ij}(1-x_iy_j)$.  Substituting the affine margins, the
$k$-dependent factor is $(\prod_ix_i)^k(\prod_jy_j)^{w'(k)}$, which
gains exactly $\Omega$ per step of the progression.
\end{proof}

\begin{proof}[Proof of Proposition~\ref{prop:atom}]
$rk^\ast+\sum_i\kappa_i=\sum_ir^\ast_i=N$ and
$\sum_j\lambda_j=N-c\,c^\ast_c$, so
$w'(k^\ast)=(N-(N-c\,c^\ast_c))/c=c^\ast_c$, an integer, and
$\bc(k^\ast)=\bc^\ast$; membership in $L$ is exactly this
integrality.  The display is Proposition~\ref{prop:shadow} evaluated
at the atom.
\end{proof}

\begin{proof}[Proof of Lemma~\ref{lem:crosscov}]
$\E_{P_\bullet}g(W)=\E\bigl[g(W)\,\ind{D_{\mathrm c}}\bigm|
D_{\mathrm r}\bigr]\big/\PP(D_{\mathrm c}\mid D_{\mathrm r})$, and
$\E[g(W)\mid D_{\mathrm r}]=0$ turns the numerator into the displayed
covariance.
\end{proof}

\begin{proof}[Proof of Corollary~\ref{cor:crossdef}]
Lemma~\ref{lem:crosscov} at each endpoint with $f\equiv1$.
\end{proof}

\begin{proof}[Proof of Theorem~\ref{thm:tablemaster}]
Theorem~\ref{thm:sysmaster} applied to $P_\bullet$, with the second
term removed by Theorem~\ref{thm:sysmatch}(i);
Corollary~\ref{cor:crossdef} identifies the endpoint inputs.
\end{proof}

\begin{proof}[Proof of Theorem~\ref{thm:system}]
Assembled from the Stein system section:
Theorems~\ref{thm:sysmaster} and~\ref{thm:sysmatch}(i),
Corollary~\ref{cor:crossdef}, and Propositions~\ref{prop:shadow}
and~\ref{prop:atom}; restated there in full as
Theorem~\ref{thm:tablemaster}.
\end{proof}

\section*{Declaration on the use of artificial intelligence}
This research was carried out in extended interactive sessions with
the large language model Claude (Anthropic), used as a research and
writing assistant.  The AI system contributed, under the author's
direction: initial derivations and first drafts of the results and
proofs; the numerical verification programs whose outputs are reported
in \S\ref{sec:computations}; and the drafting of the manuscript.  The
author formulated the research program and the questions, directed
each stage, independently checked the mathematical statements and
proofs, and takes sole and full responsibility for the correctness,
integrity, and originality of the work.  In accordance with publisher
policy, the AI system is not listed as an author and is not cited as a
source.

\end{document}